\documentclass[11pt,a4]{article}
\usepackage[T2A]{fontenc}
\usepackage[cp866]{inputenc}
\usepackage[english]{babel}
\usepackage[mathscr]{eucal}
\usepackage{amssymb}
\usepackage{amsmath}
\usepackage{fancyhea}
\renewcommand{\leq}{\leqslant}

\renewcommand{\C}{{\mathbb C}}

\newcommand{\R}{{\mathbb R}}

\newcommand{\K}{{\mathbb K}}

\renewcommand{\k}{\rule{0.7em}{0.7em}}
\begin{document}
\sloppy
\fussy
\headrulewidth = 2pt
\pagestyle{fancy}
\lhead[\scriptsize V.N.Gorbuzov]{\scriptsize V.N.Gorbuzov}
\rhead[\it 
\scriptsize 
Cylindricality
 and autonomy of integrals and last multipliers ...
]{\it \scriptsize 
Cylindricality and autonomy of integrals and last multipliers ...}
\headrulewidth=0.25pt

{\normalsize 

\thispagestyle{empty}

\mbox{}
\\[-5.5ex]
\mbox{}\hfill
{\large\bf 
CYLINDRICALITY AND AUTONOMY OF INTEGRALS 
}
\hfill\mbox{}
\\
\mbox{}\hfill
{\large\bf 
AND LAST 
MULTIPLIERS OF MULTIDIMENSIONAL
}
\hfill\mbox{}
\\
\mbox{}\hfill
{\large\bf 
DIFFERENTIAL SYSTEMS\footnote[1]{
The definitive version of this article has been published in the
monograph 
{\it Integrals of  Differential Systems}, 2006, Grodno State University, Grodno [1]; 
{\it Differential Equations}, Vol. 30 (1994), No. 6, 868-875 [2];
 {\it Differential Equations}, Vol. 34 (1998), No. 2, 149-156 [3];
{\it Vestnik of the Yanka Kupala Grodno State Univ.}, 
1999, Ser. 2, No. 1, 21-25 [4].
}
}
\hfill\mbox{}
\\[2.5ex]
\centerline{
\bf 
V.N. Gorbuzov
}
\\[2.5ex]
\centerline{
\it 
Department of Mathematics and Computer Science, 
Yanka Kupala Grodno State University,
}
\\[1ex]
\centerline{
\it 
Ozeshko 22, Grodno, 230023, Belarus
}
\\[1ex]
\centerline{
E-mail: gorbuzov@grsu.by
}
\\[7ex]
\centerline{
{\large\bf Abstract}}
\\[0.75ex]
\indent

The conditions of
cylindricality 
and autonomy  of first integrals, last multipliers and 
integral manifolds 
for  linear homogeneous systems of partial differential equations
and   total differential systems 
are established.
\\[2ex]
\\[0.75ex]
\\[5.5ex]
{\large\bf Contents} 
\\[1ex]
{\bf  Introduction} \dotfill\ 2
\\[0.5ex]
\mbox{}\hspace{1em}
0.1. Generalities  \dotfill\ 2
\\[0.5ex]
\mbox{}\hspace{1em}
0.2.  Problem definition  \dotfill\ 6
\\[1ex]
{\bf  
1. Cylindricality and autonomy of first integrals}          \dotfill\ 7
\\[0.5ex]
\mbox{}\hspace{1em}
1.1. 
Cylindricality of first integrals for
linear homogeneous system of partial 
\\
\mbox{}\hspace{3em}
differential equations
                                                                                              \dotfill\ 7
\\[0.5ex]
\mbox{}\hspace{3em}
1.1.1. Necessary condition of existence of  cylindrical first integral
                                                                                              \dotfill\ 7
\\[0.5ex]
\mbox{}\hspace{3em}
1.1.2.
Criterion  of existence of cylindrical first integral
                                                                                              \dotfill\ 8
\\[0.5ex]
\mbox{}\hspace{3em}
1.1.3.    Functionally independent cylindrical first integrals
                                                                                               \dotfill\ 10
\\[1ex]
\mbox{}\hspace{1em}
1.2.  First integrals of $s\!$-nonautonomous completely solvable 
total differential systems
                                                                                                .\, 10
\\[0.5ex]
\mbox{}\hspace{1em}
1.3.  
Autonomy  and  cylindricality of first integrals for
total differential system
                                                                                               \dotfill\ 12
\\[0.5ex]
\mbox{}\hspace{3em}
1.3.1. Necessary condition of existence of  
$s\!$-nonautonomous  $(n-k)\!$-cylindrical
\\
\mbox{}\hspace{5.7em}
first integral
                                                                                                \dotfill\ 12
\\[0.5ex]
\mbox{}\hspace{3em}
1.3.2.
Criterion  of existence of 
$s\!$-nonautonomous  $(n-k)\!$-cylindrical
\\
\mbox{}\hspace{5.7em}
first integral
					   \dotfill\ 13
\\[0.5ex]
\mbox{}\hspace{3em}
1.3.3.    Functionally independent 
$s\!$-nonautonomous  $(n-k)\!$-cylindrical
\\
\mbox{}\hspace{5.7em}
first integrals
					\dotfill\ 15
\\[1ex]
{\bf  
2. Cylindricality and autonomy of last multipliers}          \dotfill\ 15
\\[0.5ex]
\mbox{}\hspace{1em}
2.1. 
Cylindricality of last multipliers for
linear homogeneous system of partial 
\\
\mbox{}\hspace{3em}
differential equations
  				                       \dotfill\ 15
\\[0.5ex]
\mbox{}\hspace{3em}
2.1.1. Necessary condition of existence of  cylindrical 
last multiplier
			                                          \dotfill\ 16
\\[0.5ex]
\mbox{}\hspace{3em}
2.1.2.
Criterion  of existence of cylindrical last multiplier
					\dotfill\ 16
\\[0.5ex]
\mbox{}\hspace{3em}
2.1.3.    Functionally independent cylindrical 
last multipliers
					\dotfill\ 18
\\[1ex]
\mbox{}\hspace{1em}
2.2.  
Autonomy  and cylindricality of last multipliers
for  total differential system 
					\dotfill\ 18
\\[0.5ex]
\mbox{}\hspace{3em}
2.2.1. Necessary condition of existence of  
$s\!$-nonautonomous  $(n-k)\!$-cylindrical
\\
\mbox{}\hspace{5.7em}
last multiplier
					\dotfill\ 19
\\[0.5ex]
\mbox{}\hspace{3em}
2.2.2.
Criterion  of existence of 
$s\!$-nonautonomous  $(n-k)\!$-cylindrical
\\
\mbox{}\hspace{5.7em}
last multiplier
					\dotfill\ 20
\\[0.5ex]
\mbox{}\hspace{3em}
2.2.3.    Functionally independent 
$s\!$-nonautonomous  $(n-k)\!$-cylindrical
\\
\mbox{}\hspace{5.7em}
last multipliers
					\dotfill\ 21
\\[1ex]
{\bf  
3. Cylindricality and autonomy of partial integrals}          \dotfill\ 22
\\[0.5ex]
\mbox{}\hspace{1em}
3.1. 
Cylindricality of partial integrals 
for
linear homogeneous system of partial 
\\
\mbox{}\hspace{3em}
differential equations
				                          \dotfill\ 22
\\[0.5ex]
\mbox{}\hspace{3em}
3.1.1. Necessary condition of existence of  cylindrical partial integral
					\dotfill\ 22
\\[0.5ex]
\mbox{}\hspace{3em}
3.1.2.
Criterion  of existence of cylindrical partial integral
					\dotfill\ 22
\\[0.5ex]
\mbox{}\hspace{3em}
3.1.3.    Functionally independent cylindrical partial integrals
					\dotfill\ 24
\\[1ex]
\mbox{}\hspace{1em}
3.2.  
Autonomy  and cylindricality of partial integrals for
total differential system
   					\dotfill\ 25
\\[0.5ex]
\mbox{}\hspace{3em}
3.2.1. Necessary condition of existence of  
$s\!$-nonautonomous  $(n-k)\!$-cylindrical
\\
\mbox{}\hspace{5.7em}
partial integral
					\dotfill\ 26
\\[0.5ex]
\mbox{}\hspace{3em}
3.2.2.
Criterion  of existence of 
$s\!$-nonautonomous  $(n-k)\!$-cylindrical
\\
\mbox{}\hspace{5.7em}
partial integral
					\dotfill\ 27
\\[0.5ex]
\mbox{}\hspace{3em}
3.2.3.    Functionally independent 
$s\!$-nonautonomous  $(n-k)\!$-cylindrical
\\
\mbox{}\hspace{5.7em}
partial integrals
					\dotfill\ 28
\\[1ex]
\mbox{}\hspace{1em}
3.3.  Functional relations between general solutions to 
irreducible\\
\mbox{}\hspace{3em}
 Painlev\'{e} equations
					\dotfill\ 30
\\[1ex]
{\bf   References}                                                                  \dotfill\ 37

\mbox{}
\\[3ex]
\centerline{
\large\bf  
Introduction
}
\\[1.5ex]
\indent
{\bf  
0.1. Generalities
}
\\[1ex]
\indent
The subject of investigation is 
a 
linear homogeneous system of partial 
differential equations                   
\\[2ex]
\mbox{}\hfill                                 
$
{\frak L}_j^{}(x)\;\! y = 0,
\ \  j = 1,\ldots,m,
$
\hfill (\!$\partial$\!)
\\[2.25ex]
where $x\in \K^{n},\,  m\leq n,$
linear differential operators of first order
\\[2ex]
\mbox{}\hfill                       
$
\displaystyle
{\frak L}_j^{}(x) =
\sum\limits_{i=1}^n\, u^{}_{ji}(x)\;\! \partial_{x_i}^{}
$
\  for all 
$x \in G,
\ \ j = 1,\ldots, m,
$
\hfill (0.1)
\\[2.25ex]
have  holomorphic
coordinates $u_{ji}^{}\colon G\to \K,\, 
j = 1,\ldots, m,\, i = 1,\ldots, n,$ 
a domain $G\subset \K^{n},\ \K$ is a field of real $\R$ or complex  $\C$ numbers, and 
a total differential system
\\[2ex]
\mbox{}\hfill                                   
$
dx = X(t,x)\,dt,
$
\hfill (TD)
\\[2.25ex]
where  $t\in\K^m,\, x\in\K^n,\, m\leq n, \ 
dt = {\rm colon}\,(dt_1,\ldots,dt_m),\,
dx = {\rm colon}\,(dx_1,\ldots,dx_n),\ 
n\times m$ ma\-trix $X(t,x)= \|X_{ij}^{}(t,x)\|$ has  
holomorphic elements $X_{ij}\colon \Pi\to \K,\, i = 1,\ldots, n,\, j = 
\linebreak
=1,\ldots, m,$ 
a domain  $\Pi\subset \K^{m+n}.$ 
Under  $m=1$ the system (TD) is an
ordinary differential  system of $n\!$-th order.

With  a purpose of
unambiguous interpretation of the concepts we use,
let's formulate the 
generalities
 of  integrals   theory 
for systems (\!$\partial$\!) and (TD) 
and define the terminology.

Let's consider the operators (0.1) as being
not linearly bound on the domain  $G$
[5, pp. 105 -- 115]. 
At that we proceed from definition  that
linear differential operators of first order are  linearly bound 
on a domain
if  they are linearly depended over the field $\K$
in every point of this domain
[6, pp. 113 -- 114].

A holomorphic scalar function 
$F \colon G^{\,\prime} \to \K$
is called a {\it first integral} on a domain $G^{\,\prime}\subset G$ 
of  system  (\!$\partial$\!) if 
[1, pp. 35 -- 38; 7, pp. 55 -- 94]
\\[2ex]
\mbox{}\hfill                                   
$
{\frak L}_j^{} F(x) = 0
$
\  for all 
$x \in G^{\,\prime},
\ \  j =1,\ldots,m.
$
\hfill (0.2)
\\[2.25ex]
\indent
A  holomorphic scalar function 
$\mu\colon G^{\;\!\prime}\to\K$  
is called a {\it last multiplier}
on a domain $G^{\,\prime}\subset G$ 
of  system  (\!$\partial$\!) if 
[1, pp. 121 -- 124]
\\[2ex]
\mbox{}\hfill                                           
$
{\frak L}_{_{\scriptstyle j}}\;\!\mu(x) = {}-\mu(x)\,
{\rm div}\;\! {\frak L}_ j^{}(x)
$
\  for all 
$x\in G^{\;\!\prime},
\ \  j =1,\ldots,m.
$
\hfill (0.3)
\\[2.25ex]
\indent
We'll say that a holomorphic scalar function 
$w\colon G^{\;\!\prime}\to\K$ 
(a manifold $w(x)=0$\!) 
is a {\it partial  integral}\,\footnote[1]{
For example,  the term {\it second integral} is used and there is also another terminology (see [8]).} 
on a domain $G^{\,\prime}\subset G$ 
(an {\it integral manifold})
of  system  (\!$\partial$\!) if 
\\[2ex]
\mbox{}\hfill                                           
$
{\frak L}_{j}^{}\;\! w(x) =
\Phi_{j}^{}(x)
$
\  for all 
$x\in G^{\;\!\prime},
\ \ j=1,\ldots, m,
$
\hfill (0.4)
\\[2.25ex]
where  $\Phi_{j}\colon G^{\;\!\prime}\to\K,\,
j=1,\ldots,m,$ are such functions that 
\\[2.5ex]
\mbox{}\hfill                                           
$
\Phi_{j}^{}(x)_{\displaystyle |_{w(x)=0}} = 0
$
\  for all 
$x\in G^{\;\!\prime},
\ \ j=1,\ldots,m.
$
\hfill (0.5)
\\[2.25ex]
\indent
The  system (\!$\partial$\!) is called a {\it complete} system  if the
Poisson bracket of any two its operators (0.1) 
can be represented as a linear combination of this operators
[9, p. 117]
\\[2ex]
\mbox{}\hfill                                     
$
\displaystyle
\bigl[{\frak L}_j^{}(x),
{\frak L}_{\zeta}^{}(x)\bigr] =
\sum\limits_{l=1}^m\, A_{j\zeta l}^{}(x)
{\frak L}_l^{}(x)
$
\  for all 
$x \in G,
\ \ 
j=1,\ldots,m,\  \zeta=1,\ldots,m,
$
\hfill (0.6)
\\[2.5ex]
 with holomorphic coefficients
$A_{j\zeta l}^{}\colon G\to \K,\  
j=1,\ldots,m,\, \zeta=1,\ldots,m,\,l=1,\ldots,m.$

If the Poisson brackets of  operators (0.1) are symmetric, that is,
\\[2ex]
\mbox{}\hfill                                     
$
\displaystyle
\bigl[{\frak L}_j^{}(x),{\frak L}_{\zeta}^{}(x)\bigr] =
\bigl[{\frak L}_{\zeta}^{}(x),{\frak L}_j^{}(x)\bigr]
$
\  for all 
$x \in G,
\ \ j=1,\ldots,m,\  \zeta=1,\ldots,m,
$
\hfill (0.7)
\\[2.23ex]
then the system  (\!$\partial$\!)  is called  a \textit{jacobian} system [7, p. 62].

The symmetry (0.7) of the Poisson brackets 
of  operators (0.1) is equivalent to
\\[2ex]
\mbox{}\hfill                                     
$
\bigl[{\frak L}_j^{}(x),
{\frak L}_{\zeta}^{} (x)\bigl] = {\frak O}
$
\  for all 
$x \in G,
\ \ 
j=1,\ldots,m,\ \zeta=1,\ldots, m,
$
\hfill (0.8)
\\[2.25ex]
where ${\frak O}$ is the  null operator. 
The identity  (0.8) is 
the  representation
(0.6) with the coefficients  $A_{j\zeta l}^{}(x)=0$ for all 
$x \in G, \  j=1,\ldots,m,\, \zeta=1,\ldots,m,\,l=1,\ldots,m.$
Therefore a jacobian system
 (\!$\partial$\!)  is  complete [7, p. 62].

A differential system
\\[2ex]
\mbox{}\hfill                                      
$
\partial_{_{\scriptstyle x_j}} y = {\frak M}_j^{}(x)\;\! y,
\ \ j =1,\ldots,m,
$
\hfill (N\!$\partial$\!)
\\[2.35ex]
where  $x\in \K^n,\, m<n,$
linear differential operators of first order
\\[2ex]
\mbox{}\hfill                                   
$
\displaystyle
{\frak M}_j^{}(x) = \sum\limits_{s=m+1}^n\,
u_{js}^{}(x)\;\!  \partial_{x_s}^{}
$
\  for all 
$
x \in G, \ \ j = 1,\ldots,m,
\hfill 
$
\\[2.5ex]
have  holomorphic coefficients 
$u_{js}^{}\colon G\to\K,\,
j=1,\ldots,m,\, s=m+1,\ldots,n,$ 
is called a  \textit{normal} 
linear homogeneous  system of partial differential equations [7, p. 64].

Let's note that a complete normal system is  jacobian [7, p. 65].

The  complete  system ($\!\partial$\!) 
by means of  linear nonsingular  on the domain $G$ 
change  of operators  (0.1)
can be reduced  to the complete normal system
(at that only the restriction of the domain $G$  may happen) [7, p. 66].

Let the  complete  system
($\!\partial$\!) be such a system that 
the square matrix $\widehat{u}$
of order  $m,$ which is formed  by  $m$ first columns
of $m\times n$ matrix 
$u(x)=\bigl\|u_{ji}^{}(x)\bigr\|$ for all $x\in G,$
is nonsingular on the domain  $G.$
Then, the complete  system  ($\!\partial$\!) 
can be reduced  to the complete normal system
in the form of (N$\!\partial$\!), 
at that, in the neighbourhood  of any point $x$ from the domain $G,$
where  ${\rm det}\, \widehat{u}(x)\ne 0,$
this systems are integral equivalent
[1, pp. 47 -- 48].

We'll call a subdomain $H$ of the domain $G$ 
a \textit{normalization domain} of system ($\!\partial$\!)
if this system  in a neighbourhood  of any point of the domain $H$ 
can be reduced  to the integral equivalent
complete normal system
[1, p.  48].

A  set of functionally independent
on a domain $G^{\;\!\prime}\subset G$ first integrals 
$
F_l \colon G^{\;\!\prime}\to \K,$
$l= 1,\ldots, k,$
of system (\!$\partial$\!) 
is called a \textit{basis of  first integrals} 
(or an \textit{integral basis})
 on the domain  $G^{\;\!\prime}$ 
of  system (\!$\partial$\!)
if any first integral $\Psi \colon G^{\;\!\prime} \to \K$
of this system can be 
represented as 
$
\Psi(x) = \Phi\bigl(F_{1}(x), \ldots,
F_{k}(x)\bigr)$ for all $x \in G^{\;\!\prime},$
where $\Phi$ is some 
holomorphic function on the codomain of vector function
$F\colon x\to \bigl(F_{1}(x), \ldots, F_{k}(x)\bigr)$ for all $x \in G^{\;\!\prime}.$
At that, the number $k$ is called the \textit{dimension} 
of the basis of  first integrals
on the domain  $G^{\;\!\prime}$ 
of   system ($\!\partial$\!)
[1, p.  38; 7, p. 70; 10, pp. 523 -- 525].

A complete linear homogeneous  system of partial differential equations 
($\!\partial$\!)
on a neighbourhood  of any point from 
its normalization domain
has 
a  basis of  first integrals of  the dimension  $n - m$  [1, p.  51].

Every incomplete system ($\!\partial$\!)
on a domain $G$
can be reduced  to an integral equivalent
complete system
[7, pp.  243 -- 245].

We'll call a number  $\delta$  a \textit{defect} of an incomplete system ($\!\partial$\!)
if this system on the domain  $G$ 
can be  reduced  to an integral equivalent
complete system by addition
of $\delta$ equations as
\\[2ex]
\mbox{}\hfill                                        
$
\bigl[
{\frak L}_{{}_{\scriptstyle j_{{}_{\scriptsize \nu}}}}(x),
{\frak L}_{{}_{\scriptstyle l_{\mu}}}(x)\bigr]\;\! y=0, \ \ \
\Bigl[{\frak L}_{{}_{\scriptstyle \alpha_{{}_{\scriptsize \xi}}}}(x),
\bigl[{\frak L}_{{}_{\scriptstyle j_{{}_{\scriptsize \nu}}}}(x),
{\frak L}_{{}_{\scriptstyle l_{\mu}}}(x)
\bigr]\Bigr]\;\! y=0,
\hfill                                        
$
\\[2.25ex]
\mbox{}\hfill                                        
$
\Bigl[{\frak L}_{{}_{\scriptstyle \beta_{{}_{\scriptsize
\zeta}}}}(x),
\Bigl[{\frak L}_{{}_{\scriptstyle
\alpha_{{}_{\scriptsize \xi}}}}(x), \bigl[{\frak L}_{{}_{\scriptstyle
j_{{}_{\scriptsize \nu}}}}(x), {\frak L}_{{}_{\scriptstyle
l_{\mu}}}(x) \bigr]
\Bigr]\Bigr]\;\!  y=0, \, \ldots\,,
\hfill
$
\\[2.75ex]
$
\nu=1,\ldots, m_1, \   \mu=1,\ldots, m_2,
\  \xi=1,\ldots, m_3, \  \zeta= 1,\ldots, m_4,
\  \ldots\;\!,\  m_s\leq m,
\ 
s  = 1, 2,\ldots\;\!, 
$
$
\{1,\ldots,m\}\owns j_{\nu}, l_{\mu}, \alpha_{\xi},
\beta_{\zeta},\ldots $
[1, pp.  42 -- 43].

Let's  agree on a complete system has the defect $\delta=0.$
Then, one can say that every  system (\!$\partial$\!)  has the  defect $\delta$ and $0\leq \delta \leq n - m.$

The system  ($\!\partial$\!)
with defect $\delta$ 
on a neighbourhood  of any point from its 
normalization domain
has a basis of first integrals 
of dimension $n - m - \delta$
[1, p.  51].

The system ($\!\partial$\!)
is complete if and only if
on a neighbourhood  of every point from its any 
normalization domain
it has a basis of first integrals 
of dimension $n - m.$

The  system (TD) is called  {\it completely solvable} on the domain  $\Pi^\prime\subset \Pi$
if in every  point  $(t_0,x_0)\in \Pi^\prime$
for system (TD) 
a solution to the Cauchy problem 
with initial conditions $(t_0,x_0)$ is unique
 [1, p.  17].
In case, when  $\Pi^\prime= \Pi,$
we'll say that system  (TD) is completely solvable.

The system (TD) is  completely solvable
if and only if the Frobenius conditions [1, pp.  17 -- 25;  11, pp. 290 -- 302]
are satisfied:
\\[2ex]
\mbox{}\hfill                     
$
\displaystyle
\partial_{t_j^{}}^{}
X_{i\zeta}^{}(t,x) +
\sum\limits_{\xi=1}^n\,
X_{\xi j}^{}(t,x)
\partial_{x_\xi^{}}^{}
X_{i\zeta}^{}(t,x) =
\partial_{t_\zeta^{}}^{}
X_{ij}^{}(t,x) +
\sum\limits_{\xi=1}^n\,
X_{\xi \zeta}^{}(t,x)
\partial_{x_\xi^{}}^{}
X_{ij}^{}(t,x)
\hfill 
$
\\
\mbox{}\hfill (0.9)                    
\\
\mbox{}\hfill                     
for all 
$(t,x)\in \Pi, \ \ i=1,\ldots, n,
\ j=1,\ldots, m,\ \zeta=1,\ldots, m.
\hfill
$
\\[2.5ex]
\indent
The system (TD) induces  $m$  
linear differential operators of first order
\\[2ex]
\mbox{}\hfill                                        
$
\displaystyle
{\frak X}_j^{}(t,x)  =
\partial_{t_j^{}}^{}
+ \sum\limits_{i=1}^{n}\,
X_{ij}^{}(t,x)
\partial_{x_i}^{}
$ 
\ for all 
$(t,x)\in \Pi, \ \ j=1,\ldots,m.
$
\hfill (0.10)
\\[2.25ex]
Then, the Frobenius conditions (0.9) 
via the Poisson brackets can be written as the system of operator identities
\\[2ex]
\mbox{}\hfill                                        
$
\bigl[{\frak X}_j^{}(t,x),
{\frak X}_\zeta^{}(t,x)\bigr] = {\frak O}
$
\ for all 
$(t,x)\in \Pi, \ \ j=1,\ldots,m,\  \zeta=1,\ldots,m.
$
\hfill (0.11)
\\[2.5ex]
\indent
The system  (TD) is the Pfaff system of equations
\\[2ex]
\mbox{}\hfill                                        
$
\omega_i^{} (t,x)=0, \ \ i=1,\ldots,n,
$
\hfill (0.12)
\\[2ex]
with linear differential forms
\\[2ex]
\mbox{}\hfill                                        
$
\displaystyle
\omega_i^{} (t,x)=dx_i^{}-
\sum\limits_{j=1}^{m}\,
X_{ij}^{}(t,x)\, dt_j^{}
$
\ for all 
$
(t,x)\in \Pi,  \ \ i=1,\ldots,n,
$
\hfill (0.13)
\\[2.25ex]
which have holomorphic  coefficients
$X_{ij}\colon \Pi\to \K,\, i=1,\ldots,n,\, j=1,\ldots,m.$

The Frobenius conditions (0.9) of the
complete solvability of  system (TD)
(the Frobenius conditions of  closure  
of the Pfaff system of equations (0.12))
[11, pp. 299 -- 301; 12, pp. 91; 13]
can be written
via differential 1-forms (0.13) 
as the system of
exterior differential identities
\\[2ex]
\mbox{}\hfill                                                  
$
d\omega_i^{}(t,x)\wedge
\Bigl(\,\mathop{\wedge}\limits_{\xi=1}^{n}
\omega_\xi^{}(t,x)\Bigr)
= 0$
\ for all 
$(t,x)\in \Pi, \  \ i = 1,\ldots,n.
\hfill 
$
\\[2.5ex]
\indent
A holomorphic scalar function 
$F\colon \Pi^{\prime}\to\K$
is called a \textit{first integral}
on the domain $\Pi^{\prime} \subset  \Pi$
\linebreak
of system  (TD) if 
the differential of the function $F$ 
by virtue of  system (TD) 
vanishes identically
on the 
domain  $\Pi^{\prime},$ that is,
\\[2ex]
\mbox{}\hfill                                           
$
dF(t,x)_{\displaystyle |_{\rm (TD)}} =0
$
\ for all
$
(t,x) \in \Pi^{\prime}.
$
\hfill (0.14)
\\[2.25ex]
\indent
By means of 
linear differential operators (0.10) 
the identity  (0.14) can be written as the system of identities
[1, p. 26]
\\[2ex]
\mbox{}\hfill                                              
$
{\frak X}_j\;\! F(t, x) = 0
$
\ for all
$
(t,x) \in \Pi^{\prime},
\ \ j = 1,\ldots, m.
$
\hfill (0.15)
\\[2ex]
\indent
A holomorphic scalar function 
$F\colon \Pi^{\prime}\to\K$ is a first integral
on the domain $\Pi^{\prime}\subset \Pi$
of the completely solvable system  (TD) if and only if
this function keeps a constant value
along any solution $x\colon t \to x(t)$ for all $t\in T\subset \K^m$ 
to  system (TD) such that 
$(t,x(t))\in \Pi^\prime$ for all $t\in T,$ that is,
\\[1ex]
\mbox{}\hfill
$
F(t,x(t))=C$
\ for all 
$
t\in T,
\ \   C={\rm const}.
\hfill
$
\\[2.25ex]
\indent
A holomorphic scalar function 
$\mu\colon \Pi^\prime \to\K$
is called
a \textit{last multiplier}
on the domain $\Pi^{\prime}\subset \Pi$
of  system  (TD) if 
[1, pp. 129 -- 131; 9, pp. 117 -- 130]
\\[2ex]
\mbox{}\hfill                                           
$
{\frak X}_j^{}\;\!\mu(t,x) = {}-\mu(t,x)\,
{\rm div}\;\!{\frak X}_j^{}(t,x)
$
\ for all
$
(t,x)\in \Pi^\prime, \  \ j =1,\ldots,m.
$
\hfill (0.16)
\\[2.25ex]
\indent
We'll call a holomorphic scalar function 
$w\colon \Pi^\prime \to\K$ 
(a manifold $w(t,x) = 0$\!) 
 a {\it partial  integral} on the domain 
$\Pi^{\prime}\subset \Pi$
(an {\it integral  manifold})
of  system  (TD) if 
[1, pp. 161 -- 163]
\\[2ex]
\mbox{}\hfill                                           
$
{\frak X}_{j}^{}\;\! w(t,x) =
\Phi_{j}^{}(t,x)$
\ for all 
$
(t,x)\in \Pi^{\prime},
\ \ j=1,\ldots, m,
$
\hfill {\rm (0.17)}
\\[2.5ex]
where  $\Phi_{j}\colon \Pi^{\;\!\prime}\to\K,\,
j=1,\ldots,m,$ are such the functions that 
\\[2.5ex]
\mbox{}\hfill                                           
$
\Phi_{j}^{}(t,x)_{\displaystyle |_{w(t,x)=0}} = 0
$
\ for all 
$
(t,x)\in \Pi^{\prime},
\ \ j=1,\ldots, m.
$
\hfill {\rm (0.18)}
\\[2.5ex]
\indent
A holomorphic scalar function 
$w\colon \Pi^\prime \to\K$ 
(a  manifold $w(t,x) = 0$\!) 
is a partial integral
on the domain $\Pi^{\prime}\subset \Pi$
(an integral  manifold)
of the completely solvable system  (TD) if and only if
the function $w$ 
vanishes identically
along any solution $x\colon t \to x(t)$ for all $t\in T\subset \K^m$ 
to  system (TD) such that 
$(t,x(t))\in \Pi^\prime$ for all $t\in T,$ that is,
\\[2ex]
\mbox{}\hfill
$
w(t,x(t))=0$
\ for all 
$
t\in T.
\hfill
$
\\[2.25ex]
\indent
If a last multiplier $\mu$
of  system  (\!$\partial$\!)
(of  system (TD)\;\!) defines a  manifold
$\mu= 0,$ then it is a  partial integral
of this system.

Indeed, the system of  identities  (0.3) is the 
system of  identities (0.4), where
$\Phi_j(x)=
\linebreak
={}-\mu(x)\, {\rm div}\,{\frak L}_j(x)$ for all 
$x\in G^{\;\!\prime},\, j=1,\ldots, m,$ and 
therefore the conditions  (0.5) 
are satisfied.
Similarly, the system of  identities (0.16) 
is the system of  identities  (0.17) with the property (0.18).
\vspace{0.25ex}

If a last multiplier $\mu$
of system  (\!$\partial$\!)
(of  system (TD)\;\!) defines a  manifold
$1/\mu= 0,$ then it is an   integral  manifold of this system.
\vspace{0.25ex}

A set of functionaly independent
on the domain $\Pi^{\prime}\subset \Pi$ first integrals 
$
F_l \colon \Pi^{\prime}\to \K,$
$l= 1,\ldots, k,$
of  system (TD)
is called a \textit{basis of  first integrals} 
(or an \textit{integral basis})
 on the domain  $\Pi^{\prime}$ 
of system (TD)
\vspace{0.25ex}
if any first integral $\Psi \colon \Pi^{\prime} \to \K$
of this system can be 
represented as 
$
\Psi(t,x) = \Phi\bigl(F_{1}(t,x), \ldots,
F_{k}(t,x)\bigr)$ for all $(t,x) \in \Pi^{\prime},
$
\vspace{0.25ex}
where $\Phi$ is some 
holomorphic function on the codomain of vector function
\vspace{0.25ex}
$F\colon (t,x)\to \bigl(F_{1}(t,x), \ldots, F_{k}(t,x)\bigr)$ for all 
$(t,x) \in \Pi^{\prime}.$
At that, the number $k$ is called the \textit{dimension} 
of the basis of  first integrals
on the domain  $\Pi^{\prime}$
of  system (TD)
[1, p.  29].
\vspace{0.25ex}

A completely solvable system  (TD)
on a neighbourhood  of any point from 
the domain  $\Pi$
has a  basis of  first integrals of dimension $n$ 
[1, p.  34].
\vspace{0.25ex}

To construct a basis of  first integrals for  system  (TD)
without taking into consideration
a solvability of  Cauchy problem 
we'll be based on integrals theory  
of 
 linear homogeneous system of partial differential equations.

A normal  linear homogeneous system of partial differential equations
\\[2ex]
\mbox{}\hfill  
$
{\frak X}_j (t,x)\;\! y=0, \ \ j=1,\ldots,m,
$
\hfill (0.19)
\\[2.25ex]
is  associated to the  total differential system  (TD).

Directly from definitions of the first integral 
both for the total differential system  and 
for the linear homogeneous system of partial differential equations
we establish the basic integral connection between the systems (TD) and (0.19)
[1, pp.  53 -- 56].
\vspace{0.2ex}

A holomorphic function 
$F\colon (t,x)\to F(t,x)$ for all $(t,x)\in \Pi^\prime$
is a first integral on the domain
$\Pi^\prime\subset \Pi$
of  system (TD) if and only if 
it is a first integral on the domain
$\Pi^\prime $ 
of  system   (0.19) which is  associated to system (TD).

In accordance with
 this property the systems (TD) and (0.19) have 
locally common basis of first integrals
of quite concrete dimension [14].
\vspace{0.2ex}

A total differential system  (TD)
in a neighbourhood  of every point from 
normalization domain of  the
associated to this system 
normal 
linear homogeneous system of partial differential equations (0.19) 
has a basis of first integrals
of dimension $n-\delta,$
where $\delta$ is the defect of  system (0.19), $0\leq \delta\leq n.$
\\[2ex]
\indent
{\bf  
0.2.  Problem definition
}
\\[1.25ex]
\indent
Ordinary differential  system of $n\!$-th order
\\[2ex]
\mbox{}\hfill                                         
$
\dfrac{dx}{dt} = f(t,x),
$
\hfill {\rm (D)}
\\[2.5ex]
where $t\!\in\! \K,\, x\!\in\!\K^n,
\ 
\vspace{1ex}
\dfrac{dx}{dt}= {\rm colon}\Bigl(\dfrac{dx_1}{dt}\,,
\ldots, \dfrac{dx_n}{dt}\Bigr),$
the vector function
$f(t,x)={\rm colon}\,(f_{1}(t,x),
\linebreak
\ldots,f_{n}(t,x))$
\vspace{0.25ex}
has holomorphic coordinates
$f_i\colon \Pi\to \K,\, i=1,\ldots, n,$
the domain $\Pi\subset \K^{n+1},$ 
has a basis of first integrals of dimension $n$
[15, pp.  156 -- 159; 16, pp. 256 -- 263].
\vspace{0.25ex}

Autonomous ordinary differential  system of $n\!$-th order
\\[2ex]
\mbox{}\hfill                                         
$
\dfrac{dx}{dt} = f(x),
$
\hfill {\rm (AD)}
\\[-2ex]

\newpage

\noindent
where
$t\in \K,\, x\in \K^n,
\
\vspace{1ex}
\dfrac{dx}{dt} = {\rm colon}\Bigl(\dfrac{dx_1}{dt}\,,\ldots, \dfrac{dx_n}{dt}\Bigr),$
the vector function  
$f(x)={\rm colon}\,(f_{1}(x),
\linebreak
\ldots,
f_{n}(x))$
has holomorphic coordinates 
$f_i\colon G\to \K,\, i=1,\ldots, n,$
the domain  $G\subset \K^{n},$
in the $n\!$-dimensional basis has
$n-1$ autonomous first integrals which are functionaly independent on the domain $G$
[17, pp.  161 -- 169].

V.I. Mironenko has studied  [18] whether nonautonomous system (D) 
can have autonomous first integrals.
The possibility of existence of autonomous 
partial integrals and autonomous last multipliers
for the nonautonomous system (D) was studied in  [2].
Moreover, in that article 
I solved the problem  whether both  the system (D)
and the system (AD)
can have the first integrals, the last multipliers
and the partial integrals which are the functions depending
on several  $r$ variables $x_{\xi_\tau}^{},\, \xi_\tau \in \{1,\ldots,n\},\, 
\tau=1,\ldots,r,\, 1\leq r \leq n,$ but not necessarily
depending on all variables $x_i,\, i=1,\ldots,n.$

Let's give specific example in which the existence of the autonomous
partial integral depending on $r<n$
variables is useful in stability theory. 

The nonautonomous real ordinary differential  system of third order
\\[2ex]
\mbox{}\hfill                     
$
\dfrac{dx}{dt} = {}-x\bigl( x^2+(2+\sin t)y^2\bigr)-y(e^t x +t^2 xz+t y^2),
\hfill                     
$
\\[2ex]
\mbox{}\hfill                     
$
\dfrac{dy}{dt} = x(e^t x +t^2 xz+t y^2)-y\bigl( x^2+(2+\sin t)y^2\bigr),
\hfill                     
$
\\[2ex]
\mbox{}\hfill                     
$
\dfrac{dz}{dt} = z\bigl( x^2+(2+\sin t)y^2+z^2\bigr)
\hfill 
$
\\[2.25ex]
has the autonomous partial integral $w\colon (x,y,z)\to x^2+y^2$ for all $(x,y,z) \in\R^3.$
This partial integral specifies the isolated point  $x=y=0$
in the subspace $\R^2$ of the phase space $\R^3$ and satisfies the hypotheses of Rumyanzev's theorem [19, p. 29]. 
Therefore the zero solution  $x=y=z=0$ to this system 
is asymptotically stable with respect to  $(x,y).$
At the same time, this solution is unstable by Chetaev's theorem
[19, pp. 19 -- 20] (one should take $V(x,y,z)={}-x^2-y^2+z^2$\!)
\vspace{0.5ex}

Let's define the problem of  existence
of  first integrals,  last multipliers
and  partial integrals for  system ($\!\partial\!$) which
are the functions depending
on several variables 
$x_i,\, i=1,\ldots,n,$
and for  system (TD)
which
are the functions depending
on several  independend variables  $t_j,
\linebreak
j=1,\ldots,m,$
and several  dependend variables $x_i,\, i=1,\ldots,n.$
At that, we'll use the terms  and theoretical foundations
from introduction.
\\[4ex]
\centerline{
\large\bf  
1. Cylindricality and autonomy of first integrals
}
\\[2.5ex]
\centerline{
\bf
1.1. 
Cylindricality of first integrals for
linear homogeneous system 
}
\\[0.2ex]
\centerline{
\bf
of partial differential equations 
}
\\[1ex]
\indent
{\bf Definition  1.1.}
{\it
We'll say  that a  first integral  $F$
on a domain $G^{\;\!\prime}\subset G$
of  system  {\rm($\!\partial\!$)}
is {\boldmath$(n - k)\!$}\textit{\textbf{-cylindrical}}
if the function  $F$ depends only on 
\vspace{0.5ex}
 $k,\, 0\leq k\leq n,$ variables  $x_1,\ldots,x_n.$
}

Let's define the problem of  existence
for  system ($\!\partial\!$) a  $(n-k)\!$-cylindrical first integral
\\[2ex]
\mbox{}\hfill                                           
$
F\colon x \to F({}^kx)
$
\ for all 
$x\in G^{\;\!\prime}\subset G,
\quad 
{}^k x=(x_1,\ldots,x_k).
$
\hfill (1.1)
\\[3ex]
\indent
{\bf
1.1.1.
Necessary condition of existence of  cylindrical first integral.}
According to the definition of first integral, the function  (1.1) 
will be the first integral on the domain $G^{\;\!\prime}\subset G$ of  system (\!$\partial$\!)
if and only if 
\\[2ex]
\mbox{}\hfill                                           
$
{}^k {\frak L}_j^{}  F({}^kx) = 0
$
\ for all $x\in G^{\;\!\prime}, \ \ 
j =1,\ldots, m,
$
\hfill (1.2)
\\[2.2ex]
where the
linear differential operators of first order
\\[2ex]
\mbox{}\hfill     
$
\displaystyle
{}^k {\frak L}_j^{}(x)=
\sum\limits_{\xi =1}^k\,
u_{j\xi}^{} (x)\;\! \partial_{x_\xi^{}}^{}
$
\ for all 
$x\in G, \ \ 
j =1,\ldots, m.
$
\hfill (1.3)
\\[2.5ex]
\indent
Concerning the sets of functions 
$
{}^k U_j =
\bigl\{ u_{j1}^{}(x),\ldots, u_{jk}^{}(x)\bigl\}, \
j =1,\ldots, m,
$
the system of identeties
(1.2)  means that the functions of each set ${}^k U_j,\,
j =1,\ldots, m,$ are linearly dependent with respect to  
variable $x_p$ on the domain $G^{\;\!\prime}$
under any fixed values of 
variables  $x_i,\,
i=1,\ldots, n,\, i\ne p.$ 
It holds true  under each fixed index $p =k+1,\ldots,n.$
Therefore the Wronskians of each set  ${}^k U_j,\, j=1,\ldots,m,$
with respect to variables $x_p ,\, p =k+1,\ldots,n,$ 
vanish identically on the domain $G^{\;\!\prime},$ that is, 
the system of identities
\\[2.25ex]
\mbox{}\hfill                                 
$
{\rm W}_{x_p}^{}
\bigl({}^k u^{j} (x)\bigr) = 0
$
\ for all 
$x\in G^{\;\!\prime},
\ \ j=1,\ldots, m,\ p =k+1,\ldots, n,
$
\hfill (1.4)
\\[2ex]
is satisfied. Here vector functions 
$
{}^k u^{j}\colon x\to
\bigl(u_{j1}^{}(x),\ldots,u_{jk}^{}(x)\bigr)$
for all $x\in G, \ j=1,\ldots,m,$
and ${\rm W}_{x_p}$ are the Wronskians with respect to $x_p,\, p =k+1,\ldots, n.$

So the 
{\sl necessary test of  existence of $(n - k)\!$-cylindrical first integral for 
linear homogeneous system of partial differential equations} 
 is proved.
 
{\bf Theorem 1.1.}
{\it
The system of identities
{\rm (1.4)} is a necessary condition of  existence  
of $(n - k)\!$-cylindrical first integral  {\rm (1.1)}  
for system} (\!$\partial$\!).
\\[2ex]
\indent
{\bf
1.1.2. Criterion  of existence of cylindrical first integral.}
Let $m\times n$ matrix  $u(x)=\|u_{ji}^{}(x)\|$ for all $x\in G$ 
of  system (\!$\partial$\!)  satisfies  the  conditions (1.4).
Let's write the functional system 
\\[2.15ex]
\mbox{}\hfill                            
$
{}^k u^j(x) \, {}^k\varphi= 0, \ 
j=1,\ldots,m,
\quad
\partial^{\xi}_{{}_{\scriptstyle x_p}}
\;{}^k u^j(x)\, {}^k\varphi = 0,
\hfill
$
\\[-0.5ex]
\mbox{}\hfill (1.5)
\\[-0.3ex]
\mbox{}\hfill
$
j =1,\ldots, m, \ \
p =k+1,\ldots,n, \ \
\xi =1,\ldots,k-1,
\hfill
$
\\[2.5ex]
where a vector function
$
{}^k\varphi\colon x\to
\bigl(\varphi_1({}^kx),\ldots,\varphi_k({}^kx)\bigr)$
for all $x\in G^{\;\!\prime}$
is unknown,
vector functions
$
{}^k u^{j}\colon x\to
\bigl(u_{j1}^{}(x),\ldots,u_{jk}^{}(x)\bigr)$ for all $x\in G, \ j=1,\ldots,m.$
Let's introduce a Pfaffian equation 
\\[2ex]
\mbox{}\hfill                                       
$
{}^k\varphi({}^kx)\, d\ {}^kx = 0.
$
\hfill (1.6)
\\[2.5ex]
\indent
{\bf Theorem  1.2} ({\sl 
criterion  of existence of $(n - k)\!$-cylindrical first integral
for linear homogeneous system of partial differential equations}).
{\it
For  system  {\rm (\!$\partial$\!)} 
to have $(n - k)\!$-cylindrical first integral
{\rm(1.1)}
it is necessary and sufficient that there exists a vector function
${}^k\varphi,$ 
satisfying functional system {\rm(1.5)}, 
such that the function {\rm(1.1)} is the general integral of the Pfaffian equation 
{\rm(1.6)} on the domain $\widetilde{G}^{\;\!k}$
which is the natural projection of domain  $G^{\;\!\prime}$ on 
coordinate subspace $O\ {}^kx.$
}

{\sl Proof. Necessity}. Let system  (\!$\partial$\!)
has  $(n - k)\!$-cylindrical first integral (1.1)
on the domain  $G^{\;\!\prime}.$ Then the identities (1.2) are satisfied:
\\[2ex]
\mbox{}\hfill
$
\displaystyle
\sum\limits_{\xi =1}^k\,
u_{j\xi}^{}(x)\;\!
\partial_{x_\xi^{}}^{} F({}^kx) = 0
$
\ for all 
$x\in G^{\;\!\prime}, \ \ j=1,\ldots,m.
\hfill
$
\\[2ex]
\indent
By differentiating this identities  $k-1$ times with respect to $x_{p},\, p=k+1,\ldots,n,$  
we conclude  that an 
extension 
on the domain $G^{\;\!\prime}$ of the function
\\[2ex]
\mbox{}\hfill
$
{}^k\varphi\colon {}^kx \to
\bigl(
\partial_{x_1}^{}F({}^kx),\ldots, \partial_{x_k}^{}F({}^kx)\bigr)
$
\ for all ${}^k x\in \widetilde{G}^{\;\!k}
\hfill
$
\\[2ex]
is a solution to the functional system (1.5).
From this it also follows that 
the function  (1.1) is    a 
general integral on the domain $\widetilde{G}^{\;\!k}\subset \K^k$
of the  Pfaffian equation   (1.6).
\vspace{0.3ex}

{\sl Sufficiency}. Let the  vector function
$
{}^k\varphi\colon x\to {}^k\varphi({}^kx)
$ for all $x\in G^{\;\!\prime}$
be the solution to the functional system (1.5) and the Pfaffian equation  (1.6)
which is constructed on its base has a 
general integral (1.1)
 on the domain  $\widetilde{G}^{\;\!k}\subset \K^k.$
Then, the system of  identities 
\\[2ex]
\mbox{}\hfill                                     
$
\partial_{x_\xi^{}}^{} F({}^kx) -
\mu({}^kx)\;\!
\varphi_{\xi}^{}({}^kx) = 0
$
\ for all 
${}^k x\in \widetilde{G}^{\;\!k}, \ \ 
\xi =1,\ldots, k,
$
\hfill (1.7)
\\[2ex]
is satisfied,
where  $\mu\colon {}^k x\to \mu({}^kx)$ for all ${}^k x\in \widetilde{G}^{\;\!k}$ is 
a  holomorphic integrating multiplier of the Pfaffian equation  (1.6) 
which corresponds to its general integral (1.1) 
on the domain  $\widetilde{G}^{\;\!k}.$

Hence, taking into account that the function ${}^k \varphi$ is
the solution to the functional system (1.5)
we receive the system of identities (1.2).
Therefore the function (1.1) is   
an $(n - k)\!$-cylindrical first integral
on the dimain $G^{\;\!\prime}$  of   system (\!$\partial$\!). \k
\vspace{0.5ex}

{\bf Example 1.1.}
Let's consider a
linear homogeneous system of partial differential equations
\\[2ex]
\mbox{}\hfill            
$
{\frak L}_1(x)y=0,
\ \
{\frak L}_2(x)y=0,
$
\hfill (1.8)
\\[2ex]
which is constructed on the base of 
being not holomorphically linearly bound on the space $\K^4$
linear differential operators of first order
\\[2ex]
\mbox{}\hfill
$
{\frak L}_1(x)= x_1x_2\partial_{x_1}^{}
- x_1^2\partial_{x_2}^{}  + (x_1 + x_2^2 + x_3^2 - x_4^2)\partial_{x_3}^{}  +
(x_3^2 + x_4^2)\partial_{x_4}^{}
$
\ for all 
$x\in\K^4,
\hfill
$
\\[2ex]
\mbox{}\hfill
$
{\frak L}_2(x)=
x_2^2\partial_{x_1}^{}  -  x_1x_2\partial_{x_2}^{}  +
(x_3^2 + x_4^2)\partial_{x_3}^{} +
(x_1 + x_2^2 + x_3^2 + x_4^2)\partial_{x_4}^{}
$
\ for all 
$x\in\K^4.
\hfill
$
\\[2.75ex]
\indent
Let's find for system (1.8)
a 2-cylindrical first integral
\\[2ex]
\mbox{}\hfill            
$
F\colon x\to F (x_1,x_2)
$
\ for all 
$
x\in G^{\;\!\prime}\subset  \K^4.
$
\hfill (1.9)
\\[2ex]
\indent
The Wronskians of the sets of functions
$
{}^2 U_1=\{x_1x_2, {}-x_1^2\}
$
and 
$
{}^2 U_2=\{x_2^2, {}-x_1x_2\}
$
with respect to $x_3$ and  $x_4$  
vanish  identically on the space $\K^4.$ 
Therefore the  necessary conditions (Theorem 1.1)
of  existence of 2-cylindrical first integral (1.9)  for 
 system (1.8) are satisfied.

Let's write the functional system  (1.5):
\\[2ex]
\mbox{}\hfill
$
x_1x_2\,\varphi_1 +  ({}-x_1^2)\varphi_2 =  0,
\ \ \
x_2^2\,\varphi_1 + ({}-x_1x_2)\varphi_2 =  0,
\hfill
$
\\[2ex]
\mbox{}\hfill
$
\partial_{x_3}^{} (x_1x_2)\,\varphi_1 +  \partial_{x_3}^{}
({}-x_1^2)\varphi_2 =  0,
\ \ \
\partial_{x_3}^{} x_2^2\,\varphi_1 + \partial_{x_3}^{} ({}-x_1x_2)\varphi_2 =  0,
\hfill
$
\\[2ex]
\mbox{}\hfill
$
\partial_{x_4}^{} (x_1x_2)\,\varphi_1 +  \partial_{x_4}^{}({}-x_1^2)\varphi_2 =  0,
\ \ \
\partial_{x_4}^{} x_2^2\,\varphi_1 + \partial_{x_4}^{} ({}-x_1x_2)\varphi_2 =  0.
\hfill
$
\\[2.5ex]
\indent
From this we get the system 
\\[1.5ex]
\mbox{}\hfill
$
x_1x_2\,\varphi_1 -x_1^2\,\varphi_2 =  0,
\ \ \
x_2^2\,\varphi_1 -x_1x_2\,\varphi_2 =  0
\hfill
$
\\[2ex]
and  receive from it the equation
\\[0ex]
\mbox{}\hfill
$
x_2\,\varphi_1 -x_1\,\varphi_2 =  0.
\hfill
$
\\[2.5ex]
\indent
The solution on the space   $\K^4$ 
to this equation are, for example, the functions
\\[2ex]
\mbox{}\hfill
$
\varphi_1\colon  x\to x_1
$
\
for all 
$x\in\K^4
$
\ \ \ and \ \ \
$
\varphi_2\colon  x\to x_2
$
\
for all 
$x\in\K^4.
\hfill
$
\\[2ex]
\indent
Let's write on the base of this solution the Pfaffian equation
\\[2ex]
\mbox{}\hfill
$
x_1\, dx_1 +x_2\, dx_2= 0,
\hfill
$
\\[2ex]
which has the general integral $F\colon (x_1,x_2)\to  x_1^2 + x_2^2$\, for all  $(x_1,x_2)\in\K^2.$

Therefore the system  (1.8) on the space  $\K^4$ has the
2-cylindrical first integral  
\\[2ex]
\mbox{}\hfill       
$
F\colon x\to  x_1^2 + x_2^2
$
\ for all  
$
x\in\K^4.
$
\hfill (1.10)
\\[2.25ex]
\indent
Since the Poisson bracket
\\[2ex]
\mbox{}\hfill
$
[{\frak L}_{1}(x),{\frak L}_{2}(x)] =
{}-x_{2}(x_{1}^{2}+x_{2}^{2})\partial_{x_1}^{}+
x_{1}(x_{1}^{2}+x_{2}^{2})\partial_{x_2}^{} \,+
\hfill
$
\\[2ex]
\mbox{}\hfill
$
+\,
(x_{1}+2x_{1}x_{4}-x_{3}^{2}-x_{4}^{2}+
2x_{1}x_{2}^{2}+2x_{2}^{2}x_{4}+4x_{3}^{2}x_{4}-
2x_{3}x_{4}^{2}+4x_{4}^{3})\partial_{x_3}^{} \,+
\hfill
$
\\[2ex]
\mbox{}\hfill
$
+\,
(x_{1}x_{2}+2x_{1}x_{3}-2x_{1}x_{4}-2x_{3}^{2}-
2x_{1}^{2}x_{2}+2x_{2}^{2}x_{3} -
2x_{2}^{2}x_{4}+2x_{3}^{3}-4x_{3}x_{4}^{2})\partial_{x_4}^{}
\
$
for all  
$x\in\K^4
\hfill
$
\\[2.35ex]
is not a linear combination on the space  $\K^4$ of operators
${\frak L_1}$ and ${\frak L_2}$
the system  (1.8) is not complete.
Then, a basis of first integrals of  system  (1.8)
consists of no more than one first integral  
(accurate within the functional expression).

Thus  the 
2-cylindrical first integral   (1.10) of  system  (1.8)
forms its integral basis on the space $\K^4.$
\\[2ex]
\indent
{
\bf
1.1.3.    Functionally independent cylindrical first integrals.
}
\vspace{0.5ex}

{\bf Theorem  1.3.}
{\it
Let the functional system  {\rm(1.5)} has 
$q$ 
not  linearly bound
on the domain $G^{\;\!\prime}\subset G$  solutions
}
\\[1.75ex]
\mbox{}\hfill                                           
$
{}^k\varphi^\gamma\colon x\to
{}^k\varphi^\gamma({}^kx)
$
\ for all 
$x\in G^{\;\!\prime},
\ \  \gamma=1,\ldots q,
$
\hfill (1.11)
\\[2.25ex]
{\it
where the vector ${}^k\varphi^\gamma=(
\varphi^\gamma_1,\ldots, \varphi^\gamma_k),$
and the Pfaffian equations
}
\\[2.25ex]
\mbox{}\hfill                                           
$
{}^k\varphi^\gamma({}^kx)\, d\ {}^kx = 0,
\ \ \gamma =1,\ldots, q,
$
\hfill (1.12)
\\[2ex]
{\it
which are constructed on the base of this solutions
have correspondingly general integrals}
\\[2ex]
\mbox{}\hfill                                           
$
F_{\gamma}\colon
{}^kx\to F_{\gamma}({}^kx)
$
\ for all 
${}^kx \in \widetilde{G}^{\;\!k}\subset \K^k,
\ \ \gamma =1,\ldots, q,
$
\hfill (1.13)
\\[2ex]
{\it
on the domain $\widetilde{G}^{\;\!k}$ which is the 
natural projection of  domain $G^{\;\!\prime}$ on  coordinate
subspace $O\ {}^kx.$ Then, this general integrals are functionally independent
on the domain  $\widetilde{G}^{\;\!k}.$
}
\vspace{0.3ex}

{\sl Proof.} By virtue of the system of identities  (1.7) 
\\[2ex]
\mbox{}\hfill
$
\partial_{x_\xi^{}}^{}F_\gamma({}^kx) =
\mu_\gamma^{}({}^kx)\;\!
\varphi_{{}_{\scriptstyle \xi }}^\gamma({}^kx)
$
\ for all ${}^kx \in \widetilde{G}^{\;\!k},
\ \ \xi =1,\ldots, k,
 \ \gamma =1,\ldots, q.
\hfill
$
\\[2.3ex]
Therefore the Jacobi's matrix
\\[2.2ex]
\mbox{}\hfill
$
J\bigl(F_{\gamma}({}^kx);\, {}^kx\bigr)=
\bigl\|	
\mu_\gamma^{}({}^kx)\;\!\varphi_{{}_{\scriptstyle \xi }}^\gamma({}^kx)\bigr\|_{q\times k}
$
for all 
${}^kx \in \widetilde{G}^{\;\!k}.
\hfill
$
\\[2.3ex]
\indent
Since the  vector functions (1.11)
are not  linearly bound 
on the domain $\widetilde{G}^{\;\!k}$
the rank of Jacobi's matrix
$
{\rm rank}\, J\bigl(F_{\gamma}( {}^kx);\,{}^kx\bigr) = q
$
for all ${}^kx$
from the domain $\widetilde{G}^{\;\!k}$ 
perhaps with the exception of  point set of 
$k\!$-dimensional zero measure.
So
the  general integrals 
(1.13) of the Pfaffian equations  (1.12)
 are functionally independent
on the domain  $\widetilde{G}^{\;\!k}.$  \k
\vspace{0.5ex}

The Theorem 1.3 (taking into account the Theorem 1.2) let us to find
a quantity of  functionally independent
$(n-k)\!$-cylindrical first integrals of  system {\rm($\!\partial\!$)}.
\\[3ex]
\centerline{
\bf
1.2.  First integrals of {\boldmath $s$}-nonautonomous completely solvable 
}
\\[0ex]
\centerline{
\bf
 total differential systems
}
\\[-1ex]

{\bf Definition 1.2.}
{\it
We'll say that system  {\rm (TD)} is 
{\boldmath$s\!$}\textit{\textbf{-nonautonomous}}
if all functions-elements   $X_{ij}\colon \Pi\to \K,\,
i=1,\ldots,n,\,j=1,\ldots,m,$ of the matrix  $X$ depend on  $x$ and only on  $s,\,0\leq s\leq m,$ 
independent variables $t_1,\ldots,t_m.$
}
\vspace{0.5ex}

Without loss of generality 
we'll consider that
the $s\!$-nonautonomous system (TD)
has  such the functions-elements
$X_{ij}\colon \Pi\to \K,\, i=1,\ldots,n,\,j=1,\ldots,m,$
of matrix  $X$ that 
depend only on  $x$ and on first  $s$ 
independent variables $t_1,\ldots,t_s,$ that is,
\\[2ex]
\mbox{}\hfill                                        
$
dx = X({}^st,x)\, dt,
$
\hfill (TDs)
\\[2.25ex]
where
\vspace{0.3ex}
$
{}^st = (t_1,\ldots ,t_s), \
dt = {\rm colon}\,(dt_1,\ldots ,dt_m), \ 
dx = {\rm colon}\,(dx_1,\ldots ,dx_n),\
n\times m$
matrix 
$X({}^st,x)=
\left\|X_{ij}({}^st,x)\right\|
$
\vspace{0.3ex}
has holomorphic elements
$X_{ij}\colon ({}^st,x)\to X_{ij}({}^st,x)\;\; \forall ({}^st,x)\in \Pi^{s+n},
$
$
i=1,\ldots,n,\,j=1,\ldots,m,$
the domain  $\Pi^{s+n}\subset \K^{s+n}, \,  0\leq s\leq m\leq n.$
\vspace{0.3ex}

Under 
$s=0$ the system  (TDs)  is the autonomous system (ATD).
\vspace{0.3ex}

\newpage

{\bf Definition 1.3.}
{\it
We'll say that a first integral $F$
on a domain  $\Pi^\prime\subset \Pi$ of  
system {\rm(TD)} is
{\boldmath$s\!$}\textit{\textbf{-nonautonomous}}
if the  function   $F$  depends on  $x$ and only on  $s,\,0\leq s\leq m,$ 
independent variables $t_1,\ldots,t_m.$
}
\vspace{0.4ex}

If $s=0,$ then an $s\!$-nonautonomous first integral of  system (TD) 
is an autonomous first integral of  system (TD).

Let ${}^s X$ be  the matrix which is formed from  the $n\times m$ matrix 
$X({}^s t,x)$ by deletion  of first $s$
columns.
\vspace{0.3ex}

{\bf Theorem  1.4.}
{\it
If the rank of  matrix ${}^sX({}^st,x)$ of the completely solvable system {\rm (TDs)} 
is equal to $r$ on the domain $\Pi^{s+n},$ then it has on this domain exactly $n-r$ 
functionally independent $s\!$-nonautonomous first integrals
$
F_\gamma\colon \Pi^{s+n}\to \K,
\,\gamma = 1,\ldots,n-r.
$
}
\vspace{0.5ex}

{\sl Proof}. Let $x\colon t\to x(t;C)$ for all $t\in T\subset \K^m$
be the  solutions to the completely solvable system (TDs).

Without loss of generality 
we'll consider that the first  $r$  rows
of the functional matrix ${}^sX$
form the matrix of rank $r$ (one can always get that by renumbering  of dependent variables).
Then, the first  $r$ components $x_{_l}, \, l =1,\ldots,r,$
of the solutions are functionally independent
on the domain $T$ relative to variables $t_{\xi}, \, \xi =s+1,\ldots, m,$ and 
the rest components $x_{\rho}, \, \rho =
\linebreak
=r+1,\ldots, n,$ 
of the solutions are functionally dependent on first $r$  components
on the domain $T$
relative to variables $t_{\xi}, \, \xi=s+1,\ldots, m.$
So
\\[2ex]
\mbox{}\hfill
$
x_\rho (t) = \Phi_\rho \bigl({}^st,{}^kx(t);C\bigr)
$
\ for all 
$t\in T, \ \ 
\rho =r+1,\ldots, n,
\hfill
$
\\[2.25ex]
where ${}^k x = (x_1,\ldots,x_k),$ and the functions
$\Phi_\rho\colon T\to \K,\, \rho=r+1,\ldots, n,$ are holomorphic.

Taking into account the functional independence 
of functions $x_l^{},\, l=1,\ldots,r,$
relative to $t_{s+1},\ldots,t_{m}$ 
on the domain $T$ 
from  the equalities
$
x_{l}^{} = x_{l}^{}(t;C),\, l =1,\ldots,r,
\,\  
\vspace{0.3ex}
x_\rho = 
\linebreak
=\Phi_\rho({}^st,{}^kx;C),
\ \rho=r+1,\ldots, n,$
\vspace{0.5ex}
by fixation of  arbitrary vector  $C$  by means of vectors
$
C^{i}=(\delta_{i1}C_1, \ldots,\delta_{in}C_n),
\
i=1,\ldots, n,
$
\vspace{0.3ex}
where  $\delta_{ij}$ is  a Kronecker symbol,
we find  $r$ not $s\!$-nonautonomous and $n-r$ 
$s\!$-nonautonomous
functionally independent on the domain  $\Pi^{\prime}\subset\Pi$ 
first  integrals of the system (TDs). \k
\vspace{0.5ex}

Let's  pay attention to  coordination relative to independent variables $t_1,\ldots,t_s$  
between $s\!$-nonautonomy of  system (TDs)
and $s\!$-nonautonomy of first integrals in  Theorem 1.4.
\vspace{0.2ex}

Theorem 1.5 and Theorem 1.6 are the corollaries of  Theorem 1.4
for the autonomous total differential systems. 
\vspace{0.2ex}

{\bf Theorem 1.5.}
{\it
The completely solvable system  {\rm (ATD)} has 
exactly $n-r,$
where   $r =
\linebreak
={\rm rank}\,X(x)$ 
for all $x\in G^{\,\prime},$ 
functionally independent on the domain $G^{\,\prime}\subset G$ autonomous 
\vspace{0.5ex}
first integrals.
}

{\bf  Theorem 1.6.}
{\it
The completely solvable system {\rm (ATD)}  does not have 
the
autonomous first integrals if and only if 
${\rm rank}\,X(x) = n$
for all $x$ from the domain 
$G$ 
perhaps with the exception of  point set of 
$n\!$-dimensional zero measure.
}
\vspace{0.75ex}

{\bf Example 1.2.}
The completely solvable  total differential system 
\\[2ex]
\mbox{}\hfill                                 
$
dx_1 = dt_1, \ \ \
dx_2 = dt_2, \ \  \
dx_3 = \partial_{x_1}^{}g(x_1,x_2)\,dt_1
+ \partial_{x_2}^{}g(x_1,x_2)\,dt_2,
$
\hfill (1.14)
\\[2.5ex]
where a
\vspace{0.3ex}
holomorphic scalar function $g\colon \widetilde{G}\to \K,\,
 \widetilde{G}\subset \K^2,$ 
has $n-r=3-{\rm rank}\, X(x)=3-2=
\linebreak
=1$
autonomous first integral
$F\colon x \to g(x_1,x_2) - x_3$ for all $x\in G=
\vspace{0.5ex}
 \widetilde{G}\times \K$ (Theorem 1.5).

The functionally independent 
first integrals
\\[2ex]
\mbox{}\hfill                                   
$
F_1 \colon (t,x) \to t_1 - x_1$\,
for all $(t,x)\in\Pi,
\quad
F_2 \colon (t,x) \to t_2 - x_2$\,
for all $(t,x)\in\Pi,
\hfill            
$
\\[0.75ex]
\mbox{}\hfill  (1.15)          
\\[-1.75ex]
\mbox{}\hfill            
$
F_3 \colon (t,x) \to g(x_1,x_2) - x_3$\,
for all 
$(t,x)\in\Pi
\hfill 
$
\\[2.5ex]
forms an integral basis on the domain
$\Pi =\K^3\times  \widetilde{G}$ of system (1.14).

\newpage

\mbox{}
\\[-2ex]
\centerline{
\bf
1.3.
Autonomy  and  cylindricality of first integrals for
  total differential system}
\\[2ex]
\indent
{\bf Definition  1.4.}
{\it
We'll say  that a first integral  $F$
on a domain $\Pi^{\;\!\prime}\subset \Pi$
of  system  {\rm(TD)}
is {\boldmath$(n - k)\!$}\textit{\textbf{-cylindrical}}
if the function  $F$ depends on $t$ and only on 
 $k,\, 0\leq k\leq n,$ dependent variables  $x_1,\ldots,x_n.$
}
\vspace{0.5ex}

Let's define the problem of  existence
for  system (TD) an  
$s\!$-nonautonomous
$(n-k)\!$-cylindrical first integral
\\[2ex]
\mbox{}\hfill                                           
$
F\colon (t,x) \to F({}^st,{}^kx)
$
\ for all 
$(t,x)\in \Pi^{\prime}\subset \Pi,
\quad
{}^s t=(t_1,\ldots,t_s),\ \ {}^k x=(x_1,\ldots,x_k).
$
\hfill (1.16)
\\[3.5ex]
\indent
{\bf
1.3.1.
Necessary condition of existence of 
{\boldmath$s\!$}-nonautonomous
{\boldmath$(n-k)\!$}-cylindrical first integral.}
According to the definition of first integral, the function  (1.16) 
will be the first integral on the domain $\Pi^{\;\!\prime}\subset \Pi$ of  system (TD)
if and only if 
\\[2ex]
\mbox{}\hfill                                           
$
{}^{sk}{\frak X}_j^{}F({}^st,{}^kx) = 0
$
\ for all 
$(t,x)\in \Pi^{\prime}, \ \ \
j =1,\ldots,m,
$
\hfill (1.17)
\\[2.25ex]
where the
linear differential operators of first order
\\[2ex]
\mbox{}\hfill         
$
\displaystyle
{}^{sk}{\frak X}_\theta^{}(t,x) =
\partial_{ t_\theta^{}}^{} +
\sum\limits_{\xi=1}^k\,
X_{\xi \theta}^{}(t,x)\partial_{ x_\xi^{}}^{}
$
\ for all 
$(t,x)\in \Pi, \ \ \
\theta =1,\ldots,s,
\hfill
$
\\
\mbox{}\hfill (1.18)
\\[-0.5ex]
\mbox{}\hfill
$
\displaystyle
{}^{sk} {\frak X}_\nu^{}(t,x) =
\sum\limits_{\xi=1}^k\,
X_{\xi \nu}^{}(t,x)
\partial_{x_\xi^{}}^{}
$
\ for all 
$(t,x)\in \Pi, \ \ \
\nu =s+1,\ldots,m.
\hfill
$
\\[2.5ex]
\indent
Concerning the sets of functions 
\vspace{0.5ex}
$
{}^k M_\theta^{} =
\bigl\{1,X_{1\theta}^{}(t,x),\ldots,
X_{k\theta}^{}(t,x)\bigl\}, \,
\theta =1,\ldots, s, \,\
{}^k M_{\nu}^{} =
\linebreak
=
\bigl\{X_{1\nu}^{}(t,x),\ldots,X_{k\nu}^{}(t,x)\bigr\}, \,
\nu=s+1,\ldots,m,
$
\vspace{0.5ex}
the system of identeties
(1.17)  means that: 
the functions of each set ${}^k M_j,\,
j =1,\ldots, m,$ are linearly dependent with respect to  
independent variable $t_\zeta$ on the domain $\Pi^{\;\!\prime}$
\vspace{0.3ex}
under any fixed values of 
independent variables  $t_\gamma,\, \gamma=
\linebreak
=1,\ldots, m,\, \gamma\ne \zeta,$
\vspace{0.4ex}
and  dependent variables  $x_i,\, i=1,\ldots, n;$
and the functions of each set ${}^k M_j,\,
j =1,\ldots, m,$ are 
\vspace{0.3ex}
linearly dependent with respect to  
dependent variable $x_p$ on the domain $\Pi^{\;\!\prime}$
under any fixed values of 
independent variables  $t_\gamma,\, \gamma=1,\ldots, m,$
and  dependent variables  $x_i,\, i=1,\ldots, n,\, i\ne p.$
It holds true  under each fixed index $\zeta =s+1,\ldots,m$
and under each fixed index $p =k+1,\ldots,n.$

Therefore the Wronskians of each set  ${}^k M_j,\, j=1,\ldots,m,$
with respect to 
independent variables  $t_\zeta,\, \zeta=s+1,\ldots, m,$
and  dependent variables  $x_p,\, p=k+1,\ldots, n$
vanish identically on the domain $\Pi^{\;\!\prime},$ that is, 
the system of identities
\\[2.5ex]
\mbox{}\qquad\quad                                 
$
{\rm W}_{ t_\zeta}^{}
\bigl(1,{}^k X^{\theta}(t,x)\bigr) = 0
$
\ for all 
$(t,x)\in \Pi^\prime, \ \ \
\theta =1,\ldots, s,\ \zeta =s+1,\ldots,m,
\hfill
$
\\[1.85ex]
\mbox{}\qquad\quad
$
{\rm W}_{t_\zeta}^{}
\bigl({}^k X^{\nu}(t,x)\bigr) = 0
$
\ for all 
$(t,x)\in \Pi^\prime, \ \ \
\nu=s+1,\ldots, m,
\ \zeta =s+1,\ldots, m,
\hfill
$
\\
\mbox{}\hfill (1.19)
\\[-0.35ex]
\mbox{}\qquad\quad
$
{\rm W}_{x_p}^{}
\bigl(1,{}^k X^{\theta}(t,x)\bigr) = 0
$
\ for all 
$(t,x)\in \Pi^\prime, \ \ \
\theta =1,\ldots,s,\ p =k+1,\ldots, n,
\hfill
$
\\[1.85ex]
\mbox{}\qquad\quad
$
{\rm W}_{x_p}^{}
\bigl({}^k X^{\nu}(t,x)\bigr) = 0
$
\ for all 
$(t,x)\in \Pi^\prime, \ \ \
\nu =s+1,\ldots,m,\ p =k+1,\ldots,n,
\hfill
$
\\[2.5ex]
is satisfied. Here vector functions 
\vspace{0.3ex}
$
{}^k\;\!\!X^{j}\colon (t,x)\to \bigl(X_{1j}(t,x),\ldots,X_{kj}(t,x)\bigr)$
for all $(t,x)\in \Pi, $
$ j=1,\ldots,m,$
${\rm W}_{t_\zeta}$ and ${\rm W}_{x_ p}$
are 
\vspace{0.3ex}
correspondingly the Wronskians with respect to $t_\zeta$ and $x_p,\, 
 \zeta =s+1,\ldots, m,\, p =k+1,\ldots, n.$
\vspace{0.3ex}

So the 
{\sl necessary test of  existence of 
$s\!$-nonautonomous
$(n - k)\!$-cylindrical first integral for  total differential system} is proved.
\vspace{0.2ex}
 
{\bf Theorem 1.7.}
{\it
The system of identities
{\rm (1.19)} is a necessary condition of  existence  
of  $s\!$-nonautonomous $(n - k)\!$-cylindrical first integral  {\rm (1.16)}  
for system} (TD).
\\[2ex]
\indent
{\bf
1.3.2. Criterion  of existence of 
{\boldmath$s\!$}-nonautonomous
{\boldmath$(n-k)\!$}-cylindrical first integral.}
Let $n\times m$ matrix $X$ of  system (TD)  satisfies  the  conditions (1.19).
Let's write the functional system 
\\[1ex]
\mbox{}\hfill                            
$
\psi_\theta^{} +
{}^k\;\!\!X^\theta(t,x)\, {}^k\varphi = 0,
\ \ \theta=1,\ldots,s,
\hfill
$
\\[2ex]
\mbox{}\hfill
$
\partial^{\xi}_{{}_{\scriptstyle t_\zeta}}
\;{}^k\;\!\!X^\theta(t,x)\ 
{}^k\varphi = 0,
\ \ \theta =1,\ldots,s, \
\zeta =s+1,\ldots,m, \
\xi =1,\ldots,k,
\hfill
$
\\[2ex]
\mbox{}\hfill
$
\partial^{\xi}_{{}_{\scriptstyle x_p}}
\;{}^k\;\!\!X^\theta(t,x)\ {}^k\varphi = 0,
\ \ \theta =1,\ldots,s, \
p =k+1,\ldots,n, \
\xi =1,\ldots,k,
\hfill
$
\\[2ex]
\mbox{}\hfill
$
{}^k\;\!\!X^\nu(t,x)\, 
{}^k\varphi = 0,
\ \ \nu =s+1,\ldots,m,
$
\hfill (1.20)
\\[2ex]
\mbox{}\hfill
$
\partial^{\xi}_{{}_{\scriptstyle t_\zeta}}
\,{}^k\;\!\!X^\nu(t,x)\ {}^k\varphi = 0,
\ \ \nu =s+1,\ldots,m,\
\zeta =s+1,\ldots,m,\
\xi =1,\ldots,k-1, 
\hfill
$
\\[2ex]
\mbox{}\hfill
$
\partial^{\xi}_{{}_{\scriptstyle x_p}}
\,{}^k\;\!\!X^\nu(t,x) \ {}^k\varphi
= 0,
\ \
\nu=s+1,\ldots, m,  \
p =k+1,\ldots,n, \
\xi =1,\ldots,k-1,
\hfill
$
\\[2.5ex]
where the vector functions
\vspace{0.5ex}
$
{}^s \psi\colon (t,x)\to
\bigl( \psi_1^{}({}^st,{}^kx),
\ldots,
\psi_s^{}({}^st,{}^kx)\bigr)
$
for all $(t,x)\in \Pi^{\prime}
$
and 
$
{}^k \varphi\colon (t,x)\to
\bigl( \varphi_1^{}({}^st,{}^kx),
\ldots,
\varphi_k^{}({}^st,{}^kx)\bigr)
$
for all $(t,x)\in \Pi^{\prime}$
are unknown,
\vspace{0.75ex}
the vector functions
$
{}^k\;\!\!X^j\colon (t,x)\to
\bigl( X_{1j}(t,x),\ldots,X_{kj}(t,x)\bigr)
$
for all
\vspace{0.5ex}
$(t,x)\in \Pi, \
j =1,\ldots,m,\, 0\leq k\leq n.
$
Let's introduce a Pfaffian equation 
\\[2ex]
\mbox{}\hfill                                       
$
{}^s\psi({}^st,{}^kx)\, d\ {}^st +
{}^k\varphi({}^st,{}^kx)\, d\ {}^kx = 0.
$
\hfill (1.21)
\\[2.35ex]
\indent
{\bf Theorem  1.8} ({\sl 
criterion  of existence of 
$s\!$-nonautonomous
$(n - k)\!$-cylindrical first integral for 
total differential system}).
{\it
For  system  {\rm (TD)} 
to have 
$s\!$-nonautonomous
$(n - k)\!$-cylindrical first integral {\rm(1.16)}
it is necessary and sufficient that there exist the vector functions
${}^s\psi$ and ${}^k\varphi,$ 
satisfying functional system {\rm(1.20)}, 
such that the function {\rm(1.16)} is the general integral of the Pfaffian equation 
{\rm(1.21)} on the domain $\widetilde{\Pi}^{\;\!s+k}$
which is the natural projection of domain  $\Pi^{\;\!\prime}$ on 
coordinate subspace $O\ {}^st\  {}^kx.$
}

{\sl Proof. Necessity}. Let system  (TD)
has the   
$s\!$-nonautonomous
$(n - k)\!$-cylindrical  first integral (1.16)
on the domain  $\Pi^{\;\!\prime}.$ Then, the identities (1.17) are satisfied:
\\[2ex]
\mbox{}\hfill
$
\displaystyle
\partial_{t_\theta^{}}^{} F({}^st,{}^kx) +
\sum\limits_{\xi =1}^k\,
X_{\xi \theta}^{}(t,x)
\partial_{x_\xi^{}}^{}F({}^st,{}^kx) = 0
$
\ for all 
$(t,x)\in \Pi^\prime, \ \
\ \theta=1,\ldots,s,
\hfill
$
\\[1ex]
\mbox{}\hfill
$
\displaystyle
\sum\limits_{\xi=1}^k\,
X_{\xi \nu}^{}(t,x)
\partial_{x_\xi^{}}^{}F({}^st,{}^kx) = 0
$
\ for all 
$(t,x)\in \Pi^\prime, \ \
\ \nu=s+1,\ldots,m.
\hfill
$
\\[2.5ex]
\indent
By differentiating the first $s$ of this identities  $k$ times with respect to 
$t_{s+1},\ldots,t_m$ and $k$ times with respect to $x_{k+1},\ldots,x_n$
and by differentiating the rest $m-s$  identities
$k-1$ times with respect to 
$t_{s+1},\ldots,t_m$ and $k-1$ times with respect to $x_{k+1},\ldots,x_n$
\vspace{0.3ex}
we conclude  that the 
extensions 
on the domain $\Pi^{\;\!\prime}$ of the functions
\vspace{0.5ex}
$
{}^s\psi\colon ({}^st,{}^kx)\to
\partial_{{}^st}^{} F({}^st,{}^kx)
$ for all $({}^st,{}^kx)\in \widetilde{\Pi}^{s+k}
$
and 
$
{}^k\varphi\colon ({}^st,{}^kx)\to
\partial_{{}^kx}^{}F({}^st,{}^kx)
$
for all
$({}^st,{}^kx)\in \widetilde{\Pi}^{s+k}$
\vspace{0.3ex}
is a solution to the functional system (1.20), where
operators
$
\partial_{{}^st}^{} =
(\partial_{t_1},\ldots,\partial_{t_s}),
\ 
\partial_{{}^kx}^{} =
(\partial_{x_1},\ldots, \partial_{x_k}).
$

From this it also follows that 
the function  (1.16) is a  
general integral on the domain $\widetilde{\Pi}^{\;\!s+k}$
of the  Pfaffian equation   (1.21).
\vspace{0.2ex}

{\sl Sufficiency}. Let  vector functions
$
{}^s\psi\colon (t,x) \to  {}^s\psi({}^st,{}^kx), \
{}^k\varphi\colon (t,x)\to {}^k\varphi({}^st,{}^kx)$ for all $(t,x)\in \Pi^{\prime}$
be the solution to the  system (1.20) and the Pfaffian equation  (1.21)
which is constructed on its base has the 
general integral (1.16)
 on the domain  $\widetilde{\Pi}^{\;\!s+k}\subset \K^{s+k}.$
Then, the system of  identities 
\\[2ex]
\mbox{}\hfill                                     
$
\partial_{ {}^st}^{} F({}^st,{}^kx) -
\mu({}^st,{}^kx)\ 
{}^s\psi({}^st,{}^kx) = 0\quad
\forall ({}^st,{}^kx)\in \widetilde{\Pi}^{s+k},
\hfill
$
\\[-0.35ex]
\mbox{}\hfill(1.22)
\\[-0.35ex]
\mbox{}\hfill
$
\partial_{{}^kx}^{} F({}^st,{}^kx) -
\mu({}^st,{}^kx)\ 
{}^k\varphi({}^st,{}^kx) = 0
\quad
\forall ({}^st,{}^kx)\in \widetilde{\Pi}^{s+k},
\hfill
$
\\[2.75ex]
is satisfied,
where  $\mu\colon ({}^st,{}^kx)\to \mu({}^st,{}^kx)$ for all $({}^st,{}^kx)\in \widetilde{\Pi}^{s+k}$  
is  a  holomorphic integrating multiplier of the Pfaffian equation  (1.21) 
which corresponds to its general integral (1.16) 
on the domain  $\widetilde{\Pi}^{\;\!s+k}.$

Taking into account that the functions 
${}^s\psi, \,{}^k\varphi$
are the solution to the functional system (1.20)
we receive the system of identities (1.17).
Therefore the function (1.16) is   
an 
\linebreak
$s\!$-nonautonomous
$(n - k)\!$-cylindrical  first integral
on the dimain $\Pi^{\;\!\prime}$  of   system (TD). \k
\vspace{0.5ex}

For example, the integral basis (1.15)
of the  autonomous completely solvable total differential system  (1.14) (Example 1.2)
consists of the autonomous first integral $F_3$
and two 1-nonautonomous 1-cylindrical  first integrals $F_1$ and $F_2.$
\vspace{0.5ex}

{\bf Example 1.3.}
The autonomous  total differential system 
\\[2ex]
\mbox{}\hfill                                            
$
dx_1 = x_1\,dt_1 + 3x_1\,dt_2,
\quad
dx_2 = (1 + x_1 + 2x_2)\,dt_1 + (x_1 +3x_2)\,dt_2
$
\hfill (1.23)
\\[2ex]
is not completely solvable since  
the 
Poisson bracket 
\\[2ex]
\mbox{}\hfill
$
\bigl[{\frak X}_1(t,x), {\frak X}_2(t,x)\bigr] =
(3-x_1)\partial_{x_2}^{}=
{\frak X}_{12}(t,x)
$
\ for all 
$(t,x)\in\K^4
\hfill
$
\\[2.15ex]
of induced by system (1.23)  linear differential operators 
\\[2ex]
\mbox{}\hfill                                            
$
{\frak X}_1(t,x)=
\partial_{t_1}^{} +
x_1\partial_{x_1}^{}+
(1 + x_1 + 2x_2)\partial_{x_2}^{}
$ \ for all $(t,x)\in\K^4,
\hfill
$
\\[2ex]
\mbox{}\hfill                                            
$
{\frak X}_2(t,x)=
\partial_{t_2}^{} +
3x_1\partial_{x_1}^{}+
(x_1 + 3x_2)\partial_{x_2}^{}
$ \ for all 
$(t,x)\in\K^4
\hfill
$
\\[2ex]
is not the null operator on the any domain from the space
$\K^4.$

By the Frobenius' theorem   the system (1.23) doesn't
have the solutions.

The associated to the system (1.23) 
incomplete normal 
linear homogeneous system of partial differential equations
\\[1ex]
\mbox{}\hfill
$
{\frak X}_1(t,x) y=0,\ \ \ 
 {\frak X}_2(t,x) y=0
\hfill
$
\\[2ex]
can be reduced to the complete system on the space 
$\K^4$ by the addition of single equation ${\frak X}_{12}(t,x)y=0.$
Therefore this system has the defect 
$\delta=1$ and an integral basis of system (1.23) 
consists of $n-\delta=2-1=1$ first integral.

The system (1.23) is autonomous, but  
according to Theorem 1.4 
it has no autonomous first integral.
Indeed, the system of identities (1.19) are not satisfied, because, for example,
the Wronskian
\\[2ex]
\mbox{}\hfill
$
{\rm W}_{x_1} (1+x_1+2x_2,\, x_1+3x_2)=
\left|\!
\begin{array}{cc}
1+x_1+2x_2 & x_1+3x_2
\\[0.5ex]
1 & 1
\end{array}
\!
\right|
=1-x_2
$
\ for all 
$(t,x)\in\K^4
\hfill
$
\\[2ex]
does not  vanish identically on the any domain from the space  $\K^4.$

Let's find a 1-cylindrical  first integral
\\[2ex]
\mbox{}\hfill    
$
F\colon (t,x)\to F(t,x_1)
$
\ for all 
$(t,x)\in \Pi^\prime\subset \K^4
$
\hfill (1.24)
\\[2ex]
of system  (1.23).

The Wronskians of the sets of functions  
${}^1M_1=\{1,x_1\}$ and ${}^1M_2=\{1,3x_1\}$
 with respect to   $x_2$
vanish identically on the space  $\K^4.$
Therefore the  necessary conditions (Theorem 1.7)
of  existence of 1-cylindrical first integral (1.24)  for 
 system (1.23) are satisfied.

The functional system  (1.20) consists of two equations
\\[2ex]
\mbox{}\hfill    
$
\psi_1+x_1\, \varphi_1=0, \ \ \
\psi_2+3x_1\, \varphi_1=0.
\hfill
$
\\[2ex]
\indent
Its solution, for example, is 
\\[2ex]
\mbox{}\hfill   
$
\psi_1\colon (t,x)\to x_1,
 \ \
\psi_2\colon (t,x)\to 3x_1,
\ \ 
\varphi_1\colon (t,x)\to {}-1$
\ for all 
$(t,x)\in  \K^4.
\hfill 
$
\\[2ex]
\indent
The Pfaffian equation which
is constructed  on the base of this solution 
\\[2ex]
\mbox{}\hfill    
$
x_1\, dt_1+3x_1\,dt_2 -dx_1=0
\hfill
$
\\[3ex]
has the general integral
$
F\colon (t,x_1) \to x_1\,e^{^{\scriptstyle {}-(t_1 + 3t_2)}}$
for all $(t,x_1)\in\K^3.$
\vspace{0.5ex}

Therefore the system  (1.23) on the space  $\K^4$ has 
1-cylindrical first integral  
\\[2ex]
\mbox{}\hfill    
$
F\colon (t,x) \to x_1\,e^{^{\scriptstyle {}-(t_1 + 3t_2)}}
$
\ for all $(t,x)\in\K^4.
\hfill
$
\\[2ex]
\indent
This 1-cylindrical first integral  
forms the integral basis on the space $\K^4$
of system (1.23).
\\[2ex]
\indent
{
\bf
1.3.3.    Functionally independent 
{\boldmath$s\!$}-nonautonomous
{\boldmath$(n-k)\!$}-cylindrical first integrals.
}
\vspace{0.5ex}

{\bf Theorem  1.9.}
{\it
Let the functional system  {\rm(1.20)} has 
$q$ 
not  linearly bound on the domain $\Pi^{\;\!\prime}\subset \Pi$  solutions
}
\\[2ex]
\mbox{}\hfill                                           
$
{}^s\psi^\gamma\colon (t,x)\to {}^s\psi^\gamma({}^st,{}^kx)
$
\ for all 
$(t,x)\in\Pi^{\prime},\ \  \ \gamma=1,\ldots,q,
\hfill
$
\\[-.35ex]
\mbox{}\hfill {\rm (1.25)}
\\[-.35ex]
\mbox{}\hfill
$
{}^k\varphi^\gamma\colon (t,x)\to
{}^k\varphi^\gamma({}^st,{}^kx)
$
\ for all 
$(t,x)\in\Pi^{\prime},\ \  \ \gamma=1,\ldots,q,
\hfill
$
\\[2.25ex]
{\it
and the Pfaffian equations
}
\\[2ex]
\mbox{}\hfill                                           
$
{}^s\psi^\gamma({}^st,{}^kx)\  d\,{}^st +
{}^k\varphi^\gamma({}^st,{}^kx)\  d\,{}^kx = 0,
\  \ \gamma=1,\ldots,q,
$
\hfill (1.26)
\\[2ex]
{\it
which are constructed on the base of this solutions
have correspondingly general integrals}
\\[2ex]
\mbox{}\hfill                                           
$
F_{\gamma}\colon
({}^st,{}^kx)\to F_{\gamma}({}^st,{}^kx)
$
\ for all
$({}^st,{}^kx) \in \widetilde{\Pi}^{s+k}\subset \K^{s+k},
\  \ \gamma=1,\ldots,q,
$
\hfill (1.27)
\\[2ex]
{\it
on the domain $\widetilde{\Pi}^{\;\!s+k}$ which is the 
natural projection of  domain $\Pi^{\;\!\prime}$ on  coordinate
subspace $O\ {}^s t\ {}^k x.$ Then, this general integrals are functionally independent
on the domain  $\widetilde{\Pi}^{\;\!s+k}.$
}
\vspace{0.3ex}

{\sl Proof.} By virtue of the system of identities  (1.22) 
\\[2ex]
\mbox{}\hfill
$
\partial_{{}^st}^{}F_\gamma({}^st,{}^kx) =
\mu_\gamma^{}({}^st,{}^kx)\ 
{}^s\psi^\gamma({}^st,{}^kx)
$
\ for all
$({}^st,{}^kx) \in \widetilde{\Pi}^{s+k},
\  \ \gamma=1,\ldots,q,
\hfill
$
\\[2ex]
\mbox{}\hfill
$
\partial_{{}^kx}^{}F_\gamma({}^st,{}^kx) =
\mu_\gamma^{}({}^st,{}^kx)\ 
{}^k\varphi^\gamma({}^st,{}^kx)
$
\ for all
$({}^st,{}^kx) \in \widetilde{\Pi}^{s+k},
\  \ \gamma=1,\ldots,q.
\hfill
$
\\[2.5ex]
\indent
Therefore the Jacobi's matrix
\\[2.35ex]
\mbox{}\hfill
$
J\bigl(F_{\gamma}({}^st,{}^kx);\, {}^st,{}^kx\bigr)=
\bigl\|\Psi({}^st,{}^kx)\Phi({}^st,{}^kx)\bigr\|
$
\ for all
$({}^st,{}^kx) \in \widetilde{\Pi}^{s+k},
\hfill
$
\\[2.25ex]
where the matrix  $\|\Psi\Phi\|$ consists of 
\vspace{0.75ex}
$q\times s$ matrix
$
\Psi({}^st,{}^kx) =
\bigl\|\mu_\gamma^{}({}^st,{}^kx)\;\!
\psi_j^{\gamma}({}^st,{}^kx)\bigr\|
$
for all $({}^st,{}^kx) \in \widetilde{\Pi}^{s+k}$
and  $q\times k$ matrix
$
\Phi({}^st,{}^kx) =
\bigl\|\mu_\gamma^{}({}^st,{}^kx)\;\!
\varphi_i^{\gamma}({}^st,{}^kx)\bigr\|$
\vspace{0.75ex}
for all $({}^st,{}^kx) \in \widetilde{\Pi}^{s+k}.$

Since the vector functions (1.15)
\vspace{0.3ex}
are not  linearly bound 
on the domain $\widetilde{\Pi}^{\;\!s+k}$
the rank of Jacobi's matrix
\vspace{0.3ex}
$
{\rm rank}\, J\bigl(F_{\gamma}({}^st,{}^kx);\,
{}^st,{}^kx\bigr) = q
$
for all $({}^st,{}^kx)$
from the domain
 $\widetilde{\Pi}^{s+k}$ 
perhaps with the exception of  point set of 
$(s+k)\!$-dimensional zero measure.
\vspace{0.3ex}
So
the  general integrals 
(1.27) of the Pfaffian equations  (1.26)
\vspace{0.3ex}
 are functionally independent
on the domain  $\widetilde{\Pi}^{\;\!s+k}.$  \k
\vspace{0.5ex}

The Theorem 1.9 (taking into account the Theorem 1.8) let us to find
a quantity of  functionally independent
$s\!$-nonautonomous $(n-k)\!$-cylindrical  first integrals of  system (TD).
\\[4ex]
\centerline{
\large\bf  
2. Cylindricality and autonomy of last multipliers  
}
\\[2.5ex]
\centerline{
\bf 
2.1. Cylindricality of last multipliers for
linear homogeneous system 
}
\\[.2ex]
\centerline{
\bf 
of partial differential equations
}
\\[1.5ex]
\indent
{\bf Definition  2.1.}
{\it
We'll say  that a  last multiplier  $\mu$
on a domain $G^{\;\!\prime}\subset G$
of  system  {\rm($\!\partial\!$)}
is {\boldmath$(n - k)\!$}\textit{\textbf{-cylindrical}}
if the function  $\mu$ depends only on 
 $k,\, 0\leq k\leq n,$ variables  $x_1,\ldots,x_n.$
}

\newpage

Let's define the problem of  existence
for  system ($\!\partial\!$) an  $(n-k)\!$-cylindrical last multiplier
\\[2ex]
\mbox{}\hfill                                           
$
\mu\colon x \to \mu({}^kx)
$
\ for all 
$x\in G^{\;\!\prime}\subset G,
\quad {}^k x=(x_1,\ldots,x_k).
$
\hfill (2.1)
\\[3ex]
\indent
{\bf
2.1.1.
Necessary condition of existence of  cylindrical last multiplier.}
According to the definition of last multiplier,  the function  (2.1) 
will be the last multiplier on the domain $G^{\;\!\prime}\subset G$ of  system (\!$\partial$\!)
if and only if 
\\[2ex]
\mbox{}\hfill                                          
$
{}^k {\frak L}_{j}^{}\mu ({}^k x) +
\mu ({}^k  x)\, {\rm div}\, u^j(x) = 0
$
\ for all 
$x\in G^{\;\!\prime},
\ \ j = 1,\ldots,m,
$
\hfill (2.2)
\\[2ex]
where the linear differential operators of first order
\vspace{0.25ex}
${}^k {\frak L}_j,\,j = 1,\ldots,m,$
are defined by means of  (1.3),
the vector functions
$u^{j}\colon x\to\bigl(u_{j1}^{}(x),\ldots,u_{jn}^{}(x)\bigr)$ for all $x\in G,\,j = 1,\ldots,m.$
\vspace{0.5ex}

The system of identities (2.2) in the coordinates 
is given by 
\\[2ex]
\mbox{}\hfill     
$
\displaystyle
\sum\limits_{\xi=1}^k\,
u_{j\xi}^{}(x)
\partial_{x_\xi^{}}^{}
\mu ({}^k x) +
\mu({}^k x)\,{\rm div}\, u^j(x) = 0
$ \ for all 
$x\in G^{\;\!\prime},
\ \ j = 1,\ldots,m.
$
\hfill (2.3)
\\[2.5ex]
\indent
Concerning the sets of functions 
\vspace{0.3ex}
$
{}^k D_j =
\bigl\{ u_{j1}^{}(x),\ldots, u_{jk}^{}(x),\;\! {\rm div}\, u^j(x) \bigl\}, \,
j =1,\ldots, m,
$
the system of identeties
(2.3)  means that the functions of each set 
\vspace{0.2ex}
${}^k D_j,\,
j =1,\ldots, m,$ are linearly dependent with respect to  
variable $x_p$ on the domain $G^{\;\!\prime}$
\vspace{0.2ex}
under any fixed values of 
variables  $x_i,\,
i=1,\ldots, n,\, i\ne p.$ 
It holds true  under each fixed index $p =k+1,\ldots,n.$
\vspace{0.3ex}
Therefore the Wronskians of each set  ${}^k D_j,\, j=1,\ldots,m,$
\vspace{0.2ex}
with respect to variables $x_p ,\, p =k+1,\ldots,n,$ 
vanish identically on the domain $G^{\;\!\prime},$ that is, 
the system of identities
\\[2ex]
\mbox{}\hfill                                 
$
{\rm W}_{ x_p}^{}
\bigl({}^k u^{j} (x),\;\! {\rm div}\, u^{j} (x) \bigr) = 0
$
\ for all 
$x\in G^{\;\!\prime},
\ \ j = 1,\ldots,m, \
p =k+1,\ldots,n,
$
\hfill (2.4)
\\[2.25ex]
is satisfied. Here the vector functions 
\vspace{0.2ex}
$
{}^k u^{j}\colon x\to
\bigl(u_{j1}^{}(x),\ldots,u_{jk}^{}(x)\bigr)$
for all $x\in G, \ j=
\linebreak
=1,\ldots,m,$
and ${\rm W}_{x_p}$ are the Wronskians with respect to $x_p,\, p =k+1,\ldots, n.$
\vspace{0.3ex}

So the 
{\sl necessary test of  existence of $(n - k)\!$-cylindrical last multiplier for 
linear homogeneous system of partial differential equations} is proved.
\vspace{0.2ex}
 
{\bf Theorem 2.1.}
{\it
The system of identities
{\rm (2.4)} is a necessary condition of  existence  
of $(n - k)\!$-cylindrical last multiplier  {\rm (2.1)}  
for system} (\!$\partial$\!).
\\[2ex]
\indent
{\bf
2.1.2. Criterion  of existence of cylindrical last multiplier.}
Let the $m\times n$ matrix  $u(x)=\|u_{ji}^{}(x)\|$ for all $x\in G$ 
of  system (\!$\partial$\!)  satisfies  the  conditions (2.4).
Let's write the functional system 
\\[1.5ex]
\mbox{}\hfill                            
$
{}^k u^j(x)\ {}^k\varphi = {}- {\rm div}\,u^{j}(x),
\ \ j = 1,\ldots,m, 
\hfill 
$
\\[-0.5ex]
\mbox{}\hfill (2.5)
\\
\mbox{}\hfill
$
\partial^{\xi}_{{}_{\scriptstyle x_p}}
\;{}^k u^j(x) \
{}^k\varphi
= {}- \partial^{\xi}_{{}_{\scriptstyle x_p}}\,
{\rm div}\,u^{j}(x),
\ \ j = 1,\ldots,m, 
\ p =k+1,\ldots,n, \
\xi =1,\ldots,k-1,
\hfill 
$
\\[2.75ex]
where a vector function
\vspace{0.5ex}
$
{}^k\varphi\colon x\to
\bigl(\varphi_1({}^kx),\ldots,\varphi_k({}^kx)\bigr)$
for all $x\in G^{\;\!\prime}$
is unknown,
the vector functions
$
{}^k u^{j}\colon x\to
\bigl(u_{j1}^{}(x),\ldots,u_{jk}^{}(x)\bigr)$ for all $x\in G, \ j=1,\ldots,m.$
\vspace{0.75ex}

{\bf Theorem  2.2} ({\sl 
criterion  of existence of $(n - k)\!$-cylindrical last multiplier for 
linear homogeneous system of partial differential equations}).
{\it
For  system  {\rm (\!$\partial$\!)} 
to have $(n - k)\!$-cylindrical last multiplier
{\rm(2.1)}
it is necessary and sufficient that there exists a vector function
${}^k\varphi,$ 
satisfying functional system {\rm(2.5)}, 
such that the Pfaffian equation  {\rm(1.6)} 
\vspace{0.25ex}
which is constructed on the base of this vector function
is exact 
on the domain $\widetilde{G}^{\;\!k}$
which is the natural projection of domain  $G^{\;\!\prime}$ on 
coordinate subspace $O\ {}^kx.$
At that,  the last multiplier {\rm(2.1)} of 
system  {\rm (\!$\partial$\!)} is
}
\\[2ex]
\mbox{}\hfill                                          
$
\displaystyle
\mu\colon x\to \exp
\int  {}^k\varphi ({}^k x)\  d\ {}^kx
$
\ for all 
$x\in G^{\;\!\prime}.
$
\hfill (2.6)
\\[2.5ex]
\indent
{\sl Proof. Necessity}. Let system  (\!$\partial$\!)
has the  $(n - k)\!$-cylindrical last multiplier (2.1)
on the domain  $G^{\;\!\prime}.$ Then, the system of identities (2.3) is satisfied.
By means of termwise division of every identity (2.3) by 
$\mu({}^k x)$ we get a new system of identities
\\[2ex]
\mbox{}\hfill
$
\displaystyle
\sum\limits_{\xi=1}^k\,
u_{j\xi}^{}(x)
\partial_{x_\xi^{}}^{}
\ln \mu ({}^k x) +
{\rm div}\, u^j(x) = 0
$ \ for all 
$x\in G_0^{\;\!\prime}\subset G^{\;\!\prime} , \ \  
\ j=1,\ldots,m.
\hfill
$
\\[2.5ex]
\indent
By differentiating this identities  $k-1$ times with respect to $x_{p},\, p=k+1,\ldots,n,$  
we conclude  that an 
extension 
on the domain $G_0^{\;\!\prime}$ of the vector function
\\[2ex]
\mbox{}\hfill                                    
$
{}^k\varphi\colon  {}^k x\to
\bigl(
\partial_{x_1}^{} \ln \mu ({}^k x),\ldots, \partial_{x_k}^{} \ln \mu ({}^k x)\bigr)
$
for all  
${}^k x\in \widetilde{G}_0^{\;\!k}\subset \K^k.
$
\hfill (2.7)
\\[2ex]
is a solution to the functional system (2.5).

The Pfaffian equation  {\rm(1.6)} 
\vspace{0.25ex}
which is constructed on the base of the  vector function (2.7)
is exact on the domain $\widetilde{G}_0^{\;\!k}.$

From (2.7) it follows that 
$(n - k)\!$-cylindrical last multiplier $\mu$ of system (\!$\partial$\!)
is constructing on the domain $G_0^{\;\!\prime}$
on the base of solutions to the system (2.5) by formula (2.6).

By restriction
the domain $G^{\;\!\prime}$ to its codomain 
$G_0^{\;\!\prime}$ we conclude that the
necessary condition of Theorem 2.2 is satisfied.
\vspace{0.3ex}

{\sl Sufficiency}. Let the  vector function
\vspace{0.3ex}
${}^k\varphi$
be a solution to the functional system (2.5) and the Pfaffian equation  (1.6)
which is constructed on its base 
is exact on the domain $\widetilde{G}^{\;\!k}\subset \K^k.$
Then, the   identities 
\\[2ex]
\mbox{}\hfill
$
\displaystyle
\partial_{x_\xi^{}}^{}
\int  {}^k\varphi ({}^k x)\  d\ {}^kx=
\varphi_{\xi}^{} ({}^k x)
$
\ for all  
${}^k x\in \widetilde{G}^{\;\!k},
\ \ \xi=1,\ldots,k,
\hfill
$
\\[2ex]
are satisfied.

Taking into account that the vector function ${}^k \varphi$ is
a solution to the functional system (2.5)
we receive the system of identities (2.2) for the function (2.6).

Hence,  the function (2.6) is an   
 $(n - k)\!$-cylindrical last miltiplier  of   system (\!$\partial$\!). \k
\vspace{0.5ex}

{\bf Example 2.1.}
Consider the 
linear homogeneous system of partial differential equations
\\[2ex]
\mbox{}\hfill         
$
{\frak L}_1(x)y =0, \ \
{\frak L}_2(x)y =0,
$
\hfill (2.8)
\\[2ex]
where the linear differential operators of first order
\\[2ex]
\mbox{}\hfill         
$
{\frak L}_1(x)= x_1x_2\partial_{x_1}^{} +
x_1x_3\partial_{x_2}^{}+x_1x_4\partial_{x_3}^{} +
x_2^2\;\!\partial_{x_4}^{}
$
\ for all 
$x\in \R^4,
\hfill         
$
\\[2ex]
\mbox{}\hfill         
$
{\frak L}_2(x) = x_1x_3\partial_{x_1}^{} +
x_1x_4 \partial_{x_2}^{} + x_1^2\;\!\partial_{x_3}^{}+ x_2^2\;\!\partial_{x_4}^{}
$
\  for all 
$x\in \R^4.
\hfill 
$
\\[2.75ex]
\indent
Let's find for system (2.8)
a 3-cylindrical last multiplier 
\\[2ex]
\mbox{}\hfill         
$
\mu\colon
x\to\mu(x_1)
$ \ for all 
$x\in G^{\;\!\prime}\subset \R^4.
$
\hfill (2.9)
\\[2ex]
\indent
The divergences 
\\[2.1ex]
\mbox{}\hfill
$
{\rm div}\, u^1(x)=
{\rm div}\, {\frak L}_1(x) = \partial_{x_1}^{}(x_1x_2) +
\partial_{x_2}^{}(x_1x_3)+\partial_{x_3}^{}(x_1x_4)+
\partial_{x_4}^{}x_2^2=x_2
$
\ for all 
$x\in \R^4,
\hfill
$
\\[2ex]
\mbox{}\hfill
$
{\rm div}\, u^2(x)=
{\rm div}\, {\frak L}_2(x)= \partial_{x_1}^{}(x_1x_3) +
\partial_{x_2}^{}(x_1x_4)+ \partial_{x_3}^{}x_1^2 + \partial_{x_4}^{}x_2^2 =
x_3
$
\ for all 
$x\in \R^4.
\hfill
$
\\[2.25ex]
\indent
The Wronskians of the sets of functions
$
{}^1 D_1=\{x_1x_2, x_2\}
$
and 
$
{}^1 D_2=\{x_1x_3, x_3\}
$
with respect to $x_2,\,x_3,$ and  $x_4$  
vanish  identically on the space $\R^4\colon$ 
\\[2.5ex]
\mbox{}\hfill
$
{\rm W}_{x_2}^{}(x_1x_2,\, x_2)
=
\left|\!
\begin{array}{cc}
x_1x_2  & x_2
\\[0.5ex]
x_1     & 1
\end{array}
\!\!
\right|
=  0,
\,\ \ 
{\rm W}_{x_3}^{}(x_1x_2,\, x_2)=
{\rm W}_{x_4}^{}(x_1x_2,\, x_2) = 0
$
\ for all 
$x\in \R^4,
\hfill
$
\\[2ex]
\mbox{}\hfill
$
{\rm W}_{x_2}^{}(x_1x_3,\, x_3) = 0,
\ \
{\rm W}_{x_3}^{}(x_1x_3,\,x_3)=
\left|\!
\begin{array}{cc}
x_1x_3 & x_3
\\[0.5ex]
x_1    & 1
\end{array}
\!\!\right|
 = 0,
\ \
{\rm W}_{x_4}^{}(x_1x_3,\, x_3)
= 0
$
\ for all 
$x\!\!~\in\!\!~ \R^4.
\hfill
$
\\[3ex]
\indent
Therefore the  necessary conditions (Theorem 2.1)
of  existence of 3-cylindrical last multiplier  (2.9)  for 
 system (2.8) are satisfied.

Let's write the functional system 
\\[2ex]
\mbox{}\hfill
$
x_1x_2\,\varphi_1 = {}-x_2,
\ \
x_1x_3\,\varphi_1 = {}-x_3,
\ \
x_1\;\!\varphi_1 = {}-1.
\hfill
$
\\[2.5ex]
From this system we find 
$
\varphi_1\colon x\to {}-1/x_1$ for all $x\in\R^4\backslash\{x\colon x_1= 0\}.$

Since 
\\[1.5ex]
\mbox{}\hfill
$
\displaystyle
\exp\int\dfrac{dx_1}{{}-x_1} = \dfrac{C}{|x_1|}
$
\ for all 
$x_1\in\R\backslash\{0\}
\quad
(C>0),
\hfill
$
\\[2.75ex]
the function 
$
\mu \colon x \to {} - 1/x_1
$ for all 
$x\in  G^{\;\!\prime}\subset  \R^4\backslash \{x \colon x_1= 0\}$
is a 3-cylindrical last multiplier  of  system (2.8) (Theorem 2.2).
\\[2ex] 
\indent
{\bf 2.1.3.    Functionally independent cylindrical 
last multipliers.}
The method which is proposed in Theorem 2.2
can be used to construct  
the functionally independent $(n-k)\!$-cylindrical
last multipliers of system (\!$\partial$\!).
\vspace{0.3ex}

{\bf Theorem  2.3.}
{\it
Let the functional system  {\rm(2.5)} has 
$q$ 
not  linearly bound 
on the domain $G^{\;\!\prime}\subset G$  solutions
{\rm (1.11)}
and the corresponding  Pfaffian equations {\rm (1.12)}
are exact 
on the domain $\widetilde{G}^{\;\!k}$
which is the natural projection of domain  $G^{\;\!\prime}$ on 
coordinate subspace $O\ {}^kx.$
Then,  the $(n-k)\!$-cylindrical
last multipliers  of 
system  {\rm (\!$\partial$\!)} 
}
\\[2ex]
\mbox{}\hfill
$
\displaystyle
\mu_\gamma\colon x\to
\exp\int {}^k \varphi^\gamma ({}^k x )\  d\ {}^k x
$
\ for all 
$x\in G^{\;\!\prime},
\ \ \ \gamma=1,\ldots,q,
\hfill
$
\\[2.25ex]
{\it 
are functionally independent
on the domain  $G^{\;\!\prime}.$
}
\vspace{0.2ex}

{\sl Proof.} From Theorem 2.2 it follows that the last multipliers
$\mu_\gamma^{},\, \gamma=1,\ldots,q,$ of system  (\!$\partial$\!)  
are of indicated structure.
 
From representations
\\[2ex]
\mbox{}\hfill
$
\partial_{x_\xi^{}}^{}
\ln\mu_\gamma ({}^k x ) =
\varphi^{\gamma}_{{}_{\scriptstyle \xi}} ({}^k x)
$
\ for all 
${}^k x\in \widetilde{G}^{\;\!k},
\ \ \ \xi =1,\ldots,k, \ \gamma =1,\ldots,q,
\hfill
$
\\[2.25ex]
it follows that the  Jacobi's matrix
\\[2ex]
\mbox{}\hfill
$
J\bigl(\ln\mu_\gamma({}^k x);\,{}^k x \bigr) =
\bigl\| \varphi^{\gamma}_{{}_{\scriptstyle \xi}} ({}^k x)
\bigr\|_{q\times k}
$
\ for all 
${}^k x\in \widetilde{G}^{\;\!k}.
\hfill
$
\\[2.5ex]
\indent
Since the solutions (1.11) 
to the functional system  (2.5)
are not  linearly bound 
on the domain $G^{\;\!\prime}$
the rank of Jacobi's matrix
$
{\rm rank}\, J\bigl(\ln\mu_\gamma({}^kx);\,{}^kx\bigr) = q
$
nearly everywhere on the domain  $\widetilde{G}^{\;\!k}.$

So  the $(n-k)\!$-cylindrical
last multipliers  
$\mu_\gamma,\,\gamma=1,\ldots,q,$
of system  (\!$\partial$\!) 
are functionally independent
on the domain  $G^{\;\!\prime}.$  \k
\\[3ex]
\centerline{
\bf 
2.2.  
Autonomy  and cylindricality of last multipliers
for  
}
\\
\centerline{
\bf 
total differential system 
}
\\[1.25ex]
\indent
{\bf Definition 2.2.}
{\it
We'll say that a last multiplier  $\mu$
on a domain  $\Pi^\prime\subset \Pi$ of  
system {\rm(TD)} is
{\boldmath$s\!$}\textit{\textbf{-nonautonomous}}
if the  function   $\mu$  depends on  $x$ and only on  $s,\,0\leq s\leq m,$ 
independent variables $t_1,\ldots,t_m.$
If $s=0,$ then  a last multiplier
$\mu\colon (t,x) \to \mu(x)$ for all $(t,x) \in \Pi^{\prime}$
of  system {\rm (TD)} is \textit{\textbf{autonomous}}.
}
\vspace{0.5ex}

{\bf Definition  2.3.}
{\it
We'll say  that a  last multiplier  $\mu$
on a domain $\Pi^{\;\!\prime}\subset \Pi$
of  system {\rm (TD)}
is {\boldmath$(n - k)\!$}\textit{\textbf{-cylindrical}}
if the function  $\mu$ depends 
on  $t$ and only on 
 $k,\, 0\leq k\leq n,$ dependent variables   $x_1,\ldots,x_n.$
}
\vspace{0.5ex}

Let's define the problem of  existence
for  system (TD) an  
$s\!$-nonautonomous
$(n-k)\!$-cylindrical last multiplier
\\[2ex]
\mbox{}\hfill                                           
$
\mu \colon (t,x) \to \mu({}^st,{}^kx)
$
\ for all 
$(t,x)\in \Pi^{\prime}\subset \Pi,
\quad
{}^s t=(t_1,\ldots,t_s),\ \ {}^k x=(x_1,\ldots,x_k).
$
\hfill (2.10)
\\[3.5ex]
\indent
{\bf
2.2.1.
Necessary condition of existence of 
{\boldmath$s\!$}-nonautonomous
{\boldmath$(n-k)\!$}-cylindrical last multiplier.}
According to the definition of last multiplier, the function  (2.10) 
will be the last multiplier on the domain $\Pi^{\;\!\prime}\subset \Pi$ of  system (TD)
if and only if 
\\[2.25ex]
\mbox{}\hfill                                          
$
{}^{sk}{\frak X}_j^{}
\mu ({}^s t,{}^k x) +
\mu ({}^s t,{}^k  x)\, {\rm div}_x^{} X^j(t,x) = 0
$
\ for all 
$(t,x)\in \Pi^{\prime},
\ \  j = 1,\ldots,m,
$
\hfill (2.11)
\\[2.2ex]
where the 
\vspace{0.5ex}
linear differential operators of first order
${}^{sk}{\frak X}_j,\,j = 
1,\ldots,m,$
are defined by means of  (1.18),
the vector functions
\vspace{0.5ex}
$X^{j}\colon (t,x)\to \bigl(X_{1j}(t,x),\ldots,X_{nj}(t,x)\bigr)$ for all $(t,x)\in \Pi,\   j=
\linebreak
=1,\ldots,m,$
the divergence
${\rm div}_x^{} X^j(t,x)=
\sum\limits_{i=1}^n \partial_{x_i}^{} X_{ij}^{}(t,x)$ for all $(t,x)\in \Pi,\   j=1,\ldots,m.$
\vspace{1.5ex}

The system of identities (2.11) in the coordinates 
is given by 
\\[2ex]
\mbox{}\hfill       
$
\displaystyle
\partial_{t_\theta^{}}^{}
\mu({}^s t,{}^k x) +
\sum\limits_{\xi=1}^k\,
X_{\xi \theta}^{}(t,x)
\partial_{x_\xi^{}}^{}
\mu ({}^s t,{}^k x) +
\mu({}^s t,{}^k x)\,{\rm div}_x^{} X^\theta(t,x) = 0$
\ for all 
$(t,x) \in \Pi^{\prime},
\hfill
$
\\[1.5ex]
\mbox{}\hfill 
$
\displaystyle
\sum\limits_{\xi=1}^k\, X_{\xi \nu}^{}(t,x)
\partial_{x_\xi^{}}^{}
\mu ({}^s t,{}^k x\bigr) +
\mu ({}^s t,{}^k x )\,{\rm div}_x^{}X^\nu(t,x) = 0$ 
\ for all $(t,x) \in \Pi^{\prime},
$
\hfill  (2.12)     
\\[1.5ex]
\mbox{}\hfill       
$
\theta =1,\ldots,s,
\ \
\nu =s+1,\ldots,m.
\hfill
$
\\[3ex]
\indent
Concerning the sets of functions 
\vspace{1ex}
$
{}^k B_\theta^{} =
\bigl\{1,X_{1\theta}^{}(t,x),\ldots,
X_{k\theta}^{}(t,x),\;\!
{\rm div}_x^{}X^\theta(t,x) \bigl\}, \,
\theta =
\linebreak
=1,\ldots, s, \
{}^k B_{\nu}^{} =
\bigl\{X_{1\nu}^{}(t,x),\ldots,X_{k\nu}^{}(t,x),
\;\!
{\rm div}_x^{}X^\nu(t,x)\bigr\}, \,
\nu=s+1,\ldots,m,
$
\vspace{0.75ex}
the system of identeties
(2.12)  means that: 
the functions of each set ${}^k B_j,\,
j =1,\ldots, m,$ are linearly dependent with respect to  
independent variable $t_\zeta$ on the domain $\Pi^{\;\!\prime}$
under any fixed values of 
independent variables  $t_\gamma,\, \gamma=1,\ldots, m,\, \gamma\ne \zeta,$
and  dependent variables  $x_i,\, i=1,\ldots, n;$
and 
the functions of each set ${}^k B_j,\,
j =1,\ldots, m,$ are linearly dependent with respect to  
dependent variable $x_p$ on the domain $\Pi^{\;\!\prime}$
under any fixed values of 
independent variables  $t_\gamma,\, \gamma=1,\ldots, m,$
and  dependent variables  $x_i,\, i=1,\ldots, n,\, i\ne p.$
It holds true  under each fixed index $\zeta =s+1,\ldots,m$
and under each fixed index $p =k+1,\ldots,n.$

Therefore the Wronskians of each set  ${}^k B_j,\, j=1,\ldots,m,$
with respect to 
independent variables  $t_\zeta,\, \zeta=s+1,\ldots, m,$
and  dependent variables  $x_p,\, p=k+1,\ldots, n$
vanish identically on the domain $\Pi^{\;\!\prime},$ that is, 
the system of identities holds:
\\[2.75ex]
\mbox{}           
$
{\rm W}_{t_\zeta^{}}^{}
\bigl(1,{}^k\;\!\!X^\theta(t,x),\;\!
{\rm div}_x X^\theta(t,x)\bigr) = 0
$
for all 
$(t,x) \in \Pi^{\prime},\ \
\theta =1,\ldots,s,\ \zeta =s+1,\ldots,m,
\hfill
$
\\[2ex]
\mbox{}
$
{\rm W}_{t_\zeta^{}}^{}
\bigl({}^k\;\!\!X^\nu(t,x),{\rm div}_xX^\nu(t,x)\bigr) = 0
$
for all 
$(t,x) \in \Pi^{\prime},\ \
\nu =s+1,\ldots,m,\ \zeta =s+1,\ldots,m,
\hfill
$
\\[0.5ex]
\mbox{}\hfill  (2.13)
\\[0ex]
\mbox{}
$
{\rm W}_{x_p}^{}
\bigl(1,{}^k\;\!\!X^\theta(t,x),\;\!
{\rm div}_xX^\theta(t,x)\bigr)= 0
$ for all 
$(t,x) \in \Pi^{\prime},\ \
\theta =1,\ldots,s, \ p =k+1,\ldots,n,
\hfill
$
\\[2.25ex]
\mbox{}
$
{\rm W}_{x_p}^{}
\bigl({}^k\;\!\!X^\nu(t,x),\;\!
{\rm div}_xX^\nu(t,x)\bigr) = 0
$ for all 
$(t,x) \in \Pi^{\prime},\ \
\nu =s+1,\ldots, m,\  p =k+1,\ldots,n,
\hfill
$
\\[3ex]
where the vector functions 
\vspace{0.3ex}
$
{}^k\;\!\!X^{j}\colon (t,x)\to \bigl(X_{1j}(t,x),\ldots,X_{kj}(t,x)\bigr)$
for all $(t,x) \in \Pi, \
 j=
\linebreak
=1,\ldots,m,$
${\rm W}_{t_\zeta}$ and ${\rm W}_{x_ p}$
are 
\vspace{0.3ex}
correspondingly the Wronskians with respect to $t_\zeta$ and $x_p,\, 
 \zeta =s+1,\ldots, m,\, p =k+1,\ldots, n.$
\vspace{0.3ex}

So the 
{\sl necessary test of  existence of 
$s\!$-nonautonomous
$(n - k)\!$-cylindrical last multiplier  for 
total differential system} is proved.
\vspace{0.2ex}
 
{\bf Theorem 2.4.}
{\it
The system of identities
{\rm (2.13)} is a necessary condition of  existence  
of  $s\!$-nonautonomous $(n - k)\!$-cylindrical last multiplier  {\rm (2.10)}  
for system} (TD).
\\[2ex]
\indent
{\bf
2.2.2. Criterion  of existence of 
{\boldmath$s\!$}-nonautonomous
{\boldmath$(n-k)\!$}-cylindrical last multiplier.}
Let the $n\times m$ matrix $X$ of  system (TD)  satisfies  the  conditions (2.13).
Let's write the functional system 
\\[2ex]
\mbox{}\hfill                    
$
\psi_\theta +  {}^k\;\!\!X^\theta(t,x)\ {}^k\varphi =
{}-{\rm div}_xX^\theta(t,x),
\ \ \theta =1,\ldots, s,
\hfill
$
\\[2ex]
\mbox{}\hfill
$
\partial_{{}_{\scriptstyle t_\zeta}}^\xi\;
{}^k\;\!\!X^\theta(t,x)\ {}^k\varphi =
{}- \partial_{{}_{\scriptstyle t_\zeta}}^\xi{\rm div}_x
X^\theta(t,x), \ \
\theta =1,\ldots, s,
\ \zeta =s+1,\ldots, m,
\ \xi=1,\ldots, k+1,
\hfill
$
\\[2ex]
\mbox{}\hfill
$
\partial_{{}_{\scriptstyle x_p}}^\xi\;
{}^k \;\!\!X^\theta(t,x)\ {}^k\varphi =
{}- \partial_{{}_{\scriptstyle x_p}}^\xi{\rm div}_xX^\theta(t,x),
\ \ 
\theta =1,\ldots, s,
\ p =k+1,\ldots, n,
\ \xi=1,\ldots, k+1,
\hfill
$
\\[2ex]
\mbox{}\hfill
$
 {}^k \;\!\!X^\nu(t,x)\ {}^k\varphi  =
{}- {\rm div}_xX^\nu(t,x),
 \ \ \nu =s+1,\ldots, m,
$
\hfill (2.14)
\\[2ex]
\mbox{}\hfill
$
\partial_{{}_{\scriptstyle t_\zeta}}^\xi\;{}^k\;\!\!
X^\nu(t,x)\ {}^k\varphi 
= {}- \partial_{{}_{\scriptstyle t_\zeta}}^\xi{\rm div}_xX^\nu(t,x),
\ \ 
\nu =s+1,\ldots, m,
\ \zeta =s+1,\ldots, m,
\ \xi=1,\ldots, k,
\hfill
$
\\[2ex]
\mbox{}\hfill
$
\partial_{{}_{\scriptstyle x_p}}^\xi\;{}^k \;\!\!
X^\nu(t,x) \ {}^k \varphi=
{}- \partial_{{}_{\scriptstyle x_p}}^\xi{\rm div}_xX^\nu(t,x),
\ \ 
\nu =s+1,\ldots, m,
\ p =k+1,\ldots, n,
\ \xi=1,\ldots, k,
\hfill
$
\\[3.5ex]
where the vector functions
\vspace{0.75ex}
$
{}^s \psi\colon (t,x)\to
\bigl( \psi_1^{}({}^st,{}^kx),
\ldots,
\psi_s^{}({}^st,{}^kx)\bigr)
$
for all $(t,x)\in \Pi^{\prime}
$
and 
$
{}^k \varphi\colon (t,x)\to
\bigl( \varphi_1^{}({}^st,{}^kx),
\ldots,
\varphi_k^{}({}^st,{}^kx)\bigr)
$
for all $(t,x)\in \Pi^{\prime}$
are unknown,
\vspace{1ex}
the vector functions
$X^{j}\colon (t,x)\to \bigl(X_{1j}(t,x),\ldots,X_{nj}(t,x)\bigr)$ for all $(t,x)\in \Pi,\   j=1,\ldots,m,$
\vspace{1ex}
and the vector functions
$
{}^k\;\!\!X^j\colon (t,x)\to
\bigl( X_{1j}(t,x),\ldots,X_{kj}(t,x)\bigr)
$
for all
$(t,x)\in \Pi, \
j =1,\ldots,m.
$
\vspace{1.5ex}

{\bf Theorem  1.8} ({\sl 
criterion  of existence of 
$s\!$-nonautonomous
$(n - k)\!$-cylindrical  last multiplier for 
total differential system}).
{\it
For  system  {\rm (TD)} 
to have 
$s\!$-nonautono\-mous
$(n - k)\!$-cylindrical  last multiplier {\rm(2.10)}
it is necessary and sufficient that there exist the vector functions
${}^s\psi$ and ${}^k\varphi,$ 
satisfying functional system {\rm(2.14)}, 
such that 
the Pfaffian equation  {\rm(1.21)} 
\vspace{0.25ex}
which is constructed on the base of this vector functions
is exact 
on the domain $\widetilde{\Pi}^{\;\!s+k}$
which is the natural projection of domain  $\Pi^{\;\!\prime}$ on 
coordinate subspace $O\ {}^st\  {}^kx.$
At that,  the last multiplier {\rm(2.10)} of 
system  {\rm (TD)}  is
}
\\[2ex]
\mbox{}\hfill                                          
$
\mu\colon (t,x)\to \exp g({}^s t,{}^k x)
$
\ for all 
$(t,x)\in \Pi^{\prime},
$
\hfill (2.15)
\\[1ex]
{\it where}
\\[2ex]
\mbox{}\hfill                                          
$
\displaystyle
g\colon ({}^s t,{}^k  x)\to 
\int {}^s \psi ({}^s t,{}^k x)\ 
d\ {}^s t + {}^k\varphi ({}^st,{}^k x)\ 
d\ {}^kx
$
\ for all 
$({}^s t,{}^k  x)\in \widetilde{\Pi}^{s+k}.
$
\hfill (2.16)
\\[3ex]
\indent
{\sl Proof. Necessity}. Let system  (TD)
has  the
$s\!$-nonautonomous
$(n - k)\!$-cylindrical  last multiplier (2.10)
on the domain  $\Pi^{\;\!\prime}.$ Then, the identities (2.12) are satisfied.
By means of termwise division of every identity (2.12) by 
$\mu({}^s t,{}^k x)$ we get a new system of identities
\\[2ex]
\mbox{}\hfill
$
\!\!
\displaystyle
\partial_{t_\theta^{}}^{}\!
\ln\mu ({}^s t,{}^k x ) +
\sum\limits_{\xi=1}^k\;\!
X_{\xi \theta}^{}(t,x)
\partial_{x_\xi^{}}^{}\!
\ln\mu ({}^s t,{}^k x) +\;\!
{\rm div}_x X^\theta(t,x) = 0
$
for\,all
$(t,x)\!\in\! \Pi^{\prime}_0,
 \
\theta\!=\!1,\ldots,s,
\hfill
$
\\[1.5ex]
\mbox{}\hfill
$
\displaystyle
\sum\limits_{\xi=1}^k\;\!
X_{\xi\nu}^{}(t,x)
\partial_{x_\xi^{}}^{}\!
\ln\mu ({}^s t,{}^k x ) +
\;\!{\rm div}_xX^\nu(t,x) = 0
$ \ for all
$(t,x)\in \Pi^{\prime}_0\subset \Pi^\prime,
\ \  \nu=s+1,\ldots,m.
\hfill
$
\\[2.5ex]
\indent
By differentiating the first $s$ of this identities  $k$ times with respect to 
$t_{s+1},\ldots,t_m$ and $k$ times with respect to $x_{k+1},\ldots,x_n$
and by differentiating the rest $m-s$  identities
$k-1$ times with respect to 
$t_{s+1},\ldots,t_m$ and $k-1$ times with respect to $x_{k+1},\ldots,x_n$
\vspace{0.3ex}
we conclude  that the 
extensions 
on the domain $\Pi_0^{\;\!\prime}$ of the functions

\newpage

\mbox{}
\\[-2.25ex]
\mbox{}\hfill                                    
$
{}^s \psi\colon  ({}^s t,{}^k x ) \to
\partial_{_{\scriptstyle  {}^s t}}
\ln\mu ({}^s t,{}^k x)
$
\ for all 
$({}^s t,{}^k x )\in\widetilde{\Pi}_0^{s+k}\subset \K^{s+k},
\hfill
$
\\
\mbox{}\hfill(2.17)
\\
\mbox{}\hfill
$
{}^k\varphi\colon  ({}^s t,{}^k x)\to
\partial_{{}^k x}^{}
\ln\mu ({}^s t,{}^k x)$
\ for all 
$({}^s t,{}^k x )\in\widetilde{\Pi}_0^{s+k},
\hfill
$
\\[2.75ex]
is a solution to the functional system (2.14).

The Pfaffian equation  {\rm(1.21)} 
\vspace{0.25ex}
which is constructed on the base of the functions (2.17)
is exact 
on the domain $\widetilde{\Pi}_0^{\;\!s+k}.$

From (2.17) it follows that 
$s\!$-nonautonomous
$(n - k)\!$-cylindrical  last multiplier  $\mu$ of system (TD)
is constructing on the domain $\Pi_0^{\;\!\prime}$
on the base of solutions to the functional system (2.14) by formula (2.15)
with (2.16).

By restriction
the domain $\Pi^{\;\!\prime}$ to its codomain 
$\Pi_0^{\;\!\prime}$ we conclude that the
necessary condition of Theorem 2.5 is satisfied.
\vspace{0.3ex}

{\sl Sufficiency}. Let the vector functions ${}^s\psi$ and ${}^k\varphi$
be a solution to the  functional system (2.14) and 
\vspace{0.3ex}
the Pfaffian equation  (1.21)
which is constructed on its base 
is exact 
on the domain $\widetilde{\Pi}_0^{\;\!s+k}\subset \K^{s+k}.$
Then, 
\\[1ex]
\mbox{}\hfill
$
\partial_{_{\scriptstyle  {}^s t}}\,
g ({}^s t,{}^k x ) =
{}^s\psi ({}^s t,{}^k x)
$ 
\ for all 
$({}^s t,{}^k x )\in\widetilde{\Pi}^{s+k},
\hfill
$
\\[2ex]
\mbox{}\hfill
$
\partial_{ {}^k x}^{}\;\!
g ({}^s t,{}^k x ) =
{}^k \varphi ({}^s t,{}^k x)$ 
\ for all 
$({}^s t,{}^k x )\in\widetilde{\Pi}^{s+k}.
\hfill
$
\\[2.5ex]
\indent
Taking into account that the vector functions 
${}^s\psi$ and  ${}^k\varphi$
are a solution to the functional system (2.14)
we receive that the system of identities (2.11)
is satisfied ralative to the function (2.15) with (2.16).

Therefore the function (2.15) with (2.16) is   
an $s\!$-nonautonomous
$(n - k)\!$-cylindrical  last multiplier  of   system (TD). \k
\\[2ex]
\indent
{\bf
2.2.3.    Functionally independent 
{\boldmath$s\!$}-nonautonomous
{\boldmath$(n-k)\!$}-cylindrical last multipliers.}
The method which is proposed in Theorem 2.5
can be used to construct  
the functionally independent 
$s\!$-nonautonomous $(n-k)\!$-cylindrical last multipliers
\vspace{0.5ex}
of system (TD).

{\bf Theorem  2.6.}
{\it
Let the functional system  {\rm(2.14)} has 
$q$ 
not  linearly bound 
on the domain $\Pi^{\;\!\prime}$  solutions {\rm (1.25)}
and the corresponding  Pfaffian equations {\rm (1.26)}
are exact 
on the domain $\widetilde{\Pi}^{\;\!s+k}$
which is the natural projection of domain  $\Pi^{\;\!\prime}$ on 
coordinate subspace $O\ {}^st\ {}^kx.$
Then,  the 
$s\!$-nonautonomous $(n-k)\!$-cylindrical last multipliers
of  system  {\rm (TD)} 
}
\\[2ex]
\mbox{}\hfill
$
\displaystyle
\mu_\gamma\colon (t,x)\to
\exp\int {}^s \psi^\gamma ({}^s t,{}^k x )\  d\ {}^s t  \,+\,
{}^k  \varphi^\gamma ({}^s t,{}^k x )\  d\ {}^k x
$
\ for all
$(t,x)\in \Pi^{\prime},
 \ \
\gamma=1,\ldots,q,
\hfill
$
\\[2.25ex]
{\it 
are functionally independent
on the domain  $\Pi^{\;\!\prime}.$
}
\vspace{0.2ex}

{\sl Proof.} 
From Theorem 2.5 it follows that the last multipliers
$\mu_\gamma^{},\, \gamma=1,\ldots,q,$ of system  (TD)  
are of indicated structure.

From representations
\\[2.25ex]
\mbox{}\hfill
$
\partial_{{}_{\scriptstyle t_\theta^{}}}\!
\ln\mu_\gamma ({}^s t,{}^k x ) =
\psi_{{}_{\scriptstyle \theta}}^{\gamma} ({}^s t,{}^k x )
$
\ for all
$({}^s t,{}^k x)\in \widetilde{\Pi}^{s+k},
\ \ \theta = 1,\ldots, s, \ 
\gamma =1,\ldots, q,
\hfill
$
\\[2ex]
\mbox{}\hfill
$
\partial_{x_\xi^{}}^{}\!
\ln\mu_\gamma ({}^s t,{}^k x) =
\varphi_{{}_{\scriptstyle \xi}}^{\gamma} ({}^s t,{}^k x)
$
\ for all
$({}^s t,{}^k x)\in \widetilde{\Pi}^{s+k},
\ \ \xi = 1,\ldots, k, \ 
\gamma =1,\ldots, q,
\hfill
$
\\[2.5ex]
it follows that the  Jacobi's matrix
\\[2.25ex]
\mbox{}\hfill
$
J\bigl(\ln\mu_\gamma({}^s t,{}^k x); {}^s t,{}^k x \bigr) =
\bigl\| \Psi({}^s t,{}^k x)\ \Phi({}^s t,{}^k x)\bigr\|_{q\times(s+k)}
$
\ for all
$({}^s t,{}^k x)\in \widetilde{\Pi}^{s+k},
\hfill
$
\\[2.75ex]
where the matrix  $\|\Psi\Phi\|$ consists of 
\vspace{0.75ex}
$q\times s$ matrix
$
\Psi({}^st,{}^kx) =
\bigl\|\psi_{_{\scriptstyle \theta}}^{\gamma}({}^s t,{}^k x)\bigr\|$
for all $({}^st,{}^kx) \in 
\linebreak
\in\widetilde{\Pi}^{s+k}$
and  $q\times k$ matrix
$
\Phi({}^st,{}^kx) =
\bigl\|\varphi_{_{\scriptstyle \xi}}^{\gamma}({}^s t,{}^k x)\bigr\|$
\vspace{0.75ex}
for all $({}^st,{}^kx) \in \widetilde{\Pi}^{s+k}.$

Since  the solutions (1.25) 
\vspace{0.3ex}
to the functional system  (2.14)
are not  linearly bound 
on the domain $\Pi^{\;\!\prime}$
the rank of Jacobi's matrix
$
{\rm rank}\, J
\bigl(\ln\mu_\gamma({}^st,{}^kx); {}^st,{}^kx\bigr) = q
$
nearly everywhere on the domain  $\widetilde{\Pi}^{\;\!s+k}.$
Therefore  the 
$s\!$-nonautonomous $(n-k)\!$-cylindrical last multipliers
$\mu_\gamma,\,\gamma=1,\ldots,q,$
of system  (TD) 
are functionally independent
on the domain  $\Pi^{\;\!\prime}.$  \k
\\[4ex]
\centerline{
\large\bf  
3. Cylindricality and autonomy of partial integrals
}
\\[2.25ex]
\centerline{
\bf 
3.1. Cylindricality of partial integrals for  linear homogeneous   system 
}
\\[0.2ex]
\centerline{
\bf 
of partial differential equations
}
\\[1.25ex]
\indent
{\bf Definition  3.1.}
{\it
We'll say  that a  partial integral  $w$
on a domain $G^{\;\!\prime}\subset G$
of  system  {\rm($\!\partial\!$)}
is {\boldmath$(n - k)\!$}\textit{\textbf{-cylindrical}}
if the function  $w$ depends only on 
\vspace{0.5ex}
 $k,\, 0\leq k\leq n,$ variables  $x_1,\ldots,x_n.$
}

Let's define the problem of  existence
for  system ($\!\partial\!$) an  $(n-k)\!$-cylindrical partial integral
\\[2ex]
\mbox{}\hfill                                           
$
w\colon x \to w({}^kx)
$
\ for all 
$x\in G^{\;\!\prime}\subset G,
\quad 
{}^k x=(x_1,\ldots,x_k).
$
\hfill (3.1)
\\[2.75ex]
\indent
{\bf
3.1.1.
Necessary condition of existence of  cylindrical
partial integral.}
According to the definition of 
partial integral,  the function  (3.1) 
will be the partial integral on the domain $G^{\;\!\prime}\subset G$ of  system (\!$\partial$\!)
if and only if 
\\[2ex]
\mbox{}\hfill                                           
$
{}^k {\frak L}_j^{}  w({}^kx) = \Phi_j(x)
$
\ for all 
$x\in G^{\;\!\prime}, \ \ \
j =1,\ldots, m,
$
\hfill (3.2)
\\[2.2ex]
where the linear differential operators of first order
\vspace{0.25ex}
${}^k {\frak L}_j,\,j = 1,\ldots,m,$
are defined by means of  (1.3),
the
scalar  functions $\Phi_j\colon  G^{\;\!\prime}\to \K,\, j =1,\ldots, m,$
are such that 
\\[2ex]
\mbox{}\hfill                                           
$
\Phi_j(x)_{\displaystyle |_{\scriptstyle w({}^k x)=0}}
=0
$
\ for all 
$x\in G^{\;\!\prime}, \ \ \
j =1,\ldots, m.
$
\hfill (3.3)
\\[2.5ex]
\indent
The system of identities (3.2) in the coordinates  is given by 
\\[2ex]
\mbox{}\hfill                                           
$
\displaystyle
\sum\limits_{\xi =1}^k\,
u_{j\xi}^{} (x)\;\! \partial_{x_\xi^{}}^{}w({}^k x)=
\Phi_j(x)
$
\ for all 
$x\in G^{\;\!\prime}, \ \ \
j =1,\ldots, m.
$
\hfill (3.4)
\\[2.25ex]
\indent
Concerning the sets of functions 
\vspace{0.3ex}
$
{}^k U_j =
\bigl\{ u_{j1}^{}(x),\ldots, u_{jk}^{}(x)\bigl\}, \,
j =1,\ldots, m,$
the system of identeties
(3.4) with (3.3)  means that the functions of each set 
\vspace{0.3ex}
${}^k U_j,\,
j =1,\ldots, m,$ are linearly dependent with respect to  
variable $x_p$ on the integral manifold  
$w({}^k x)=0$
under any fixed values of 
variables  $x_i,\,
i=1,\ldots, n,\, i\ne p.$ 
It holds true  under each fixed index $p =k+1,\ldots,n.$
Therefore the Wronskians of each set  ${}^k U_j,\, j=1,\ldots,m,$
with respect to variables $x_p ,\, p =k+1,\ldots,n,$ 
vanish identically on the integral manifold  
$w({}^k x)=0,$ that is,  the system of identities holds:
\\[2.25ex]
\mbox{}\hfill                                 
$
{\rm W}_{x_p}^{}
\bigl({}^k u^{j} (x)\bigr) = \Xi_{jp}^{}(x)
$
\ for all 
$x\in G^{\;\!\prime},
\ \ \ j=1,\ldots, m,\ p =k+1,\ldots, n,
$
\hfill (3.5)
\\[2.25ex]
where the vector functions 
\vspace{0.3ex}
$
{}^k u^{j}\colon x\to
\bigl(u_{j1}^{}(x),\ldots,u_{jk}^{}(x)\bigr)$
for all $x\in G, \ j=1,\ldots,m,$
${\rm W}_{x_p}$ are the Wronskians with respect to $x_p,\, p =k+1,\ldots, n,$
\vspace{0.3ex}
the scalar functions $\Xi_{jp}^{}\colon G^{\;\!\prime}\to\K,$
$ j=1,\ldots, m,\, p =k+1,\ldots, n,$
are such that
\\[2ex]
\mbox{}\hfill                                           
$
\Xi_{jp}(x)_{\displaystyle |_{\scriptstyle w({}^k x)=0}}
=0
$
\ for all 
$x\in G^{\;\!\prime}, \ \ \
j =1,\ldots, m, 
\  p =k+1,\ldots, n.
$
\hfill (3.6)
\\[2.5ex]
\indent
So the 
{\sl necessary test of  existence of $(n - k)\!$-cylindrical 
partial integral  for 
linear homogeneous   system of partial differential equations}
 is proved.
 
{\bf Theorem 3.1.}
{\it
The system of identities
{\rm (3.5)} with {\rm (3.6)} is a necessary condition of  existence  
of $(n - k)\!$-cylindrical partial integral {\rm (3.1)}  
for system} (\!$\partial$\!).
\\[1.75ex]
\indent
{\bf
3.1.2. Criterion  of existence of cylindrical partial integral.}
Let the $m\times n$ matrix  $u(x)=\|u_{ji}^{}(x)\|$ for all $x\in G$ 
of  system (\!$\partial$\!)  satisfies  the  conditions (3.5)$\cup$(3.6).
Let's write the functional system 
\\[1ex]
\mbox{}\hfill                            
$
{}^k u^j(x) \, {}^k\varphi= H_j(x), \ 
\ j=1,\ldots,m,
\hfill
$
\\[-0.75ex]
\mbox{}\hfill (3.7)
\\
\mbox{}\hfill
$
\partial^{\xi}_{{}_{\scriptstyle x_p}}
\;{}^k u^j(x)\, {}^k\varphi = 
\partial^{\xi}_{{}_{\scriptstyle x_p}}
H_j(x),
\ \ 
j=1,\ldots,m,
\ \xi =1,\ldots,k-1, \
p =k+1,\ldots,n, 
\hfill
$
\\[2.5ex]
where a vector function
\vspace{0.5ex}
$
{}^k\varphi\colon x\to
\bigl(\varphi_1({}^kx),\ldots,\varphi_k({}^kx)\bigr)$
for all $x\in G^{\;\!\prime}$
is unknown,
the vector functions
$
{}^k u^{j}\colon x\to
\bigl(u_{j1}^{}(x),\ldots,u_{jk}^{}(x)\bigr)$ for all $x\in G, \ j=1,\ldots,m,$
\vspace{0.3ex}
the scalar functions $H_{j}\colon G^{\;\!\prime}\to\K,\,  j=1,\ldots, m,$
are such that
\\[2ex]
\mbox{}\hfill                                           
$
H_j(x)_{\displaystyle |_{\scriptstyle w({}^k x)=0}}
=0
$
\ for all 
$x\in G^{\;\!\prime}, \ \ \
j =1,\ldots, m.
$
\hfill (3.8)
\\[2.75ex]
\indent
{\bf Theorem  3.2} ({\sl 
criterion  of existence of $(n - k)\!$-cylindrical partial integral
for linear homogeneous   system of partial differential equations}).
{\it
For  system  {\rm (\!$\partial$\!)} 
to have $(n - k)\!$-cylindrical partial integral
{\rm(3.1)}
it is necessary and sufficient that there exists a vector function
${}^k\varphi$
and  scalar functions $H_j,\,j =1,\ldots, m,$ with {\rm(3.8),}
satisfying functional system {\rm(3.7)}, 
\vspace{0.2ex}
such that the Pfaffian equation  {\rm(1.6)} 
has the general integral
$w\colon {}^k x \to w({}^k x)$ for all ${}^kx \in \widetilde{G}^{\;\!k},$
\vspace{0.3ex}
where domain $\widetilde{G}^{\;\!k}$
is the natural projection of domain  $G^{\;\!\prime}$ on 
coordinate subspace $O\ {}^kx.$
}
\vspace{0.3ex}

{\sl Proof. Necessity}. Let system  (\!$\partial$\!)
has  the $(n - k)\!$-cylindrical  partial integral (3.1)
on the domain  $G^{\;\!\prime}.$ Then, the system of identities (3.4) with (3.3) is satisfied.
By differentiating this identities  $k-1$ times with respect to $x_{p},\, p=k+1,\ldots,n,$  
we conclude  that an 
extension 
on the domain $G^{\;\!\prime}$ of the  function
\\[2ex]
\mbox{}\hfill
$
{}^k\varphi\colon {}^kx \to
\bigl(
\partial_{x_1}^{} w({}^kx),\ldots, \partial_{x_k}^{}w({}^kx)\bigr)
$
\ for all  ${}^k x\in \widetilde{G}^{\;\!k}
\hfill
$
\\[2.25ex]
is a solution to the functional system (3.7) with (3.8).
\vspace{0.3ex}
From this  it also follows that 
the function (3.1) is the general integral
on the domain $\widetilde{G}^{\;\!k}\subset \K^k$
of the Pfaffian equation  (1.6). 
\vspace{0.3ex}

{\sl Sufficiency}. Let  the vector function
$
{}^k\varphi\colon x\to {}^k\varphi({}^kx)$ for all $x\in G^{\;\!\prime}$
be a solution to the functional system (3.7) with (3.8)
and the Pfaffian equation  (1.6)
which is constructed on its base 
has the general integral
$w\colon {}^k x \to w({}^k x)$ for all ${}^kx \in \widetilde{G}^{\;\!k}.$
Then, the system of  identities 
\\[2.25ex]
\mbox{}\hfill                                     
$
\partial_{x_\xi^{}}^{} w({}^kx) -
\mu({}^k x) \varphi_{\xi}^{}({}^kx) = 0
$
\ for all 
${}^k x\in \widetilde{G}^{\;\!k}, \ \ \
\xi =1,\ldots, k,
$
\hfill (3.9)
\\[2.5ex]
is  satisfied. Here
$\mu\colon {}^k x \to \mu({}^k x)$ for all ${}^kx \in \widetilde{G}^{\;\!k}$ 
 is  the  holomorphic integrating multiplier of the Pfaffian equation  (1.6) 
which corresponds to its general integral 
$w\colon {}^k x \to w({}^k x)$ for all ${}^kx \in \widetilde{G}^{\;\!k}.$

Taking into account that the vector function ${}^k \varphi$ is
the solution to the functional system (3.7) with (3.8)
we receive the system of identities (3.4), where 
\\[2ex]
\mbox{}\hfill                                    
$
\Phi_j(x)=\mu({}^k x) H_j(x)
$
\ for all  
$x\in G^{\;\!\prime}, \ \ \
j =1,\ldots, m.
\hfill 
$
\\[2.5ex]
\indent
Therefore  the function (3.1) is   
an $(n - k)\!$-cylindrical partial integral   
on the domain  $G^{\;\!\prime}$ of   system (\!$\partial$\!). \k
\vspace{0.75ex}

{\bf Example 3.1.}
Consider the
linear homogeneous   system of partial differential equations
\\[2ex]
\mbox{}\hfill            
$
{\frak L}_1(x)y=0,
\ \
{\frak L}_2(x)y=0,
$
\hfill (3.10)
\\[2ex]
where the linear differential operators of first order
\\[2ex]
\mbox{}\hfill            
$
{\frak L}_1(x)= 
x_1(x_2 + x_3)\partial_{x_1}^{}  +
x_2(x_2 + x_3)\partial_{x_2}^{} + 
(x_1^2 + x_2^2 + x_3^2 + x_4^2)\partial_{x_3}^{}+ 
(x_1^2 - x_2^2 + x_3^2 - x_4^2)\partial_{x_4}^{}
\hfill            
$
\\[1.5ex]
\mbox{}\hfill
$            
\text{for all} \ \ x\in \K^4, 
\hfill
$ 
\\[2.5ex]
\mbox{}\hfill            
$
{\frak L}_2 = x_1(x_3 + x_4)\partial_{x_1}^{}+ 
x_2(x_3 + x_4)\partial_{x_2}^{}+ 
(x_1^2 - x_2^2 + x_3^2 - x_4^2)\partial_{x_3}^{} +(x_1^2 + x_2^2 + x_3^2 + x_4^2)\partial_{x_4}^{}
\hfill            
$
\\[1.5ex]
\mbox{}\hfill            
$            
\text{for all} \ \ x\in \K^4. 
\hfill
$ 
\\[2.5ex]
\indent
Let's find for system (3.10)
a 2-cylindrical partial  integral
\\[2ex]
\mbox{}\hfill            
$
w\colon x \to w(x_1,x_2)
$
\ for all  
$x\in G^{\;\!\prime}\subset  \K^4.
$
\hfill (3.11)
\\[2.5ex]
\indent
The Wronskians of the sets of functions
\vspace{0.3ex}
${}^2U_1=\{ x_1(x_2 + x_3), x_2(x_2 + x_3)\}$
and 
${}^2U_2=
\linebreak
=\{ x_1(x_3+ x_4), x_2(x_3 + x_4)\}$
with respect to $x_3$ and  $x_4$  
vanish  identically on the space $\K^4\colon$ 
\\[2.25ex]
\mbox{}\hfill            
$
{\rm W}_{x_3}^{} \bigl(x_1(x_2 + x_3), x_2(x_2 + x_3)\bigr) =
\left|\!
\begin{array}{cc}
x_1(x_2 + x_3) & x_2(x_2 + x_3)
\\[0.5ex]
x_1 & x_2
\end{array}
\!\right| = 0 
$ \ for all 
$x \in \K^4,
\hfill            
$
\\[2.5ex]
\mbox{}\hfill            
$
{\rm W}_{x_4}^{}\bigl( x_1(x_2 + x_3),
x_2(x_2 + x_3)\bigr) =  0
$ 
\ for all 
$x \in \K^4,
\hfill            
$
\\[2.5ex]
\mbox{}\hfill            
$
{\rm W}_{x_3}^{}\bigl(x_1(x_3 + x_4),
 x_2(x_3 + x_4)\bigr) =
\left|\!
\begin{array}{cc}
x_1(x_3 + x_4) & x_2(x_3 + x_4)
\\[0.5ex]
x_1 & x_2
\end{array}
\!\right| = 0
$ \ for all 
$x \in \K^4,
\hfill            
$
\\[2.5ex]
\mbox{}\hfill            
$
{\rm W}_{x_4}^{}\bigl( x_1(x_3 + x_4),
x_2(x_3 + x_4)\bigr) = 
\left|\!
\begin{array}{cc}
x_1(x_3 + x_4) & x_2(x_3 + x_4)
\\[0.5ex]
x_1 & x_2
\end{array}
\!\right| = 0
$
\ for all 
$x \in \K^4.
\hfill            
$
\\[2.75ex]
\indent
Therefore the  necessary conditions (Theorem 3.1)
of  existence of 2-cylindrical partial integral (3.11)  for 
system (3.10) are satisfied.

Let's write the functional system  (3.7) with (3.8):
\\[2ex]
\mbox{}\hfill            
$
x_1(x_2 + x_3)\,\varphi_1  + x_2(x_2 + x_3)\,\varphi_2  =  
(x_1 + x_2)(x_2 + x_3),
\quad
x_1\,\varphi_1 + x_2\,\varphi_2 = x_1 + x_2,
\hfill            
$
\\[2ex]
\mbox{}\hfill            
$
x_1(x_3 + x_4)\,\varphi_1  + 
x_2(x_3 + x_4)\,\varphi_2  =  (x_1 + x_2)(x_3 + x_4),
\quad
x_1\,\varphi_1  + x_2\,\varphi_2  = x_1 + x_2,
\hfill
$
\\[2.75ex]
where $H_1(x)  =  (x_1 + x_2)(x_2 + x_3)$ for all  $x\in \K^4, \
\vspace{0.5ex}
H_2(x) = (x_1 + x_2)(x_3 + x_4)$ for all  $x\in\K^4.$

On the base of solution
$\varphi_1\colon x \to 1$ for all $x \in \K^4, \
\varphi_2\colon x \to 1$ for all $x \in \K^4$
to this system we construct the Pfaffian equation
\\[2ex]
\mbox{}\hfill            
$
dx_1 + dx_2 = 0
\hfill
$
\\[2ex]
which is exact
(the integrating multiplier
$\mu\colon (x_1,x_2) \to 1$ for all $(x_1,x_2) \in \K^2$\!)
on the plane $\K^2$
and has the general integral 
$w\colon (x_1,x_2) \to x_1  +  x_2$ for all $(x_1,x_2) \in \K^2.$
\vspace{0.3ex}

By extension 
of the general integral on the space 
$\K^4$ we get the 
2-cylindrical partial  integral
$w\colon x \to x_1  +  x_2$ for all $x \in \K^4$
of system (3.10).
\\[2ex] 
\indent
{\bf 3.1.3.    Functionally independent cylindrical 
partial  integral.}
The method which is proposed in Theorem 3.2
can be used to construct  
the functionally independent $(n-k)\!$-cylindrical
partial  integrals of system (\!$\partial$\!).
\vspace{0.3ex}

{\bf Theorem  3.3.}
{\it
Let $h$ functional systems  {\rm(3.7)} with {\rm(3.8)} has 
$q$ 
not  linearly bound 
on the domain $G^{\;\!\prime}\subset G$  solutions
{\rm (1.11)}
and for each of them the corresponding  Pfaffian equation {\rm (1.12)}
has the general integral 
}
\\[2ex]
\mbox{}\hfill                                           
$
w_{\gamma}^{}\colon
{}^kx\to w_{\gamma}^{}({}^kx)
$
\ for all  
${}^kx \in \widetilde{G}^{\;\!k}\subset \K^k,
\ \ \ \gamma =1,\ldots, q,
$
\hfill (3.12)
\\[2ex]
{\it
on the domain $\widetilde{G}^{\;\!k}$
which is the natural projection of domain  $G^{\;\!\prime}$ on 
coordinate subspace $O\ {}^kx.$
Then,  the general integrals {\rm (3.12)}
are functionally independent
on the domain $\widetilde{G}^{\;\!k}.$
}
\vspace{0.2ex}

{\sl Proof.} 
In accordance with the system of identities (3.9) we have
\\[2.2ex]
\mbox{}\hfill
$
\partial_{x_\xi^{}}^{} w_\gamma^{}({}^kx) -
\mu_\gamma^{}({}^kx) \varphi_{{}_{\scriptstyle \xi }}^\gamma({}^kx)=0
$
\ for all  
${}^kx \in \widetilde{G}^{\;\!k},
\ \ \ \xi =1,\ldots, k,
 \ \gamma =1,\ldots, q.
\hfill
$
\\[2.5ex]
Therefore the Jacobi's matrix
$
J\bigl(w_{\gamma}^{}({}^kx);\, {}^kx\bigr)=
\bigl\|	
\mu_\gamma^{}({}^kx) \varphi_{{}_{\scriptstyle \xi }}^\gamma({}^kx)\bigr\|_{q\times k}
$ for all ${}^kx \in \widetilde{G}^{\;\!k}.$
\vspace{0.3ex}

Since the vector functions (1.11) 
are not  linearly bound 
on the domain  $\widetilde{G}^{\;\!k}$
the rank of Jacobi's matrix
$
{\rm rank}\, J\bigl(w_{\gamma}^{}( {}^kx);\,{}^kx\bigr) = q
$
nearly everywhere on the domain  $\widetilde{G}^{\;\!k}.$
So  the general integrals {\rm (3.12)}
of the Pfaffian equation {\rm (1.12)}
are functionally independent
on the domain $\widetilde{G}^{\;\!k}.$  \k
\vspace{0.5ex}

{\bf Example 3.2.}
The 
linear homogeneous   system of partial differential equations
\\[2ex]
\mbox{}\hfill            
$
{\frak L}_1(x)y=0,
\ \
{\frak L}_2(x)y=0,
$
\hfill (3.13)
\\[2ex]
which is constructed on the base of linear differential operators of first order
\\[2ex]
\mbox{}\hfill
$
\displaystyle
{\frak L}_1(x)=
\sum\limits_{i=1}^5\, x_i\;\! \partial_{x_i}^{}
$
\ for all 
$x\in \K^5,
\qquad
{\frak L}_2(x)=
\sum\limits_{\nu=1}^3\, x_{\nu}\;\!\partial_{x_\nu}^{} +  x^2_4\,
\partial_{x_4}^{} + x^2_5\, \partial_{x_5}^{} 
$
\ for all 
$x\in \K^5,
\hfill
$
\\[2.5ex]
has the 4-cylindrical partial integrals
\\[2ex]
\mbox{}\hfill
$
w_\nu\colon x\to x_\nu
$
\ for all 
$x\in \K^5,
\ \ \ \nu =1,2,3,
\hfill
$
\\[2.5ex]
as  ${\frak L}_j x_\nu=x_\nu$ for all $x\in \K^5, \ j=1,2,\, \nu=1,2,3.$
\vspace{0.3ex}

Let's construct a basis of first integrals for  system (3.13)
on the base of this 4-cylindrical partial integrals.

The system (3.13) is incomplete
and can be reduced to the complete system  
by the addition of single  operator 
\\[2ex]
\mbox{}\hfill
$
{\frak L}_{12} (x)=
\left[ {\frak L}_1 (x),{\frak L}_2 (x)\right]=
x^2_4\, \partial_{x_4}^{} + x^2_5\,\partial_{x_5}^{}
$
\ for all 
$x\in\K^5.
\hfill
$
\\[2.5ex]
\indent
Therefore the incomplete system  (3.13)
has the defect  $\delta=1$
and its integral basis consists of  $n-m-\delta=5-2-1=2$
functionally independent first integrals.

Let's reduce the system ${\frak L}_1(x)y=0,\,{\frak L}_2(x)y=0,\,{\frak L}_{12} (x)y=0$
to the complete normal system
\\[0.5ex]
\mbox{}\hfill
$
\partial_{x_1}^{} y = {}-x_2 x_1^{{}-1}\,\partial_{x_2}^{} y
-x_3x_1^{{}-1}\,\partial_{x_3}^{} y,
\quad
\partial_{x_4}^{} y = 0,
\quad
\partial_{x_5}^{} y = 0
\hfill
$
\\[3ex]
on a domain $H_1\subset \{ x\colon  x_1 \ne 0\}.$ 

From this we find an integral basis 
on the domain $H_1$ of system (3.13),
which consists of  two 
functionally independent 
3-cylindrical first integrals
\\[2ex]
\mbox{}\hfill
$
F_{12}\colon x\to
x_2\;\! x_1^{{}-1}
$ 
\ for all 
$x\in H_1$
\quad and \quad
$
F_{13}\colon x\to
x_3\;\!x_1^{{}-1}
$ 
\ for all 
$x\in H_1.
\hfill
$
\\[2.5ex]
\indent
Similarly, 
\vspace{0.3ex}
the system ${\frak L}_1(x)y=0,\,{\frak L}_2(x)y=0,\,{\frak L}_{12} (x)y=0$
is normalized on the domains $H_{\xi}\subset \{x\colon  x_{\xi}^{} \ne 0\},\,\xi=2,3,$
\vspace{0.3ex}
and the corresponding integral basises
on the domains $H_2$ and $H_3$ of system (3.13)
consist of functionally independent 
3-cylindrical first integrals
\\[2.25ex]
\mbox{}\hfill
$
F_{\xi\nu}^{}\colon x\to x_{\nu}\,x_{\xi}^{{}-1}
$
\ for all 
$x\in H_{\xi},
\qquad
F_{\xi\theta}^{}\colon x\to x_{\theta}^{}\,x_{\xi}^{{}-1}$ 
\ for all $x\in H_{\xi},
\hfill
$
\\[2ex]
\mbox{}\hfill
$
\xi=2,3, \ \ \nu =1,2,3,\ \ \nu\ne \xi, \ \  
\theta=1,2,3,  \ \ \theta\ne \xi, \ \  \nu\ne \theta.
\hfill
$
\\[4.5ex]
\centerline{
\bf
3.2.
Autonomy  and  cylindricality of partial integrals for
}
\\[0ex]
\centerline{
\bf
total differential system 
}
\\[1.5ex]
\indent
{\bf Definition 3.2.}
{\it
We'll say that a partial integral  $w$
on a domain  $\Pi^\prime\subset \Pi$ of  
system {\rm(TD)} is
{\boldmath$s\!$}\textit{\textbf{-nonautonomous}}
if the  function   $w$  depends on  $x$ and only on  $s,\,0\leq s\leq m,$ 
independent variables $t_1,\ldots,t_m.$
If $s=0,$ then  a partial integral
$w\colon (t,x) \to w(x)$ for all $(t,x) \in \Pi^{\prime}$
of  system {\rm (TD)} is \textit{\textbf{autonomous}}.
}
\vspace{0.5ex}

{\bf Definition  3.3.}
{\it
We'll say  that a partial integral  $w$
on a domain $\Pi^{\;\!\prime}\subset \Pi$
of  system  {\rm(TD)}
is {\boldmath$(n - k)\!$}\textit{\textbf{-cylindrical}}
if the function  $w$ depends on $t$ and only on 
$k,\, 0\leq k\leq n,$ dependent variables  $x_1,\ldots,x_n.$
}
\vspace{0.5ex}

Let's define the problem of  existence
for  system (TD) an  
$s\!$-nonautonomous
$(n-k)\!$-cylindrical partial integral
\\[2ex]
\mbox{}\hfill                                           
$
w\colon (t,x) \to w({}^st,{}^kx)
$ \ for all $(t,x)\in \Pi^{\prime}\subset \Pi,
\ \ \ 
{}^s t=(t_1,\ldots,t_s),\ \ {}^k x=(x_1,\ldots,x_k).
$
\hfill (3.14)
\\[3.5ex]
\indent
{\bf
3.2.1.
Necessary condition of existence of 
{\boldmath$s\!$}-nonautonomous
{\boldmath$(n-k)\!$}-cylindrical partial integral.}
According to the definition of partial integral, the function  (3.14) 
will be the partial integral on the domain $\Pi^{\;\!\prime}\subset \Pi$ of  system (TD)
if and only if 
\\[2ex]
\mbox{}\hfill                                           
$
{}^{sk}{\frak X}_j^{}w({}^st,{}^kx) = \Phi_j^{}(t,x)
$ \ for all
$(t,x)\in \Pi^{\prime}, \ \ \
j =1,\ldots,m,
$
\hfill (3.15)
\\[2.25ex]
where the 
linear differential operators of first order
${}^{sk}{\frak X}_j^{},\, j =1,\ldots,m,$ are defined by means of (1.18),
the scalar functions $\Phi_j\colon \Pi^\prime\to \K,\,j =1,\ldots,m,$ are such that
\\[2ex]
\mbox{}\hfill                                           
$
\Phi_j^{}(t,x)_{\displaystyle |_{\scriptstyle 
w({}^st,{}^kx)=0}}=0
$ \ for all 
$(t,x)\in \Pi^{\prime},
 \ \ \
j =1,\ldots,m.
$
\hfill (3.16)
\\[2.5ex]
\indent
The system of identities (3.15) in the coordinates  is given by 
\\[2ex]
\mbox{}\hfill      
$
\displaystyle
\partial_{t_\theta^{}}^{} w({}^st,{}^kx) +
\sum\limits_{\xi =1}^k\,
X_{\xi \theta}^{}(t,x)
\partial_{x_\xi^{}}^{}w({}^st,{}^kx) = \Phi_{\theta}^{}(t,x)
$
\ for all 
$(t,x)\in \Pi^\prime, \ \ \ \theta=1,\ldots,s,
\hfill
$
\\
\mbox{}\hfill (3.17)
\\[-1ex]
\mbox{}\hfill
$
\displaystyle
\sum\limits_{\xi=1}^k\,
X_{\xi \nu}^{}(t,x)
\partial_{x_\xi^{}}^{}w({}^st,{}^kx) = 
 \Phi_{\nu}(t,x)
$
\ for all 
$(t,x)\in \Pi^\prime, \ \ \ \nu=s+1,\ldots,m.
\hfill
$
\\[2.65ex]
\indent
Concerning the sets of functions 
\vspace{0.5ex}
$
{}^k M_\theta^{} =
\bigl\{1,X_{1\theta}^{}(t,x),\ldots,
X_{k\theta}^{}(t,x)\bigl\}, \,
\theta =1,\ldots, s, \
{}^k M_{\nu}^{} =
\linebreak
=
\bigl\{X_{1\nu}^{}(t,x),\ldots,X_{k\nu}^{}(t,x)\bigr\}, \,
\nu=s+1,\ldots,m,
$
\vspace{0.5ex}
the system of identeties
(3.17) with (3.16)  means that: 
the functions of each set ${}^k M_j,\,
j =1,\ldots, m,$ 
\vspace{0.25ex}
are linearly dependent with respect to  
independent variable $t_\zeta$ 
\vspace{0.25ex}
on the integral manifold $w({}^st,{}^kx)=0$ 
under any fixed values of 
independent variables  $t_\gamma,\, \gamma=1,\ldots, m,\, \gamma\ne \zeta,$
\vspace{0.5ex}
and  dependent variables  $x_i,\, i=1,\ldots, n;$
and the functions of each set ${}^k M_j,\,
j =1,\ldots, m,$ 
\vspace{0.3ex}
are linearly dependent with respect to  
dependent variable $x_p$ on the integral manifold $w({}^st,{}^kx)=0$ 
\vspace{0.25ex}
under any fixed values of 
independent variables  $t_\gamma,\, \gamma=1,\ldots, m,$
and  dependent variables  $x_i,\, i=1,\ldots, n,\, i\ne p.$
\vspace{0.25ex}
It holds true  under each fixed index $\zeta =s+1,\ldots,m$
and under each fixed index $p =k+1,\ldots,n.$
\vspace{0.25ex}

Therefore the Wronskians of each set  ${}^k M_j,\, j=1,\ldots,m,$
\vspace{0.2ex}
with respect to 
independent variables  $t_\zeta,\, \zeta=s+1,\ldots, m,$
and  dependent variables  $x_p,\, p=k+1,\ldots, n$
\vspace{0.3ex}
vanish identically on the integral manifold $w({}^st,{}^kx)=0,$ 
that is, 
the system of identities
\\[2.5ex]
\mbox{}\quad                                 
$
{\rm W}_{ t_\zeta}^{}
\bigl(1,{}^k X^{\theta}(t,x)\bigr) = 
\overset{*}\Xi_{\theta\zeta}^{}(t,x)
$
\ for all
$(t,x)\in \Pi^\prime, \ \ \
\theta =1,\ldots, s,\ \zeta =s+1,\ldots,m,
\hfill
$
\\[1.85ex]
\mbox{}\quad
$
{\rm W}_{t_\zeta}^{}
\bigl({}^k X^{\nu}(t,x)\bigr) = 
\overset{*}\Xi_{\nu\zeta}^{}(t,x)
$
\ for all
$(t,x)\in \Pi^\prime, \ \ \
\nu=s+1,\ldots, m,
\ \zeta =s+1,\ldots, m,
\hfill
$
\\[0.5ex]
\mbox{}\hfill (3.18)
\\[-0.75ex]
\mbox{}\quad
$
{\rm W}_{x_p}^{}
\bigl(1,{}^k X^{\theta}(t,x)\bigr) = 
\overset{**}\Xi_{\theta p}^{}(t,x)
$
\ for all
$(t,x)\in \Pi^\prime, \ \ \
\theta =1,\ldots,s,\ p =k+1,\ldots, n,
\hfill
$
\\[1.85ex]
\mbox{}\quad
$
{\rm W}_{x_p}^{}
\bigl({}^k X^{\nu}(t,x)\bigr) = 
\overset{**}\Xi_{\nu p}^{}(t,x)
$
\ for all
$(t,x)\in \Pi^\prime, \ \ \
\nu =s+1,\ldots,m,\ p =k+1,\ldots,n,
\hfill
$
\\[3ex]
is satisfied. Here the vector functions 
\vspace{0.5ex}
$
{}^k\;\!\!X^{j}\colon (t,x)\to \bigl(X_{1j}(t,x),\ldots,X_{kj}(t,x)\bigr)$
for all $(t,x)\in 
\linebreak
\in\Pi, \,
 j=1,\ldots,m,$
${\rm W}_{t_\zeta}$ and ${\rm W}_{x_ p}$
are 
correspondingly the Wronskians with respect to $t_\zeta$ and $x_p,\, 
 \zeta =s+1,\ldots, m,\, p =k+1,\ldots, n,$
the scalar functions 
$\overset{*}\Xi_{j \zeta}\colon \Pi^\prime\to \K, \ j=1,\ldots, m,
\, \zeta =
\linebreak
=s+1,\ldots, m,$
and
$\overset{**}\Xi_{j p}\colon \Pi^\prime\to \K, \ j=1,\ldots, m,\, p =k+1,\ldots, n,$
are such that
\\[2.5ex]
\mbox{}\hfill                                           
$
\overset{*}\Xi_{j \zeta}^{}(t,x)_{\displaystyle |_{\scriptstyle 
w({}^st,{}^kx)=0}}=0
$
\ for all 
$(t,x)\in \Pi^{\prime},
 \ \ \
j =1,\ldots,m,
\ \zeta =s+1,\ldots, m,
\hfill 
$
\\
\mbox{}\hfill (3.19)
\\
\mbox{}\hfill 
$
\overset{**}\Xi_{j p}^{}(t,x)_{\displaystyle |_{\scriptstyle 
w({}^st,{}^kx)=0}}=0
$
\ for all 
$(t,x)\in \Pi^{\prime},
 \ \ \
j =1,\ldots,m,
\ p =k+1,\ldots, n.
\hfill
$
\\[2.75ex]
\indent
So the 
{\sl necessary test of  existence of 
$s\!$-nonautonomous
$(n - k)\!$-cylindrical partial integral for 
total differential system} is proved.
\vspace{0.2ex}

{\bf Theorem 3.4.}
{\it
The system of identities
{\rm (3.18)} with {\rm (3.19)} is a necessary condition of  existence  
of the  $s\!$-nonautonomous $(n - k)\!$-cylindrical partial integral  {\rm (3.14)}  
for system} (TD).
\\[2ex]
\indent
{\bf
3.2.2. Criterion  of existence of 
{\boldmath$s\!$}-nonautonomous
{\boldmath$(n-k)\!$}-cylindrical partial integral.}
Let $n\times m$ matrix $X$ of  system (TD)  satisfies  the  conditions (3.18) with (3.19).
Let's write the functional system 
\\[2.25ex]
\mbox{}\hfill                            
$
\psi_\theta^{} +
{}^k\;\!\!X^\theta(t,x)\, {}^k\varphi = H_\theta^{}(t,x),
\ \ \theta=1,\ldots,s,
\hfill
$
\\[2ex]
\mbox{}\hfill
$
\partial^{\xi}_{{}_{\scriptstyle t_\zeta}}
\;{}^k\;\!\!X^\theta(t,x)\ 
{}^k\varphi = \partial^{\xi}_{{}_{\scriptstyle t_\zeta}} H_\theta^{}(t,x),
\ \ \theta =1,\ldots,s, \
\zeta =s+1,\ldots,m, \
\xi =1,\ldots,k,
\hfill
$
\\[2ex]
\mbox{}\hfill
$
\partial^{\xi}_{{}_{\scriptstyle x_p}}
\;{}^k\;\!\!X^\theta(t,x)\ {}^k\varphi =
\partial^{\xi}_{{}_{\scriptstyle x_p}}H_\theta^{}(t,x),
\ \ \theta =1,\ldots,s, \
p =k+1,\ldots,n, \
\xi =1,\ldots,k,
\hfill
$
\\[2ex]
\mbox{}\hfill
$
{}^k\;\!\!X^\nu(t,x)\, 
{}^k\varphi = H_\nu^{}(t,x),
\ \ \nu =s+1,\ldots,m,
$
\hfill (3.20)
\\[2.2ex]
\mbox{}\hfill
$
\partial^{\xi}_{{}_{\scriptstyle t_\zeta}}
\,{}^k\;\!\!X^\nu(t,x)\ {}^k\varphi = 
\partial^{\xi}_{{}_{\scriptstyle t_\zeta}}
H_\nu^{}(t,x),
\ \ \nu =s+1,\ldots,m,\
\zeta =s+1,\ldots,m,\
\xi =1,\ldots,k-1, 
\hfill
$
\\[2ex]
\mbox{}\hfill
$
\partial^{\xi}_{{}_{\scriptstyle x_p}}
\,{}^k\;\!\!X^\nu(t,x) \ {}^k\varphi
= 
\partial^{\xi}_{{}_{\scriptstyle x_p}}
H_\nu^{}(t,x),
\ \
\nu=s+1,\ldots, m,  \
p =k+1,\ldots,n, \
\xi =1,\ldots,k-1,
\hfill
$
\\[3ex]
where the vector functions
\vspace{0.5ex}
$
{}^s \psi\colon (t,x)\to
\bigl( \psi_1^{}({}^st,{}^kx),
\ldots,
\psi_s^{}({}^st,{}^kx)\bigr)
$
for all $(t,x)\in \Pi^{\prime}
$
and 
$
{}^k \varphi\colon (t,x)\to
\bigl( \varphi_1^{}({}^st,{}^kx),
\ldots,
\varphi_k^{}({}^st,{}^kx)\bigr)
$
for all $(t,x)\in \Pi^{\prime}$
are unknown,
\vspace{0.75ex}
the vector functions
$
{}^k\;\!\!X^j\colon (t,x)\to
\bigl( X_{1j}(t,x),\ldots,X_{kj}(t,x)\bigr)
$
for all
\vspace{0.5ex}
$(t,x)\in \Pi, \
j =1,\ldots,m,\, 0\leq k\leq n,
$
the scalar functions 
$H_j\colon \Pi^\prime\to \K
,\ j=1,\ldots, m,$ are such that
\\[2.25ex]
\mbox{}\hfill         
$
H_j^{}(t,x)_{\displaystyle |_{\scriptstyle w({}^st,{}^kx)=0}}=0
$
\ for all 
$(t,x)\in \Pi^{\prime},
 \ \ \
j =1,\ldots,m.
$
\hfill (3.21)
\\[2.75ex]
\indent
{\bf Theorem  3.5} ({\sl 
criterion  of existence of 
$s\!$-non\-auto\-nomous
$(n - k)\!$-cylindrical partial integral for 
total differential system}).
{\it
For  system  {\rm (TD)} 
to have 
$s\!$-nonautonomous
$(n - k)\!$-
\linebreak
cylindrical partial integral {\rm(3.14)}
it is necessary and sufficient that there exist the vector functions
${}^s\psi,\, {}^k\varphi$
and scalar  functions $H_j,\, j=1,\ldots, m,$
with {\rm(3.21),}
satisfying functional system {\rm(3.20)}, 
\vspace{0.3ex}
such that the  Pfaffian equation 
{\rm(1.21)} has the general integral 
$w\colon ({}^st,{}^kx)\to 
w({}^st,{}^kx)$ for all $({}^st,{}^kx)\in \widetilde{\Pi}^{s+k},$
\vspace{0.3ex}
where  the domain $\widetilde{\Pi}^{\;\!s+k}$
is the natural projection of domain  $\Pi^{\;\!\prime}$ on 
coordinate subspace $O\ {}^st\  {}^kx.$
}
\vspace{0.3ex}

{\sl Proof. Necessity}. Let system  (TD)
has the  
$s\!$-nonautonomous
$(n - k)\!$-cylindrical  partial  integral (3.14)
on the domain  $\Pi^{\;\!\prime}.$ Then, the identities (3.17) 
with (3.18) are satisfied.
By differentiating the first $s$ of this identities  $k$ times with respect to 
$t_{s+1},\ldots,t_m$ and $k$ times with respect to $x_{k+1},\ldots,x_n$
and by differentiating the rest $m-s$  identities
$k-1$ times with respect to 
$t_{s+1},\ldots,t_m$ and $k-1$ times with respect to $x_{k+1},\ldots,x_n$
\vspace{0.5ex}
we conclude  that the 
extensions 
on the domain $\Pi^{\;\!\prime}$ of the functions
\vspace{0.5ex}
$
{}^s\psi\colon ({}^st,{}^kx)\to
\bigl( \partial_{t_1}^{} w({}^st,{}^kx),\ldots,
\partial_{t_s}^{} w({}^st,{}^kx)\bigr)$
for all $({}^st,{}^kx)\in \widetilde{\Pi}^{s+k}
$
and
\vspace{0.25ex}
$
{}^k\varphi\colon ({}^st,{}^kx)\to
\bigl(\partial_{x_1}^{}w({}^st,{}^kx),\ldots,
\partial_{x_k}^{}w({}^st,{}^kx)\bigr)
$ for all $({}^st,{}^kx)\in 
\linebreak
\in
\widetilde{\Pi}^{s+k}$
\vspace{0.3ex}
is a solution to the functional system (3.17) with (3.16).
From this it also follows that 
the function  (3.14) is a  
general integral on the domain $\widetilde{\Pi}^{s+k}\subset \K^{s+k}$
of the  Pfaffian equation   (1.21).

{\sl Sufficiency}. Let the  vector functions
$
{}^s\psi\colon (t,x) \to  {}^s\psi({}^st,{}^kx),\ {}^k\varphi\colon (t,x)\to {}^k\varphi({}^st,{}^kx)$ for all $(t,x)\in \Pi^{\prime}$
be a solution to the functional system (3.20) with (3.21) and the Pfaffian equation  (1.21)
which is constructed on its base has the 
general integral 
$w\colon ({}^st,{}^kx)\to 
w({}^st,{}^kx)$ for all $({}^st,{}^kx)\in \widetilde{\Pi}^{s+k}.$
Then, the system of  identities 
\\[2.4ex]
\mbox{}\hfill                                     
$
\partial_{ t_\zeta^{}}^{} w({}^st,{}^kx) -
\mu({}^st,{}^kx)\;\! \psi_\zeta^{}({}^st,{}^kx) = 0
$
\ for all
$({}^st,{}^kx)\in \widetilde{\Pi}^{s+k}, \ \
\zeta =1,\ldots,s,
\hfill
$
\\[-0.1ex]
\mbox{}\hfill(3.22)
\\[-0.1ex]
\mbox{}\hfill
$
\partial_{x_\xi^{}}^{} w({}^st,{}^kx) -
\mu({}^st,{}^kx)\;\!\varphi_\xi^{}({}^st,{}^kx) = 0
$
\ for all 
$({}^st,{}^kx)\in \widetilde{\Pi}^{s+k},
\ \ \xi =1,\ldots,k,
\hfill
$
\\[2.5ex]
is satisfied,
where  $\mu\colon ({}^st,{}^kx)\to \mu({}^st,{}^kx)$ for all $({}^st,{}^kx)\in \widetilde{\Pi}^{s+k}$  
is  a  holomorphic along the manifold  $w({}^st,{}^kx)=0$
integrating multiplier of the Pfaffian equation  (1.21) 
which corresponds to its general integral 
$w\colon ({}^st,{}^kx)\to 
w({}^st,{}^kx)$ for all $({}^st,{}^kx)\in \widetilde{\Pi}^{s+k}.$

Taking into account that the vector functions 
${}^s\psi, \,{}^k\varphi$
are the solution to the functional system (3.20) with (3.21)
we receive the system of identities (3.17)
with $\Phi_j(t,x)\!\!~=\!\!~\mu({}^st,{}^kx) H_j(t,x)$ for all $(t,x)\in \Pi^{\prime},\  j =1,\ldots,m.$

Therefore the function (3.14) is   
an $s\!$-nonautonomous
$(n - k)\!$-cylindrical  partial  integral
on the dimain $\Pi^{\;\!\prime}$  of   system (TD). \k
\vspace{0.5ex}

{\bf Example 3.3.}
The real completely solvable
autonomous   total differential system 
\\[2.25ex]
\mbox{}\hfill                  
$
\begin{array}{c}
dx_1 = {}-\bigl(x_2 + x_1(x_1^2 + x_2^2 + x_3^2)\bigr)\;\!dt_1
- x_1(x_1^2 + x_2^2 + x_3^2)\,dt_2,
\\[2ex]
dx_2 = \bigl(x_1 - x_2(x_1^2 + x_2^2 + x_3^2)\bigr)\;\!dt_1
- x_2(x_1^2 + x_2^2 + x_3^2)\,dt_2,
\\[2ex]
dx_3 = x_3(x_1^2 + x_2^2 + x_3^2)\,(dt_1 + dt_2)
\end{array}
$
\hfill (3.23)
\\[2.5ex]
has the autonomous  2-cylindrical  partial  integral
$w\colon (t, x) \to x_1^2 + x_2^2$ for all $(t, x) \in \R^5.$
On the  coordinate plane $Ox_1x_2$ of phase space $\R^3$
this partial integral specifies the isolated point  $x_1= x_2 = 0$  
and 
for this point
the hypotheses of Theorem 11 from [20] are satisfied, when
\\[2ex]
\mbox{}\hfill
$
\partial_{t_1}^{} w(t,x)_{_{\scriptstyle
\Bigl|\!\!\!\begin{array}{l} 
\scriptstyle (3.23)
\\[-.4ex]
\scriptstyle x_3 = 0
\end{array}
}} \ = \,
\partial_{t_2}^{} w(t,x)_{_{\scriptstyle
\Bigl|\!\!\!\begin{array}{l} 
\scriptstyle (3.23)
\\[-.4ex]
\scriptstyle x_3 = 0
\end{array}
}} 
\ =
{}-2 (x_1^2 + x_2^2)^2 \leq 0
$
\ for all
$(x_1,x_2) \in \R^2.
\hfill
$
\\[3ex]
\indent
So the zero solution  $x_1= x_2 = x_3 = 0$  to  system (3.23)
is stable on the plane $Ox_1x_2.$

The equilibrium point $O(0,0,0)$ of the 
induced by system (3.23)
autonomous  ordinary differential  system 
\\[2ex]
\mbox{}\hfill
$
\dfrac{dx_1}{dt_1} = {}-x_2 - x_1(x_1^2 + x_2^2 + x_3^2),
\quad
\dfrac{dx_2}{dt_1} = x_1 - x_2(x_1^2 + x_2^2 + x_3^2),
\quad
\dfrac{dx_3}{dt_1} =
x_3(x_1^2 + x_2^2 + x_3^2)
\hfill
$
\\[3ex]
is unstable by Chetaev's theorem
[19, pp. 19 -- 20] with $V(x_1, x_2, x_3) = {}-x_1^2-x_2^2 +x_3^2.$ 
\vspace{0.3ex}

Therefore the zero solution  $x_1= x_2 =x_3 = 0$  to  system (3.23)
is unstable.
\\[2ex]
\indent
{
\bf
3.2.3.    Functionally independent 
{\boldmath$s\!$}-nonautonomous
{\boldmath$(n-k)\!$}-cylindrical partial integrals.}
The method which is proposed in Theorem 3.5
can be used to construct  
the fun\-ctionally independent 
$s\!$-nonautonomous $(n-k)\!$-cylindrical partial integrals
of system (TD).
\\[0.5ex]
\indent
{\bf Theorem  3.6.}
{\it
Let $h$ functional systems  {\rm(3.20)} with {\rm(3.21)} has 
$q$ 
not  linearly bound 
on the domain $\Pi^{\;\!\prime}\subset \Pi$  solutions {\rm(1.25)}
and for each of them the corresponding  Pfaffian equation {\rm (1.26)}
has the general integral}
\\[2.35ex]
\mbox{}\hfill                                           
$
w_{\gamma}^{}\colon
({}^st,{}^kx)\to w_{\gamma}^{}({}^st,{}^kx)
$
\ for all
$({}^st,{}^kx) \in \widetilde{\Pi}^{s+k}\subset \K^{s+k},
\ \ \ \gamma=1,\ldots,q,
$
\hfill (3.24)
\\[2.25ex]
{\it
on the domain $\widetilde{\Pi}^{\;\!s+k}$ which is the 
\vspace{0.2ex}
natural projection of  domain $\Pi^{\;\!\prime}$ on  coordinate
subspace $O\ {}^s t\ {}^k x.$ Then, the general integrals 
{\rm (3.24)} are functionally independent
\vspace{0.3ex}
on the domain  $\widetilde{\Pi}^{\;\!s+k}.$}

{\sl Proof.} By virtue of the system of identities  (3.22) 
\\[2.25ex]
\mbox{}\hfill
$
\partial_{t_\zeta}^{} w_\gamma^{}({}^st,{}^kx) -
\mu_\gamma^{}({}^st,{}^kx)\;\! \psi^\gamma_{_{\scriptstyle \zeta}}({}^st,{}^kx)
=0
$
\ for all
$({}^st,{}^kx) \in \widetilde{\Pi}^{s+k},
\ \ \ \zeta=1,\ldots,s,
\ \gamma=1,\ldots,q,
\hfill
$
\\[2ex]
\mbox{}\hfill
$
\partial_{x_\xi}^{} w_\gamma^{}({}^st,{}^kx) 
-
\mu_\gamma^{}({}^st,{}^kx)\;\!
\varphi^\gamma_{_{\scriptstyle \xi}}({}^st,{}^kx)=0
$ \ for all
$({}^st,{}^kx) \in \widetilde{\Pi}^{s+k},
\ \  \ \xi=1,\ldots,k,
 \ \gamma=1,\ldots,q.
\hfill
$
\\[2.75ex]
So the Jacobi's matrix
\vspace{0.5ex}
$
J\bigl(w_{\gamma}^{}({}^st,{}^kx);\, {}^st,{}^kx\bigr)=
\bigl\|\Psi({}^st,{}^kx)\Phi({}^st,{}^kx)\bigr\|$
for all
$({}^st,{}^kx) \in 
\widetilde{\Pi}^{s+k},
$
where the matrix  $\|\Psi\Phi\|$ consists of 
\vspace{0.75ex}
$q\times s$ matrix
$
\Psi({}^st,{}^kx) =
\bigl\|\mu_\gamma({}^st,{}^kx)\;\!\psi_\zeta^{\gamma}({}^st,{}^kx)\bigr\|$
for all $({}^st,{}^kx) \in \widetilde{\Pi}^{s+k}$
and  $q\times k$ matrix
$
\Phi({}^st,{}^kx) =
\bigl\|
\mu_\gamma^{}({}^st,{}^kx)\;\! \varphi_\xi^{\gamma}({}^st,{}^kx)\bigr\|$
\vspace{0.75ex}
for all $({}^st,{}^kx) \in \widetilde{\Pi}^{s+k}.$

Since the vector functions (1.25)
\vspace{0.3ex}
are not  linearly bound 
on the domain $\widetilde{\Pi}^{s+k}$
the rank of Jacobi's matrix
$
{\rm rank}\, J\bigl(w_{\gamma}^{}({}^st,{}^kx);\,
{}^st,{}^kx\bigr) = q
$
nearly everywhere on the domain  $\widetilde{\Pi}^{s+k}.$
Therefore 
the  general integrals 
(3.24) of the Pfaffian equations  (1.26)
\vspace{0.3ex}
 are functionally independent
on the domain  $\widetilde{\Pi}^{s+k}.$  \k
\vspace{0.5ex}

{\bf Example 3.4.}
The  system 
of equations in total differentials
\\[2.3ex]
\mbox{}\hfill                                        
$
dx_i = x_i\, \biggl[ \dfrac{t_2-1}{t_1(t_2-t_1)}\ dt_1 -
\dfrac{t_1-1}{t_2(t_2-t_1)}\ dt_2\biggr],
\ \ i=1,2,3,
$
\hfill (3.25)
\\[2.5ex]
is not completely solvable since  
the expression in  square brackets
is not the exact differential under  independent
variables $t_1$ and $t_2.$

The associated  normal linear homogeneous partial system 
\\[2ex]
\mbox{}\hfill
$
{\frak X}_1(t,x)\,y=0,\ \ 
{\frak X}_2(t,x)\,y=0,                                
\hfill
$
\\[2ex]
which is constructed on the base of 
operators of differentiation
by virtue of  system (3.25) 
\\[2ex]
\mbox{}\hfill                                       
$
\displaystyle
{\frak X}_1(t,x)=
\partial_{t_1}^{} + \dfrac{t_2 - 1}{t_1(t_2 - t_1)}\;
\sum\limits_{i=1}^3\,
x_i\;\! \partial_{x_i}^{}
$
\ for all 
$(t,x)\in\Pi,
\hfill
$
\\[2.25ex]
\mbox{}\hfill
$
{\frak X}_2(t,x)=
\partial_{t_2}^{} - \dfrac{t_1 - 1}{t_2(t_2 - t_1)}\;
\sum\limits_{i=1}^3\,
x_i\;\! \partial_{x_i}^{}
$
\ for all 
$(t,x)\in\Pi,
\hfill
$
\\[3ex]
is incomplete on every domain $\Pi$
from the set $\{(t,x)\colon t_1\ne0,\, t_2\ne0,\,t_2\ne t_1\}$
and has the defect $\delta=1.$
\vspace{0.2ex}

Therefore an integral basis of system (3.25)
consists of  $n-\delta=3-1=2$
functionally independent first integrals.

The system (3.25)  has  the autonomous
2-cylindrical partial integrals
\\[2ex]
\mbox{}\hfill
$
w_i\colon (t,x)\to x_i
$
\ for all 
$(t,x)\in\Pi,
\ \ \ i=1,2,3,
\hfill
$
\\[3ex]
since
${\frak X}_1^{} x_i=
\dfrac{t_2 - 1}{t_1(t_2 - t_1)}\ x_i$ for all $(t,x)\in\Pi,  \,\ 
{\frak X}_2^{} x_i={}-\dfrac{t_1 - 1}{t_1(t_2 - t_1)}\ x_i$
for all $(t,x)\in\Pi.$
\vspace{1ex}

From system (3.25) we get 
\\[2ex]
\mbox{}\hfill
$
\dfrac{dx_1}{x_1}=\dfrac{dx_2}{x_2}=
\dfrac{dx_3}{x_3}\,.
\hfill
$
\\[2.5ex]
\indent
From this by immediate integration we get that 
on every  domain  $\Pi_i\subset 
\{(t,x)\colon t_1\ne0,$
$t_2\ne0,\,t_2\ne t_1,\, x_i\ne0\},\, i=1,2,3,$
the system (3.25) has the  integral basis which consists of two
functionally independent autonomous 1-cylindrical first integrals
\\[2.25ex]
\mbox{}\hfill
$
F_{i\xi}^{}\colon (t,x)\to x_{\xi}^{}\,x_i^{{}-1}$
\ for all 
$(t,x)\in \Pi_i,
\qquad
F_{i\theta}^{}\colon (t,x)\to x_{\theta}^{}\,x_{i}^{{}-1}$
\ for all $x\in \Pi_i,
\hfill
$
\\[2ex]
\mbox{}\hfill
$
i=1,2,3, \ \ \xi =1,2,3,\ \ \theta=1,2,3,  
\ \ \xi\ne i, \ \  
 \theta\ne i, \ \  \theta\ne \xi.
\hfill
$
\\[-2ex]

\newpage

\mbox{}
\\[-1ex]
\centerline{
\bf
3.3.  Functional relations between general solutions to irreducible
}
\\[0ex]
\centerline{
\bf
Painlev\'{e} equations
}
\\[2ex]
\indent
Let's consider the twelfth-order differential system
\\[2ex]
\mbox{}\hfill                
$
\dfrac{dx_i}{dt} = y_i ,
\ \ \
\dfrac{dy_i}{dt} =P_i(t,x_i,y_i),
\ \ \ 
i = 1,\ldots, 6,
$
\hfill (PS)
\\[2ex]
where 
\\[1ex]
\mbox{}\hfill                
$
P_{1}^{}\colon (t,x,y)\to 6x_1^2 + t
$
\  for all $(t,x,y)\in \C^{13},
\hfill                
$
\\[2ex]
\mbox{}\hfill                
$
P_{2}^{}\colon (t,x,y)\to 2x_2^3 + \alpha_2 +tx_2$
\ for all $(t,x,y)\in \C^{13},\ \ \alpha_2\in \C,
\hfill                
$
\\[3ex]
\mbox{}\hfill                
$
P_{3}^{}\colon (t,x,y)\to ({}-y_3 + \alpha_3x_3^2 +
\beta_3)t^{{}-1} +y_3^2x_3^{{}-1} + \gamma_3x_3^3 +
\delta_3x_3^{{}-1}
\hfill                
$
\\[2ex]
\mbox{}\hfill                
for all 
$(t,x,y)\in \{(t,x,y)\colon t\ne 0,\,x_3\ne 0\}, \ \ 
\alpha_3,\beta_3,\gamma_3,\delta_3\in \C,
\hfill                
$
\\[3ex]
\mbox{}\hfill                
$
P_{4}^{}\colon (t,x,y)\to
\dfrac{1}{2}\ y_4^2x_3^{{}-1} +
\dfrac{3}{2}\ x_4^3 -
2\alpha_4x_4 + \beta_4x_4^{{}-1} +4tx_4^2 + 2t^2x_4
\hfill                
$
\\[2ex]
\mbox{}\hfill                
for all 
$(t,x,y)\in \{(t,x,y)\colon x_4\ne 0\}, \ \
\alpha_4,\beta_4\in \C,
\hfill                
$
\\[3ex]
\mbox{}\hfill                
$
P_{5}^{}\colon (t,x,y)\to 
(x_5 - 1)^2 (\alpha_5x_5
+ \beta_5x_5^{{}-1} )t^{{}-2} +
({}-y_5 +\gamma_5x_5)t^{{}-1} +
\dfrac{1}{2}\ (3x_5 -1)x_5^{{}-1} (x_5 - 1)^{{}-1}y_5^2 \ +
\hfill                
$
\\[2ex]
\mbox{}\hfill                
$
\!
+\,\delta_5x_5 (x_5 + 1)(x_5 - 1)^{{}-1} 
$
for all 
$(t,x,y)\!\!~\in\!\!~ \{(t,x,y)\colon t\!\!~\ne\!\!~ 0,\;\!x_5\!\!~\ne\!\!~ 0,\;\!x_5\!\!~\ne\!\!~ 1\}, \ 
\alpha_5,\beta_5,\gamma_5,\delta_5\!\!~ \in\!\!~ \C,
\!
\hfill                
$
\\[3ex]
\mbox{}\hfill                
$
P_{6}^{}\colon (t,x,y)\to 
\Bigl({}-\dfrac{1}{2}\ y_6^2 + y_6 - \delta_6\Bigr)
(t -x_6)^{{}-1} + 
(x_6 - 1)^2 (\alpha_6x_6
+\beta_6x_6^{{}-1})(t - 1)^{{}-2} \ +
\hfill                
$
\\[2ex]
\mbox{}\hfill                
$
+\,\bigl({}-y_6 + \alpha_6x_6(x_6 - 1)(1 - 2x_6) - \beta_6(x_6 - 1) +
\gamma_6x_6 + \delta_6x_6\bigr)(t - 1)^{{}-1} \,+
\hfill                
$
\\[2ex]
\mbox{}\hfill                
$
+\,
x_6^2\bigl( \alpha_6(x_6 -1) - \gamma_6( x_6 -
1)^{{}-1}\bigr)t^{{}-2} +
\bigl({}-y_6 + \alpha_6x_6(x_6 - 1)(2x_6 - 1) +
\beta_6(x_6 - 1) - \gamma_6x_6 \,-
\hfill                
$
\\[2ex]
\mbox{}\hfill                
$
-\, \delta_6(x_6 -
1)\bigr)t^{{}-1} +
\dfrac{1}{2}\,\bigl(x_6^{{}-1} + (x_6 - 1)^{{}-1}\bigr)y_6^2
\hfill                
$
\\[2ex]
\mbox{}\hfill                
for all 
$(t,x,y)\in
\{(t,x,y)\colon t\ne 0,\,t\ne 1,\,x_6\ne 0,\,x_6\ne 1, \,x_6\ne t\},  
\ \
\alpha_6, \beta_6, \gamma_6, \delta_6 \in \C,
\hfill                
$
\\[2ex]
\mbox{}\hfill                
$
x=(x_1,\ldots,x_6),\ \
y=(y_1,\ldots,y_6).
\hfill
$
\\[3ex]
We'll consider the system   (PS) 
on  any domain $\Pi$ from the  
set
$
D=
\{(t,x,y)\colon t\ne 0,\, t\ne 1,$
$x_3\ne0,\, x_4\ne 0,\, 
x_5\ne 0,\, x_5\ne 1,\, x_6\ne0,\, x_6\ne 1,x_6\ne t\}.
$
\vspace{0.5ex}

The system (PS) is constructed on the base of six 
irreducible
Painlev\'{e} equations [21, pp. 463 -- 465]
\\[1.25ex]
\mbox{}\hfill                     
$
\dfrac{d^2x_1}{dt^2} = 6x_1^2 + t,
$
\hfill (P-1)
\\[2.5ex]
\mbox{}\hfill                     
$
\dfrac{d^2x_2}{dt^2}  =
2x_2^3 + tx_2 + \alpha_2\;\!,
$
\hfill (P-2)
\\[3ex]
\mbox{}\hfill                     
$
\dfrac{d^2x_3}{dt^2}  =
\dfrac{1}{x_3}\,\Bigl(\dfrac{dx_3}{dt}\Bigr)^{\!2}  -
\dfrac{1}{t}\ \dfrac{dx_3}{dt}  +
\dfrac{\alpha_3}{t}\ x_3^2  +
\dfrac{\beta_3}{t} + \gamma_3x_3^3  +
\dfrac{\delta_3}{x_3}\,,
$
\hfill (P-3)
\\[3.5ex]
\mbox{}\hfill                     
$
\dfrac{d^2x_4}{dt^2} =
\dfrac{1}{2x_4}\, \Bigl(\dfrac{dx_4}{dt}\Bigr)^{\!2}  +
\dfrac{3}{2}\ x_4^3  + 4tx_4^2 +
2(t^2 - \alpha_4 )x_4 + \dfrac{\beta_4}{x_4}\,,
$
\hfill  (P-4)
\\[3.5ex]
\mbox{}\hfill                     
$
\dfrac{d^2x_5}{dt^2} =
\dfrac{3x_5 -1}{2x_5(x_5 - 1)}\,
\Bigl(\dfrac{dx_5}{dt}\Bigr)^{\!2} -
\dfrac{1}{t}\ \dfrac{dx_5}{dt} +
\dfrac{\alpha_5}{t^2}\ x_5(x_5 - 1)^2\, +\ \ \ \ \ 
\hfill
$
\\[0ex]
\mbox{}\hfill (P-5)
\\[0ex]
\mbox{}\hfill
$
+ \,
\dfrac{\beta_5}{t^2}\ \dfrac{(x_5 -1)^2 }{x_5} +
\dfrac{\gamma_5}{t}\ x_5 +\delta_5\,
\dfrac{x_5(x_5 + 1)}{x_5 - 1}\,, \ \ \ 
\hfill
$
\\[3.5ex]
\mbox{}\hfill                     
$
\dfrac{d^2x_6}{dt^2}  =
\dfrac{1}{2}\,\Bigl(\dfrac{1}{x_6} + \dfrac{1}{x_6 - 1} \!+\!
\dfrac{1}{x_6 - t}\Bigr) \Bigl(\dfrac{dx_6}{dt} \Bigr)^2 -
\Bigl(\dfrac{1}{t} + \dfrac{1}{t - 1} + \dfrac{1}{x_6 -
t}\Bigr)\dfrac{dx_6}{dt} \ +
\hfill 
$
\\[0.2ex]
\mbox{}\hfill (P-6)
\\[0.2ex]
\mbox{}\hfill
$
+\,
\dfrac{x_6 (x_6 - 1) (x_6 - t) }{t^2(t -
1)^2}\,\left(\alpha_6 + \beta_6\ \dfrac{t}{x_6^2} +
\gamma_6\,\dfrac{t - 1}{(x_6 - 1 )^2} +
\delta_6\,\dfrac{t(t - 1)}{(x_6 - t)^2} \right).
\hfill
$
\\[3ex]
\indent
Therefore  we name  the system (PS) the Painlev\'{e} system.

The components $x_i,\, i = 1,\ldots, 6,$ of the general solution to system (PS) 
are the general solutions to the irreducible
Painlev\'{e} equations ($\!\text{P-\;\!}i$\!), $i = 1,\ldots, 6,$
correspondingly,
and the components $y_i,\, i = 1,\ldots, 6,$ are their derivatives.

Thus, the problem on the relations between the general solutions to the 
irreducible Painlev\'{e} equations ($\!\text{P-\;\!}i$\!), $i = 1,\ldots, 6,$
and  their derivatives
\\[2ex]
\mbox{}\hfill                     
$
w (x_1, \ldots,x_6,\;\! {\sf D}x_1,\ldots, {\sf D}x_6) = 0,
$
\hfill (3.26)
\\[2.25ex]
where  $x_i,\, i = 1,\ldots, 6,$
are the general solutions to  equations ($\!\text{P-\;\!}i$\!), $i = 1,\ldots, 6,$
correspon\-dingly,
is equivalent to the problem of finding
autonomous partial integrals 
$w\colon (x,y)\!\!~\to\!\!~ w(x,y)$ for all $(x,y)\in \widetilde{\Pi}^{12},$
where  $\widetilde{\Pi}^{12}$ is  the 
natural projection of  domain $\Pi$ on the phase 
space   $\C^{12},$
for the differential Painlev\'{e} system (PS).

The existence of 
autonomous 6-cylindrical partial integrals 
$w\colon (x)\to w(x)$ for all $(x,y) \in 
\linebreak
\in \widetilde{\Pi}^{6},$
where  $\widetilde{\Pi}^{6}$ is  the 
natural projection of  domain $\Pi$ on the phase 
subspace   $Ox,$
for the differential Painlev\'{e} system (PS)
determines the relation $w$ between the general solutions to the 
Painlev\'{e} equations ($\!\text{P-\;\!}i$\!), $i = 1,\ldots, 6.$
Otherwise, there are no 
relations given by ho\-lomor-
\linebreak
phic function $w$
between the general solutions to the 
Painlev\'{e} equations ($\!\text{P-\;\!}i$\!), $i = 1,\ldots, 6.$

Let's seek autonomous partial integrals 
$w\colon (x,y)\to w(x,y)$ for all $(x,y)\in \widetilde{\Pi}^{12}$
of  the  Painlev\'{e} system (PS).
To this end, we compose a system of the form (3.20)
following Theorem 3.5 and consider the first equation of this system
\\[2ex]
\mbox{}\hfill
$
\displaystyle
\sum_{i = 1}^{6}\,
y_i\;\!\varphi_i +
\sum_{i = 7}^{12} \,
P_{i - 6}(t,x,y)\varphi_i = H_1(t,x,y),
\hfill
$
\\[1.5ex]
where
\\[1.5ex]
\mbox{}\hfill
$
H_1\colon(t,x,y)\to  W_{-1}(x,y)(t - x_6)^{{}-1}
+ V_{-2}(x, y)(t - 1)^{{}-2} + 
V_{-1}(x,y)(t - 1)^{{}-1}\,+
\hfill
$
\\[1.8ex]
\mbox{}\hfill
$
+\, 
U_{-2}(x, y)t^{{}-2} + U_{-1}(x, y)t^{{}-1}  +
U_0(x, y)  + U_1(x, y)t + U_2(x, y)t^2
\hfill
$
\\[3ex]
and 
\vspace{1ex}
$W_{-1}(x,y)_{\displaystyle |_{\scriptstyle w(x,y)=0 }}= 0,\
V_{-2}(x,y)_{\displaystyle |_{\scriptstyle w(x,y)=0 }}=0,\
V_{-1}(x,y)_{\displaystyle |_{\scriptstyle w(x,y)=0 }}=0,\
U_j(x,y)_{\displaystyle |_{\scriptstyle w(x,y)=0 }}=
\linebreak
= 0, \ j = {}-2,\ldots,2.$
\vspace{0.5ex}

Hence, by choosing  
\\[2ex]
\mbox{}\hfill
$
t^2, \ \ t, \ \ 1, \ \ t^{{}-1},
\ \ t^{{}-2}, \ \ (t - 1)^{{}-1}, \ \ (t - 1)^{{}-2}, \ \
(t -x_6)^{{}-1}
\hfill
$
\\[2.25ex]
as a basis, we obtain the system
\\[2ex]
\mbox{}\hfill                     
$
\Bigl({}-\dfrac{1}{2}\ y_6^2 + y_6 - \delta_6\Bigr)\varphi_{12} =W_{-1}(x, y),
\hfill
$
\\[2.5ex]
\mbox{}\hfill
$
(x_6 - 1)^2(\alpha_6x_6+ \beta_6x_6^{-1})\;\!\varphi_{12}=
V_{-2}(x, y),
\hfill
$
\\[2.5ex]
\mbox{}\hfill
$
\bigl({}-y_6 + \alpha_6x_6
(x_6 -1)(1 - 2x_6) - \beta_6(x_6 -1) +
\gamma_6x_6 + \delta_6x_6\bigr)\varphi_{12}=
 V_{-1}(x, y),
\hfill
$
\\[2.25ex]
\mbox{}\hfill
$
(x_5 - 1)^2(\alpha_5x_5 + \beta_5x_5^{{}-1})\varphi_{11}
+ 
x_6^2\bigl(\alpha_6(x_6 - 1) -
\gamma_6\bigl(x_6 - 1\bigr)^{{}-1}\bigr)\varphi_{12}=
U_{-2}(x, y),
\hfill
$
\\[2.5ex]
\mbox{}\hfill
$
({}-y_3 + \alpha_3x_3^2 + \beta_3)\varphi_9 +
 ({}-y_5 + \gamma_5x_5)\varphi_{11} +
\hfill
$
\\[1.8ex]
\mbox{}\hfill
$
+\,
\bigl({}-y_6 +
\alpha_6x_6(x_6 - 1)(2x_6-1) +
\beta_6(x_6-1) - \gamma_6x_6 - \delta_6(x_6-1)\bigr)\varphi_{12} =
U_{-1}(x, y),
\hfill
$
\\[1.5ex]
\mbox{}\hfill (3.27)
\\[-1.5ex]
\mbox{}\hfill
$
\displaystyle
\sum\limits_{i = 1}^{6}\,y_i\;\!\varphi_i +
6x_1^2\varphi_7 + (2x_2^3 + \alpha_2)\;\!\varphi_8 \,+
(y_3^2x_3^{{}-1} +
\gamma_3x_3^3 + \delta_3x_3^{{}-1})\varphi_9 \,+
\hfill
$
\\[1ex]
\mbox{}\hfill
$
+\,
\Bigl(\dfrac{1}{2}\ y_4^2x_4^{{}-1} +
\dfrac{3}{2}\ x_4^3 - 2\alpha_4x_4 + \beta_4x_4^{-1}\Bigr) \varphi_{10}
\,+
\hfill
$
\\[1.8ex]
\mbox{}\hfill
$
+\,
\Bigl(\dfrac{1}{2}\
(3x_5 - 1)x_5^{{}-1}(x_5 -1)^{{}-1}y_5^2 +
\delta_5x_5(x_5 + 1)(x_5 -1)^{{}-1}\Bigr)\varphi_{11} \,+
\hfill
$
\\[1.8ex]
\mbox{}\hfill
$
+\,
\Bigl(\dfrac{1}{2}\,
\bigl(x_6^{{}-1} + (x_6 - 1)^{{}-1}\bigr)y_6^2\Bigr)\varphi_{12} = U_0(x, y),
\hfill
$
\\[3ex]
\mbox{}\hfill
$
\varphi_7(x, y) + x_2\;\!\varphi_8
+ 4x_4^2\;\!\varphi_{10} = U_1(x, y),
\hfill
$
\\[2.5ex]
\mbox{}\hfill
$
2x_4\varphi_{10} = U_2(x, y),
\hfill
$
\\[3ex]
where the functions 
$\varphi_{\xi}^{}\colon(x,y) \to  \varphi_{\xi}^{}(x,y)$
for all 
$(x,y) \in  \widetilde{\Pi}^{12},\  \xi =1,\ldots,12,$
are unknown.
\vspace{0.5ex}

{\bf Theorem  3.7.} 
{\it
There is no relation of the form {\rm (3.26)}
between the general solutions 
$x_i\colon t \to x_i(t,C_{i1},C_{i2}),\, i = 1,\ldots, 6,$ 
to the 
irreducible Painlev\'{e} equations {\rm($\!\text{P-\;\!}i$\!)}, $i = 1,\ldots, 6,$
and  their derivatives.
} 
\vspace{0.3ex}

{\sl Proof}.
Let $w\colon (t,x,y) \to w(x,y)$ for all $(t,x,y)\in \Pi$ 
be an autonomous 
nonconstant
partial integral of  the  Painlev\'{e} system (PS).

It follows from the last equation of system (3.27)
that
\\[2ex]
\mbox{}\hfill                 
$
\varphi_{10}(x,y) =
\dfrac{1}{2}\ x_4^{{}-1}\  U_2(x, y).
$
\hfill (3.28)
\\[2.5ex]
\indent
Let's consider the identity (see system  (3.22))
\\[2ex]
\mbox{}\hfill                
$
\partial_{y_4}\;\!w(x, y) - \dfrac{1}{2}\
\mu(x,y)x_4^{{}-1}\ U_2(x, y) = 0
$
\hfill (3.29)
\\[2.5ex]
corresponding to the function $\varphi_{10}.$
Since the integrating multiplier $\mu$
is holomorphic along $w(x,y) = 0$
by Theorem 3.5,
we get from identity (3.29)
that  $w$ is a partial integral of the equation  $x_4\,\partial_{y_4}\;\!w = 0.$

From this it follows that 
\\[2ex]
\mbox{}\hfill                    
$
w = w(x, y_1,y_2, y_3, y_5, y_6),
$
\hfill (3.30)
\\[2ex]
that is,  $w$ is independent of  $y_4.$

Since an autonomous  partial integral is not identical constant,
the multiplier $\mu$ don't vanish identically. 
Therefore,  
\\[1ex]
\mbox{}\hfill              
$
\varphi_{10}(x, y) \equiv 0
$
\hfill (3.31)
\\[2ex]
by virtue of  (3.29) and (3.30) and on the base of representation (3.28).

From the first equation of system (3.27) we find 
\\[2ex]
\mbox{}\hfill                        
$
\varphi_{12}(x, y) =
\Bigl({}-\dfrac{1}{2}\ y_6^2 + y_6 -
\delta_6\Bigr)^{{}-1}\ W_{-1}(x, y).
$
\hfill (3.32)
\\[2.25ex]
\indent
Similarly, considering the identity (see system  (3.22))
\\[2ex]
\mbox{}\hfill                
$
\partial_{y_6}w(x, y) - \mu(x, y)
\Bigl({}-\dfrac{1}{2}\ y_6^2 +
y_6 - \delta_6\Bigr)^{{}-1}\ W_{-1}(x, y) = 0,
$
\hfill (3.33)
\\[2.25ex] 
corresponding to the function $\varphi_{12},$
we arrive at conclusion that  
 $w$ is a partial integral of the equation 
\\[0ex]
\mbox{}\hfill
$
\Bigl({}-\dfrac{1}{2}\ y_6^2 + y_6 -
\delta_6\Bigr) \partial_{y_6}w = 0.
\hfill
$
\\[3ex]
\indent
From this  taking into account (3.30)
it follows that
\\[2ex]
\mbox{}\hfill        
$
w = w (x, y_1, y_2,y_3,y_5),
$
\hfill (3.34)
\\[2.25ex]
that is,  $w$ is independent of  $y_4$ and  $y_6.$

From (3.32), (3.33), and (3.34) it follows that
\\[2ex]
\mbox{}\hfill          
$
\varphi_{12}(x, y) \equiv 0.
$
\hfill (3.35)
\\[2.25ex]
\indent
In view of (3.31) and (3.35) the functional system (3.27)
can be rewritten in the form 
\\[2.25ex]
\mbox{}\hfill                  
$
\displaystyle
(x_5 - 1)^2 (\alpha_5x_5 +
\beta_5x_5^{{}-1})\varphi_{11}  = U_{-2}(x, y),
\hfill
$
\\[2ex]
\mbox{}\hfill
$
\displaystyle
 ({}-y_3 + \alpha_3x_3^2 +\beta_3)\varphi_9 +
({}-y_5 + \gamma_5x_5)\varphi_{11} =
U_{-1}(x,y),
\hfill
$
\\[2ex]
\mbox{}\hfill
$
\displaystyle
\sum\limits_{i = 1}^{6}\, y_i\varphi_i +
6x_1^2 \varphi_7+  (2x_2^3 + \alpha_2)\varphi_8
+
 (y_3^2x_3^{{}-1}
\gamma_3x_3^3 +\delta_3x_3^{{}-1})\varphi_9
\,+
$
\hfill (3.36)
\\[2ex]
\mbox{}\hfill
$
+\, 
\Bigl(\dfrac{1}{2}\ (3x_5 -1)x_5^{{}-1}
(x_5 - 1)^{{}-1}y_5^2 +
\delta_5x_5(x_5+ 1)(x_5 - 1)^{{}-1}\Bigr)\varphi_{11} = U_0(x, y),
\hfill
$
\\[2.5ex]
\mbox{}\hfill
$
\displaystyle
\varphi_7 + x_2\varphi_8 = U_1(x, y).
\hfill
$
\\[2.5ex]
\indent
Let $|\alpha_5| + |\beta_5| \ne 0.$ Then, 
from first  equation of system (3.36) we get 
\\[1.5ex]
\mbox{}\hfill                 
$
\varphi_{11}(x, y) = (x_5 - 1)^{{}-2}\,
(\alpha_5x_5 +\beta_5x_5^{{}-1})^{{}-1}\ U_{-2}(x, y).
$
\hfill (3.37)
\\[2.2ex]
\indent
The identity (see system  (3.22))
\\[2ex]
\mbox{}\hfill                   
$
\partial_{y_5}w(x, y) - \mu(x,y)
(x_5 - 1)^{{}-2}
(\alpha_5x_5 + 
\beta_5x_5^{{}-1})^{{}-1} U_{-2}(x, y)= 0
$
\hfill (3.38)
\\[2.25ex]
corresponds  to the function $\varphi_{11}$
and by Theorem 3.5
a function 
 $w$ is a partial integral of the equation 
\\[0ex]
\mbox{}\hfill
$
(x_5 - 1)^2(\alpha_5x_5 +\beta_5x_5^{{}-1})\;\!
\partial_{y_5}w = 0.
\hfill
$
\\[2.5ex]
\indent
From this  taking into account (3.34)
it follows that
\\[2ex]
\mbox{}\hfill               
$
w = w(x, y_1, y_2,y_3),
$
\hfill (3.39)
\\[2.25ex]
that is,  $w$ is independent of  $y_4,\, y_5,\, y_6.$

So,  in view of  (3.37), (3.38),  and (3.39)
\\[2ex]
\mbox{}\hfill                 
$
\varphi_{11}(x, y) \equiv 0.
$
\hfill  (3.40)
\\[2.25ex]
\indent
From the second equation of system (3.36) with (3.40) it follows that
\\[2ex]
\mbox{}\hfill
$
\varphi_9(x, y) = ({}-y_3 + \alpha_3x_3^2 +
\beta_3)^{{}-1}\ U_{-1}(x, y).
\hfill
$
\\[2.25ex]
\indent
Next  we consider the identity 
\\[2ex]
\mbox{}\hfill
$
\partial_{y_3}w(x, y) - \mu(x, y)
({}-y_3 + \alpha_3x_3^2 +\beta_3)^{{}-1}\ U_{-1}(x, y)
= 0
\hfill
$
\\[2.25ex]
and (by Theorem 3.5)  taking into account   (3.39) 
we establish that  
\\[2ex]
\mbox{}\hfill             
$
w = w(x, y_1, y_2),
$
\hfill (3.41)
\\[2.25ex]
that  is,  $w$ is independent of   $y_{\tau},\, \tau = 3,\ldots, 6.$
Therefore,
\\[2ex]
\mbox{}\hfill           
$
\varphi_9(x, y) \equiv  0.
$
\hfill (3.42)
\\[2.25ex]
\indent
In view of  (3.40)  and (3.42)
from the functional system (3.36) we get the system
\\[2ex]
\mbox{}\hfill           
$
\displaystyle
\sum_{i =1}^{6} \,
y_i\varphi_i + 6x_1^2\varphi_7 +
(2x_2^3 + \alpha_2)\varphi_8 = U_0(x, y),
\qquad
\varphi_7 + x_2\varphi_8 = U_1(x, y).
$
\hfill (3.43)
\\[2.25ex]
\indent
The system (3.43) has the  solution 
\\[2ex]
\mbox{}\hfill                     
$
\varphi_i = z_i(x,y), \ \ i = 1,\ldots, 6,
\hfill
$
\\[2ex]
\mbox{}\hfill
$
\displaystyle
\varphi_7 = ({}-6x_1^2x_2 + 2x_2^3 + \alpha_2)^{{}-1}
\biggl({}-x_2U_0(x, y)  +
(2x_2^3 + \alpha_2)U_1(x, y) + x_2\,
\sum\limits_{i =1}^{6}\,y_iz_i(x, y)\biggr),
$
\hfill (3.44)
\\[1.75ex]
\mbox{}\hfill
$
\displaystyle
\varphi_8 \!=\! (\!{}-6x_1^2x_2 + 2x_2^3 + \alpha_2)^{{}-1}
\Bigl(U_0(x,
y) - 6x_1^2\,U_1(x, y) - \sum\limits_{i = 1}^{6}\, y_iz_i(x, y)\Bigr).
\hfill
$
\\[2.75ex]
\indent
Taking into account (3.31), (3.35), (3.40), (3.42), and (3.44)
the system (3.22)  for 
\linebreak
$|\alpha_5| + |\beta_5| \ne 0$ assumes the form
\\[2.25ex]
\mbox{}\hfill             
$
\partial_{x_i}w(x, y) - \mu(x, y)z_i(x, y) = 0, \ \
i =1,\ldots, 6,
\hfill
$
\\[2.25ex]
\mbox{}\hfill
$
\partial_{y_1}w(x, y) -
\mu(x, y)({}-6x_1^2x_2 + 2x_2^3 +
\alpha_2)^{{}-1}\Bigl({}-x_2U_0(x, y) \, +
\hfill
$
\\[1.5ex]
\mbox{}\hfill
$
\displaystyle
+\,
(2x_2^3 +\alpha_2)U_1(x, y) +
x_2\,\sum\limits_{i = 1}^{6}\,y_iz_i(x, y)\Bigr) = 0,
$
\hfill (3.45)
\\[2ex]
\mbox{}\hfill
$
\partial_{y_2}w(x, y) - \mu(x,y)({}-6x_1^2x_2 + 2x_2^3 +
\alpha_2)^{{}-1} \Bigl(U_0(x, y)  -
6x_1^2\,U_1(x,y) - \sum\limits_{i = 1}^{6}\,y_iz_i(x, y)\Bigr) =
0,
\hfill
$
\\[2.5ex]
\mbox{}\hfill
$
\displaystyle
\partial_{y_i}w(x, y) = 0, \ \
i = 3,\ldots,6,
\hfill
$
\\[2.5ex]
\indent
From the seventh and the eighth identities of system (3.45) we get 
\\[2ex]
\mbox{}\hfill
$
\partial_{y_1}w(x, y) +
x_2\,\partial_{y_2}w(x, y) = \mu(x, y)U_1(x,y)
\hfill
$
\\[2ex]
and  (by Theorem 3.5)  a function  $w$ is a 
partial integral of equation 
\\[2ex]
\mbox{}\hfill
$
\partial_{y_1} w +
x_2\;\!\partial_{y_2}\,w = 0.
\hfill
$
\\[1ex]
\indent
Hence,
\\[0ex]
\mbox{}\hfill          
$
w = y_2 - y_1 x_2+ h(x)
$
\hfill (3.46)
\\[2ex]
in view of representation (3.41).

From the first six identities and from the eighth identity of system (3.45) we get 
\\[2ex]
\mbox{}\hfill
$
\displaystyle
\sum\limits_{i = 1}^{6}\,
y_i\;\!\partial_{x_i}w(x, y) + ({}-6x_1^2x_2
+ 2x_2^3 + \alpha_2)\partial_{y_2}w(x, y)=
 \mu(x, y)\bigl(U_0(x, y)-
6x_1^2U_1(x,y)\bigr).
\hfill
$
\\[2.5ex]
\indent
Therefore, $w$ is a 
partial integral of equation 
\\[2ex]
\mbox{}\hfill              
$
\displaystyle
\sum\limits_{i =1}^{6}\,y_i\,\partial_{x_i}w +
({}-6x_1^2x_2 + 2x_2^3 + \alpha_2)\;\!\partial_{y_2}w = 0.
$
\hfill (3.47)
\\[2.5ex]
\indent
Since no function of the form (3.46) satisfies equation (3.47),
we conclude that system (PS)
doesn't have the autonomous  partial integrals 
other then constants for  $|\alpha_5| + |\beta_5| \ne 0.$

Let's consider the case $\alpha_5=\beta_5= 0.$

In this case  system (3.36) has the general solution
\\[2ex]
\mbox{}\hfill                 
$
\displaystyle
\varphi_i = z_i(x, y), \ \  i = 1,\ldots, 6,
\hfill
$
\\[2.5ex]
\mbox{}\hfill
$
\displaystyle
\varphi_7  = \biggl(
(2x_2^3 + \alpha_2)U_1(x,y) + x_2
\biggl(\, \sum\limits_{i = 1}^{6}\,y_iz_i(x, y) +
(y_3^2x_3^{-1} + \gamma_3x_3^3 \, +
\hfill
$
\\[2ex]
\mbox{}\hfill
$
+ \,
\delta_3x_3^{{}-1})({}-y_3 + \alpha_3x_3^2 + \beta_3)^{{}-1}
\bigl(U_{-1}(x, y) -
({}-y_5 + \gamma_5x_5)g(x, y)\bigr)
\, +
\hfill
$
\\[2ex]
\mbox{}\hfill
$
+ \,
\Bigl(\dfrac{1}{2}\, (3x_5 - 1)x_5^{{}-1}\,(x_5 -
1)^{{}-1}\,y_5^2 +
\delta_5x_5 (x_5 +1)(x_5 -1)^{{}-1}\Bigr)g(x, y) \, -
\hfill
$
\\[2ex]
\mbox{}\hfill
$
- \,
U_0(x, y)\biggr)\biggr)
({}-6x_1^2x_2 + 2x_2^3 + \alpha_2)^{{}-1},
\hfill
$
\\[-0.25ex]
\mbox{}\hfill (3.48)
\\[-0.25ex]
\mbox{}\hfill
$
\displaystyle
\varphi_8 =
{}-\biggl(\,
\sum\limits_{i = 1}^{6}\,y_iz_i(x, y) + 6x_1^2\,U_1(x,
y) - (y_3^2x_3^{{}-1} + \gamma_3x_3^3 +
\delta_3x_3^{{}-1})({}-y_3 \, +
\hfill
$
\\[2ex]
\mbox{}\hfill
$
+ \,
\alpha_3x_3^2 +\beta_3)^{{}-1}
\bigl(U_{-1}(x, y) -
({}-y_5 + \gamma_5x)g(x, y)\bigr) +
\Bigl(\dfrac{1}{2}\,(3x_5 - 1)x_5^{{}-1}y_5^2 \, +
\hfill
$
\\[2ex]
\mbox{}\hfill
$
+ \,
\delta_5x_5(x_5 + 1)(x_5 -1)^{{}-1}\Bigr)g(x, y) -
U_0(x, y)\biggr)
({}-6x_1^2x_2 + 2x_2^3 + \alpha_2)^{{}-1},
\hfill
$
\\[3ex]
\mbox{}\hfill
$
\varphi_9 = \bigl(U_{-1}(x, y) -
({}-y_5 + \gamma_5x_5)g(x, y)\bigr)
({}-y_3 + \alpha_3x_3^3 + \beta_3)^{{}-1},
\hfill
$
\\[2.5ex]
\mbox{}\hfill
$
\varphi_{11} = g(x, y).
\hfill
$
\\[3ex]
\indent
In view of  (3.31), (3.35), and (3.48) system  (3.22) for
$\alpha_5 = \beta_5 = 0$ takes the form
\\[2.25ex]
\mbox{}\hfill                  
$
\partial_{x_i}w(x, y) - \mu(x, y)z_i(x, y) = 0,
\ \ i = 1,\ldots, 6,
\hfill
$
\\[2ex]
\mbox{}\hfill
$
\partial_{y_1}w(x, y) - \mu(x,y)\varphi_7(x, y) = 0,
\hfill
$
\\[2ex]
\mbox{}\hfill
$
\partial_{y_2}w(x, y) -\mu(x, y)\varphi_8(x, y) = 0,
\hfill
$
\\[2ex]
\mbox{}\hfill
$
\partial_{y_3}w(x,y) -
\mu(x, y)({}-y_3 + \alpha_3x_3^2 + \beta_3)^{{}-1}
\bigl(U_{-1}(x, y)  -
({}-y_5 + \gamma_5x_5)g(x,y)\bigr) = 0,
$
\hfill (3.49)
\\[2ex]
\mbox{}\hfill
$
\partial_{y_4}w(x, y) = 0,
\hfill
$
\\[2ex]
\mbox{}\hfill
$
\partial_{y_5}w(x, y) - \mu(x, y)g(x, y) = 0,
\hfill
$
\\[2ex]
\mbox{}\hfill
$
\partial_{y_6}w(x, y) = 0.
\hfill
$
\\[2.5ex]
\indent
From the seventh and the eighth identities of system (3.49)
we get 
\\[2ex]
\mbox{}\hfill
$
\partial_{y_1}w(x, y) + x_2\partial_{y_2}w(x,
y) = \mu(x, y)U_1(x, y)
\hfill
$
\\[2ex]
and much as we do in the first case
we conclude that $w$ is a 
partial integral of the equation  $\partial_{y_1}w + x_2\;\!\partial_{y_2}w =0.$

From this in view of  (3.44) we obtain that
\\[2ex]
\mbox{}\hfill                 
$
w = y_2 - y_1x_2 + h(x, y_3, y_5).
$
\hfill (3.50)
\\[2ex]
\indent
It follows from the first six, eighth, and eleventh identities of system (3.49) that 
\\[2ex]
\mbox{}\hfill
$
\displaystyle
\sum\limits_{i =1}^{6}\,y_i\partial_{x_i}w(x, y) +
({}-6x_1^2x_2 +\alpha_2)\partial_{y_2}w(x, y) \,+
\hfill
$
\\[2ex]
\mbox{}\hfill
$
+\,\biggl(\Bigl(
\dfrac{1}{2}\  (3x_5 -1)x_5^{{}-1}\,(x_5
- 1)^{{}-1}\,y_5^2 +
\delta_5x_5(x_5 +1)(x_5 - 1)^{{}-1}\Bigr) \,-
\hfill
$
\\[2ex]
\mbox{}\hfill
$
-\,
({}-y_5+ \gamma_5x_5) (y_3^2x_3^{{}-1} + \gamma_3x_3^3 +
\delta_3x_3^{{}-1}) ({}-y_3 +
\alpha_3x_3^2 +
\beta_3)^{{}-1} \biggr) \partial_{y_5} w(x,y) =
\hfill
$
\\[2ex]
\mbox{}\hfill
$
=
 \mu(x, y)
\bigl({}-(y_3^2x_3^{{}-1} + \gamma_3x_3^3 +
\delta_3x_3^{{}-1})({}-y_3 +
\alpha_3x_3^2 +\beta_3)^{{}-1}\,U_{-1}(x, y) +
U_0(x, y) - 6x_1^2\,U_1(x, y)\bigr).
\hfill
$
\\[3ex]
\indent
Therefore (by Theorem 3.5)
a
scalar function  $w$ is a partial integral of the equation
\\[2ex]
\mbox{}\hfill              
$
\displaystyle
\sum\limits_{i= 1}^{6}\,
y_i\;\!\partial_{x_i}w + ({}-6x_1^2x_2 +
2x_2^3 + \alpha_2)\;\!\partial_{y_2}w \,
+
\hfill
$
\\[2ex]
\mbox{}\hfill
$
+ 
\biggl(\Bigl(
\dfrac{1}{2}\ (3x_5 -1)x_5^{{}-1}\,
(x_5 - 1)^{{}-1}\,y_5^2 + \delta_5x_5(x_5 +1)
(x_5 - 1)^{{}-1}\Bigr)  -
$
\hfill (3.51)
\\[2ex]
\mbox{}\hfill
$
- \,
({}-y_5 + \gamma_5x_5) (y_3^2x_3^{{}-1} + \gamma_3x_3^3 +
\delta_3^{{}-1}) ({}-y_3 + \alpha_3x_3^2 + \beta_3)^{{}-1}\biggr)
 \partial_{y_5}w = 0.
\hfill
$
\\[2.75ex]
\indent
Since no function of the form (3.50) satisfies equation (3.51),
the  Painlev\'{e} system (PS)
doesn't have the autonomous  partial integrals 
other then constants for  $\alpha_5=\beta_5= 0.$ \k
\vspace{0.5ex}

Thus, we have given answers to the questions concerning the  
 relations between the general solutions to the 
irreducible Painlev\'{e} equations ($\!\text{P-\;\!}i$\!), $i = 1,\ldots, 6.$
These answers are as follows.
\vspace{0.25ex}

{\bf 1.}
There is no functional relation 
$w(x_1, \ldots , x_6) = 0$
with  holomorphic function $w$ 
between the general solutions $x_i,\, i = 1,\ldots, 6,$ 
\vspace{0.5ex}
to the  Painlev\'{e} equations {\rm($\!\text{P-\;\!}i$\!)}, $i = 1,\ldots, 6.$

{\bf 2.}
There is no functional relation of the form (3.26)
with  holomorphic function $w$ 
between the general solutions 
$x_i,\, i = 1,\ldots, 6,$ 
to the  Painlev\'{e} equations {\rm($\!\text{P-\;\!}i$\!)}, $i = 1,\ldots, 6,$
and  their derivatives ${\sf D}x_i.$
\vspace{0.5ex}

At the same time, the differential Painlev\'{e} system (PS)
has nonautonomous  partial integrals and therefore we can 
assert  as follows.
\vspace{0.25ex}

{\bf 3.}
There exists a  functional relation of the form 
$w(t, x_1, \ldots , x_6,\;\! {\sf D}x_1, \ldots ,
{\sf D}x_6) = 0$
between the general solutions 
$x_i,\, i = 1,\ldots, 6,$ 
to the  Painlev\'{e} equations {\rm($\!\text{P-\;\!}i$\!)}, $i = 1,\ldots, 6,$
and  their first derivatives ${\sf D}x_i.$
\vspace{0.5ex}

In particular, if the general  solution $x_k,\, k \in \{1, \ldots, 6\},$ 
to the $k\!$-th Painlev\'{e} equation {\rm($\!\text{P-\;\!}k$\!)}
and  its  derivative ${\sf D}x_k$ are known, then 
we can assert  as follows.
\vspace{0.25ex}

{\bf 4.}
There exists a first-order ordinary differential equation 
of the form 
\\[2ex]
\mbox{}\hfill
$
w_s(t, z, {\sf D}z, x_k, {\sf D}x_k) = 0,
\hfill
$
\\[2ex]
 where  $t$ is the independent variable,
$z$ is an  unknown function, 
$x_k$ is the general  solution to the Painlev\'{e} equation {\rm($\!\text{P-\;\!}k$\!)},
whose general  solution
is the general  solution to the $s\!$-th Painlev\'{e} equation {\rm($\!\text{P-\;\!}s$\!)},\,$s \in \{1, \ldots , 6\},\, s \ne k,$ 
that is, $z = x_s,\, s \ne k.$

}
\end{document}